\documentclass[11pt,reqno]{amsart}

\usepackage[a4paper,left=23mm,right=23mm,top=23mm,bottom=17mm]{geometry}

\usepackage{amsmath,amssymb,mathtools,bbm}
\usepackage{enumitem}
\usepackage{comment}
\usepackage{cite}
\usepackage{subcaption}
\usepackage[normalem]{ulem}
\usepackage{microtype}
\usepackage[dvipsnames,table]{xcolor}
\usepackage{graphicx}
\usepackage{placeins}
\usepackage{fix-cm}
\usepackage{tikz}

\newtheorem{theorem}{Theorem}[section]
\newtheorem{lemma}[theorem]{Lemma}
\newtheorem{corollary}[theorem]{Corollary}

\newtheorem{proposition}[theorem]{Proposition}

\theoremstyle{definition}
\newtheorem{definition}[theorem]{Definition}
\newtheorem{remark}[theorem]{Remark}

\numberwithin{equation}{section}

\usepackage{etoolbox}
\makeatletter
\patchcmd{\@settitle}
  {\uppercasenonmath\@title}
  {}
  {}{}
\patchcmd{\@settitle}
  {\bfseries}
  {\fontsize{15.5}{16}\selectfont}
  {}{}
\def\@setauthors{
  \begingroup
  \trivlist
  \centering
  \footnotesize
  \@topsep30\p@\relax
  \advance\@topsep by -\baselineskip
  \item\relax
  \fontsize{12}{14}\selectfont
  \authors
  \endtrivlist
  \endgroup
}
\def\section{\@startsection{section}{1}%
  \z@{.7\linespacing\@plus\linespacing}{.5\linespacing}%
  {\normalfont\bfseries\fontsize{12}{14}\selectfont\centering}}
\makeatother

\usepackage{stmaryrd}
\SetSymbolFont{stmry}{bold}{U}{stmry}{m}{n}

\allowdisplaybreaks[4]
\setlist{topsep=4pt,itemsep=2pt}

\let\AMSnormal\normalsize
\renewcommand{\normalsize}{%
  \AMSnormal
  \baselineskip=\dimexpr\baselineskip\relax
    plus .25pt minus .05pt\relax
  \parskip=0pt plus .3pt\relax
  \jot=.25\baselineskip
  \abovedisplayskip=.35\baselineskip plus 3pt minus 2pt\relax
  \belowdisplayskip=.35\baselineskip plus 3pt minus 2pt\relax
  \abovedisplayshortskip=0pt plus 3pt\relax
  \belowdisplayshortskip=.175\baselineskip plus 2pt minus 1pt\relax
}

\normalsize
\usepackage{hyperref}
\hypersetup{colorlinks=true,linkcolor=blue!90!black,citecolor=red!95!black,urlcolor=blue!65!black,pdfstartview=Fit,bookmarksnumbered=true,breaklinks=true,hypertexnames=false}

\begin{document}
\title[Anisotropic propagation for reaction-diffusion equations]
{The geometry of anisotropic propagation for multidimensional reaction-diffusion equations in shifting habitats: A zero-sum stopping game approach}

\author[Tao, Li and Ruan]{Wen Tao$^{1}$, Wan-Tong Li$^{1,*}$ and Shigui Ruan$^{2}$}
\thanks{\hspace{-0.53cm}
$^1$School of Mathematics and Statistics, Lanzhou University, Lanzhou, Gansu 730000, P. R. China.\\
$^2$Department of Mathematics, University of Miami, Coral Gables, FL 33146, USA.\\
$^*${\sf Corresponding author} (wtli@lzu.edu.cn)}

\date{\today}

\begin{abstract}
We study the propagation geometry of multidimensional reaction-diffusion equations in shifting habitats.
Using geometric optics techniques and viscosity solution theory, we prove that the associated rate functions converge locally uniformly to the unique viscosity solution of a variational inequality with singular initial data.
We identify this limit as the value of a zero-sum stopping game whose zero set determines the Wulff shape.
This zero set coincides with that of the positive part of the value function of an unconstrained optimal control problem if and only if Freidlin's condition holds. When Freidlin's condition fails, we derive a pointwise characterization of the stopping-game zero set, yielding a convex hull representation. These results then yield explicit Wulff shapes and directional spreading speeds for all habitat speed. As the habitat speed increases, the Wulff shape evolves from a large ball to a truncated ball. It then develops conical, cylindrical, and paraboloidal boundary components attached to spherical parts, before eventually becoming a smaller ball. Some of these geometric transitions occur while both axial spreading speeds remain unchanged. Moreover, by characterizing the optimal trajectories, we distinguish two mechanisms of nonlocal speed selection and show that such selection can occur in oblique directions even when it is absent in both axial directions.



\end{abstract}

\maketitle

\section{Introduction}
\label{section:introduction}

Changes in climate can shift the spatial distribution of suitable habitats, and reaction-diffusion models provide a natural framework for studying how populations respond through growth and dispersal. Potapov and Lewis \cite{potapov2004} studied persistence and invasion in habitats with moving boundaries, while Berestycki et al. \cite{berestycki2009} established conditions under which a species can persist in a shifting favorable patch.
Two central concepts in this theory are asymptotic spreading speeds and forced waves.
Spreading speeds describe the large time expansion of a population from an initial distribution, as in the work of Aronson and Weinberger \cite{aronson1975} and Weinberger \cite{wein1982}.
Forced waves are population profiles moving at the speed of the habitat.
Their existence can be studied using generalized principal eigenvalues \cite{berestycki2018}, building on the theory of Berestycki et al. \cite{berestycki1994}.

We consider a moving transition between two regions that both support population growth.
A population can persist in the less favorable region without following the transition, while access to the more favorable region may alter its spreading speeds in different directions.
We model this transition by a monotone growth profile that varies along a fixed spatial direction and shifts at a prescribed speed.
More precisely, we consider the multidimensional Fisher--KPP Cauchy problem:
\begin{align}
\label{eq:rd-nd-main}
  \begin{cases}
    u_t=\Delta u+u\bigl(r(\mathbf{x}\cdot\mathbf{e}-ct)-u\bigr) & \text{in }\mathbb{R}^{N}\times\mathbb{R}^+,\\
    u(\mathbf{x},0)=g(\mathbf{x}) & \text{on }\mathbb{R}^{N},
  \end{cases}
\end{align}
where $u(\mathbf{x},t)$ denotes the population density at location $\mathbf{x}$ and time $t\geqslant 0$, $N\geqslant2$ is an integer, $\mathbf{e}\in\mathbb S^{N-1}$ is a fixed unit vector specifying the direction of environmental shifting, and $c\in\mathbb{R}$ is the signed shifting speed.
The growth rate $r\in C(\mathbb{R};\mathbb{R}^{+})$ is monotone with positive limits at $\pm\infty$. 
Without loss of generality, we assume that $r$ is nondecreasing, with
\begin{align*}
  \lim_{s\to-\infty}r(s)=r_{1}, \qquad \lim_{s\to+\infty}r(s)=r_{2}, \qquad 0<r_{1}<r_{2}.
\end{align*}
For a nonincreasing profile, set $\tilde r(s)=r(-s)$, $\tilde{\mathbf{e}}=-\mathbf{e}$, and $\tilde c=-c$.
Then $\tilde r$ is nondecreasing and
\begin{align*}
  r(\mathbf{x}\cdot\mathbf{e}-ct)
  =\tilde r(\mathbf{x}\cdot\tilde{\mathbf{e}}-\tilde c\,t),
\end{align*}
so the problem reduces to the nondecreasing case.
The initial datum $g$ is continuous, nonnegative, nontrivial, and compactly supported.

We study the large time expansion of the population governed by \eqref{eq:rd-nd-main}.
In a homogeneous environment, the movement of the population is isotropic so that spreading occurs at the same speed in all directions, and the invaded region expands asymptotically like a ball.
Here, the space--time heterogeneity of the growth rate may lead to different spreading speeds in different directions.
The natural object is therefore the asymptotic expansion set, which describes the large time shape of the invaded region.
Our aim is to determine this set, referred to here as the \textit{Wulff shape} (representing the anisotropic analog of a ball), and the corresponding directional spreading speeds for every $c\in\mathbb{R}$.
In particular, we ask how the shifting speed and the difference between the two limiting growth rates affect the geometry of invasion.

\subsection{Previous work}

Most explicit results for propagation in shifting habitats are obtained in one spatial dimension. Li et al. \cite{li2014} and Hu et al.
\cite{hu2019} used the method of super- and subsolutions to study favorable--unfavorable habitats (i.e., $r_{1}<0$ and $r_{2}>0$), showing that a critical habitat speed separates persistence from extinction for the populations.
The abstract theories of monotone semiflows in \cite{liang2007,fang2017} describe spreading in translation-invariant and space--time periodic settings.
Yi and Zhao \cite{yi2020b} developed a framework without spatial translation invariance that extends spreading results for shifting habitats.
For forced waves, Hu and Zou \cite{hu2017}, Fang et al. \cite{fang2016}, and Berestycki and Fang \cite{berestycki2018} established existence, multiplicity, and attractivity results under different assumptions on the habitat.

For environments with two favorable limits, Hu et al.
\cite{hu2020} studied spreading speeds in several ranges of the habitat speed.
Subsequently, Lam and Yu \cite{lam2022a} developed a Hamilton-Jacobi framework to study spreading speeds in weakly monotone habitats with exponentially bounded initial data. They also completed the classification for compactly supported initial data, including nonlocal pulling and locked fronts.
Lam and Lou \cite[Section 6.4]{lamlou2022} also derive the nonlocally pulled spreading speed using the method of super- and subsolutions.
Tao et al. \cite{tao2026} further identified a new nonlocal pulling mechanism in decreasing favorable habitats for suitable shifting speeds and exponential initial tails with small decay rates.
Beyond weakly monotone habitats, Giletti et al.
\cite{giletti2024} combined principal eigenvalue analysis with the method of super- and subsolutions to study pulling and locking caused by a moving region of increased growth.
Lam et al. \cite{lam2025} further developed a flux-limited Hamilton-Jacobi framework, using the theory in \cite{imbert2017}, to classify spreading speeds for compactly supported initial data and positive growth profiles without assuming monotonicity. We refer to Lam et al. \cite{lam2027} for a survey of recent developments on this topic.
These results explain how environmental shifting selects a spreading speed in one dimension, including speeds that depend on growth away from the front and speeds equal to that of the environment.

In higher dimensions, Berestycki and Rossi \cite{berestyck2008} studied shifting habitats with a bounded favorable region and an unfavorable far field.
They used generalized principal eigenvalues to characterize persistence and extinction, and proved convergence to a unique forced wave in the persistence regime.
Guo et al. \cite{guo2026} used asymptotic spectral radius theory to obtain criteria for extinction and the existence of forced waves in cooperative reaction-diffusion systems on cylindrical domains with an unfavorable far field.
Yi et al. \cite{yi2025} considered a Nicholson equation with local diffusion and maturation delay in a directionally shifting favorable--unfavorable habitat.
Using iterations of traveling-wave maps, they showed that, in the spreading regimes, the invaded region is asymptotically an expanding ball truncated by the moving favorable half-space.
Vo \cite{vo2026} established truncated-ball spreading results for a time-periodic parabolic--elliptic Keller--Segel system in a shifting favorable--unfavorable habitat under suitable damping conditions.
These studies deal with habitats with an unfavorable region, whereas both limiting environments in our setting support population growth.

Propagation sets for scalar reaction-diffusion equations with general space--time heterogeneity were studied by Berestycki et al. \cite{berestycki2008}.
Berestycki and Nadin obtained spectral bounds on one-dimensional spreading speeds \cite{beres2012} and developed a theory of expansion sets using generalized principal eigenvalues, Hamilton-Jacobi equations, and variational representations \cite{berestycki2022}.
Their Proposition 3.14 describes a stationary two-dimensional transition between two positive limiting growth rates and gives a convex hull formula for the expansion set.
It provides a geometric counterpart to the one-dimensional theory, but does not describe how motion of the transition changes the front.

\subsection{Research gap and difficulties}

For the monotone environmental profile considered here, the corresponding one-dimensional spreading speeds are already well understood.
The question is to determine the expansion in directions not parallel to the environmental axis and the resulting geometry of the front.
The two axial speeds, corresponding to propagation in the directions $\mathbf e$ and $-\mathbf e$, alone cannot specify this geometry.
In particular, it is not clear from the one-dimensional classification whether a change in habitat speed can alter the shape of the front without changing either axial speed, or whether nonlocal pulling can occur away from the axis when neither axial front exhibits it.
The existing multidimensional results above do not answer these questions for a moving transition between two favorable regions.

The one-dimensional dependence of the coefficient does not resolve this difficulty.
The initial population is localized in all spatial variables, and restricting the solution to a ray does not give a closed one-dimensional equation.
Super- and subsolutions adapted to individual directions must provide compatible sharp bounds for the full expansion set, while forced waves do not by themselves determine this set for localized initial data.
In one-dimensional Hamilton-Jacobi calculations, self-similarity yields a single-variable problem whose solution can be identified by matching admissible branches.
In higher dimensions, the rotational symmetry of the limiting problem still leaves a self-similar equation in two spatial variables.
One must determine where admissible branches apply and locate the free boundary separating the invaded and uninvaded regions.

A variational representation offers another way to determine this boundary, but it introduces a question of its own.
An unconstrained optimal control problem minimizes the total action, which balances the cost of dispersal against the growth along a trajectory.
The obstacle in the limiting Hamilton-Jacobi equation, however, leads to constraints at intermediate times.
Nonpositivity of the total action does not ensure nonpositivity of all intermediate actions.
Even though the unconstrained action separates into an axial contribution and a quadratic cost in the orthogonal directions, this separation does not by itself solve the obstacle problem.
We must determine when the unconstrained action gives the expansion set and how to compute the true zero set when this simplification is not available.

\subsection{Variational approach and main ideas}

The geometric optics approach originates in Freidlin's large deviation analysis \cite{freidlin1985} and the PDE formulation of Evans and Souganidis \cite{evans1989}, with further developments for cooperative systems in \cite{barles1990}.
Within this framework, Liu et al. \cite{liu2020} analyzed a one-dimensional propagation problem with a discontinuous Hamiltonian.
In a road--field model, Henderson and Lam \cite{henderson2026} used Hamilton-Jacobi equations and optimal control to determine the expansion set and relate its geometry to minimizing trajectories.
For the present problem, we use a stopping-game formulation \cite{evans1989} of the limiting Hamilton-Jacobi obstacle problem and derive reductions that preserve the propagation set.

Under hyperbolic rescaling $u^{\varepsilon}(\mathbf{x},t)=u(\frac{\mathbf{x}}{\varepsilon},\frac{t}{\varepsilon})$, the Hopf-Cole rate functions $w^{\varepsilon}=-\varepsilon\ln u^{\varepsilon}$ converge locally uniformly in $\mathbb{R}^{N}\times\mathbb{R}^{+}$ to a nonnegative continuous viscosity solution $w$ of a variational inequality with singular initial data:
\begin{equation*}
  \left\{\begin{array}{ll}
    \min\left\{\partial_t w+|Dw|^{2}+R(\mathbf{x},t),\,w\right\} =0 \quad \text{in }\mathbb{R}^{N}\times\mathbb{R}^{+}, \\
    w(\mathbf{x},0)=
    \begin{cases}
      0, & \text{if }\mathbf{x}=\mathbf{0},\\
      +\infty, &\text{if }\mathbf{x}\neq \mathbf{0},
    \end{cases}
  \end{array}\right.
\end{equation*}
where $R(\mathbf{x},t)= r_{1}\mathbbm{1}_{\{\mathbf{x}\cdot\mathbf{e}<ct\}} +r_{2}\mathbbm{1}_{\{\mathbf{x}\cdot \mathbf{e}\geqslant ct\}}$.
The singular initial condition is understood in the relaxed boundary sense described in Remark \ref{remark:relax-initial-condition}.

We identify the rate function $w$ with the value of a zero-sum stopping game:
\begin{align*}
  W(\mathbf{x},t) =\sup_{\vartheta\in\Theta} \inf_{\gamma\in\mathbb{X}^{\mathbf{x},t}} \int_0^{t\wedge\vartheta[\gamma]} \left(\frac{|\dot\gamma(s)|^{2}}{4} -R(\gamma(s),t-s) \right)\,ds.
\end{align*}
Here $\mathbb{X}^{\mathbf{x},t}$ consists of trajectories $\gamma\in H^1_{\mathrm{loc}}([0,\infty);\mathbb{R}^{N})$ with $\gamma(0)=\mathbf{x}$ and $\gamma(t)=\mathbf{0}$.
The minimizing player controls $\gamma$, while the maximizing player chooses a stopping time $\vartheta\in\Theta$ based only
on the portion of the trajectory already observed.
Precise definitions are given in Section \ref{section:main-results}.

We prove the exact pathwise formula
\begin{align*}
  W(\mathbf{x},t) =\min_{\gamma\in\mathbb{X}^{\mathbf{x},t}} \max_{\sigma\in[0,t]} \int_0^\sigma \left(\frac{|\dot\gamma(s)|^{2}}{4} -R(\gamma(s),t-s) \right)\,ds.
\end{align*}
This identity replaces optimization over stopping times by minimization over trajectories of the maximum action over all intermediate times.
Freidlin and Lee \cite[Lemma~2.4]{freidlin1996} proved an analogous equivalence for spatially heterogeneous systems with Lipschitz coefficients.
For a single monotone transition between two favorable regions, we establish the corresponding case of the min-max representation conjectured by Lam et al. \cite[Section 2.5]{lam2025} and its multidimensional analogue; see Remark \ref{remark:minmax-lam-conjecture}.
In particular, $W(\mathbf{x},t)=0$ if and only if some admissible trajectory has nonpositive action on every interval $[0,\sigma]$, where $0\leqslant\sigma\leqslant t$.
Self-similarity gives $W(\mathbf{x},t)=tW(\frac{\mathbf{x}}{t},1)$.
We prove that the zero set of $W$ determines the expansion set of $u$.
The Wulff shape is therefore given by
\begin{align*}
  \mathcal{W} =\{\mathbf{x}\in\mathbb{R}^{N}\mid W(\mathbf{x},1)=0\}.
\end{align*}

We then develop two complementary reductions to derive the expansion set.
The first reduction connects the stopping game with unconstrained optimal control.
Let $J(\mathbf{x},t)$ denote the infimum of the total action over trajectories joining $(\mathbf{x},t)$ to $(\mathbf{0},0)$, with no stopping before time $t$ (see \eqref{eq:optimal-control-J}).

For the present model with $N\geqslant2$, we prove that Freidlin's condition (Definition~\ref{def:backward-freidlin-condition}) is necessary and sufficient for the identity
\begin{align*}
  W(\mathbf{x},t)=\max\{0,J(\mathbf{x},t)\} \quad \text{for all }(\mathbf{x},t)\in\mathbb{R}^{N}\times\mathbb{R}^{+}.
\end{align*}
This condition holds precisely when $c\leqslant-2\sqrt{r_{2}}$ or $c\geqslant2\sqrt{r_{2}-r_{1}}$.
In these regimes, the minimizing trajectories of $J$ determine the expansion set explicitly.

In the complementary non-Freidlin regime $-2\sqrt{r_{2}}<c<2\sqrt{r_{2}-r_{1}}$, the zero set of $\max\{0,J(\cdot,1)\}$ strictly contains the Wulff shape $\mathcal W$, so unconstrained optimal control overestimates the expansion set; see Theorem \ref{thm:main-freidlin-dichotomy}(iii) and Remark \ref{remark:scope-unconstrained-reduction}.
We therefore retain the intermediate-time constraints and work directly with the zero set of $W$.
We prove that every point in this zero set admits a trajectory with nonpositive running cost almost everywhere.
This reduces the determination of the expansion set to a geometric problem with pointwise speed bounds, yielding an exact convex-hull formula even when Freidlin's condition fails.
These two reductions determine the expansion set without an explicit formula for the full multidimensional viscosity solution.

\begin{figure}[htpb]
  \centering
  \includegraphics[width=0.7\textwidth]{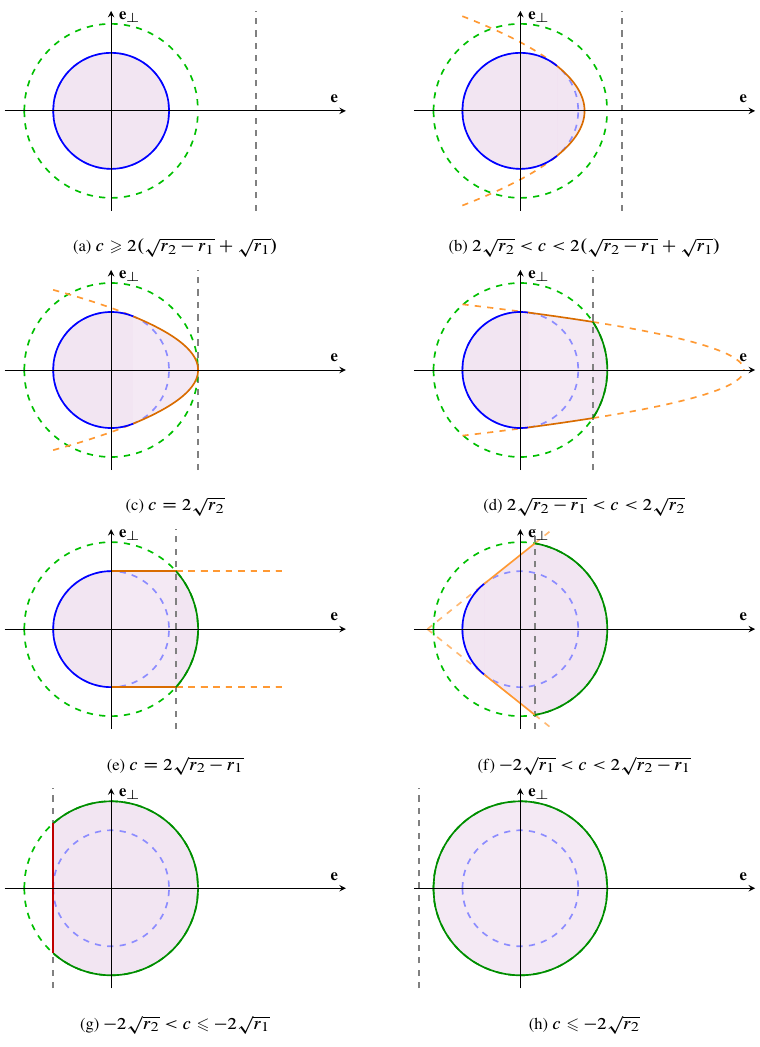}
  \caption{Propagation geometry in an orthonormal basis $(\mathbf{e},\mathbf{e}_{\perp})$ of $\mathbb{R}^{2}$, with regimes ordered by decreasing $c$.
   The lavender region represents the Wulff shape.
   Solid curves form its boundary, while dashed curves show the corresponding non-front branches of the low-growth circle, the high-growth circle, and the intermediate branches when present.
   The low-growth and high-growth circles have radius $2\sqrt{r_{1}}$ and $2\sqrt{r_{2}}$, respectively.
   Dashed gray lines indicate the shifting interface, and the red segment is the part of the front that coincides with the shifting interface.
   In higher dimensions, the Wulff shape is rotationally symmetric about the $\mathbf{e}$-axis, and these plots show its sections in planes containing this axis.}
  \label{fig:front-regimes-parameterized}
\end{figure}

\subsection{Main contributions}

Our axial spreading speeds recover the one-dimensional formulas for weakly monotone habitats and compactly supported initial data \cite{lam2022a}, see Subsection \ref{subsec:discussion-1d}.
We obtain the full expansion set for every $c\in\mathbb{R}$ by determining how the obstacle constraints can be reduced without changing the propagation set.
The geometric classification and the optimal trajectories together explain how environmental shifting changes the front and the mechanisms of directional speed selection even when both axial speeds remain fixed.
Our main contributions are as follows:

 \textbf{(a) Singular initial data and uniqueness.}
  We derive precise growth estimates for the solution of the limiting variational inequality with singular initial data.
  Heat-semigroup estimates and Laplace-type bounds yield
  \begin{align*}
  \frac{|\mathbf{x}|^{2}}{4t}-Ct \leqslant w(\mathbf{x},t) \leqslant\frac{|\mathbf{x}|^{2}}{4t}+Ct \quad \text{for some } C>0.
  \end{align*}
  These bounds specify a growth class in which we establish the uniqueness of the continuous viscosity solution; see Theorem \ref{thm:wkb-direct-convergence}.
  Earlier treatments of singular initial data use dynamic programming \cite{evans1989}, the solution-semigroup approach \cite{crandall1989}, or, in one dimension, a reduction to a boundary value problem based
  on the $1$-homogeneity property \cite{lam2022a,lam2025}.
  Our argument establishes uniqueness in the full space--time domain by a direct PDE comparison argument.
  Self-similarity then follows from uniqueness and the scaling invariance of the equation.

  \textbf{(b) Min-max representation and convex-hull geometry.}
  For every $c\in\mathbb{R}$, we identify the geometric optics limit with a stopping-game value and prove an exact min-max formula that characterizes its zero set by nonpositive action at every intermediate time; see Proposition \ref{prop:inf-action-characterization}.
  We also prove that Freidlin's condition is necessary and sufficient for recovering the expansion set from $J$; see Theorem \ref{thm:main-freidlin-dichotomy}.
  The main step in the non-Freidlin regime is to prove that whenever $W(\mathbf{x},t)=0$, there exists an admissible trajectory with nonpositive running cost almost everywhere.
  This condition is stronger for a given trajectory, but imposing it does not change the set of endpoints; see Proposition \ref{prop:safe-action-pointwise-char}.
  For a fixed coordinate along $\mathbf e$ and a fixed first contact time with the shifting interface, we obtain a sharp bound on the distance from the $\mathbf e$-axis; see Lemma \ref{lem:fixed-contact-optimization}.
  We use trajectories attaining this bound to construct trajectories with the original endpoints and nonpositive running cost.
  The resulting pointwise speed bounds, together with the geometry of the shifting interface, give the convex hull formula for the Wulff shape in Theorem \ref{thm:main-freidlin-dichotomy}:
  \begin{align*}
  \mathcal{W}=\mathrm{co}(D_{1}\cup D_{2}),
  \end{align*}
  where $D_{1}=\overline{B_{2\sqrt{r_{1}}}(\mathbf{0})}$, $D_{2}=\{\mathbf{x}\in\overline{B_{2\sqrt{r_{2}}}(\mathbf{0})}\mid \mathbf{x}\cdot\mathbf{e}\geqslant c\}$.
  At $c=0$ and $N=2$, our formula recovers the stationary geometry of Berestycki and Nadin \cite[Proposition 3.14]{berestycki2022}.
  Our result determines how this geometry changes under environmental shifting throughout the non-Freidlin regime and in every dimension $N\geqslant2$.

  \textbf{(c) Geometric transitions beyond the axial speeds.}
  Writing $c_\mathrm{dir}(\mathbf n)$ for the radial speed in direction $\mathbf n$, we have, throughout $-2\sqrt{r_{1}}< c\leqslant2\sqrt{r_{2}}$, that
  \begin{align*}
  c_\mathrm{dir}(\mathbf{e})=2\sqrt{r_{2}}, \qquad c_\mathrm{dir}(-\mathbf{e})=2\sqrt{r_{1}}.
  \end{align*}
  These fixed axial speeds coexist with a conical intermediate front for $-2\sqrt{r_{1}}< c<2\sqrt{r_{2}-r_{1}}$, a cylindrical one at $c=2\sqrt{r_{2}-r_{1}}$, and a paraboloidal one for $2\sqrt{r_{2}-r_{1}}<c\leqslant 2\sqrt{r_{2}}$.
  The cylindrical transition coincides with the transition from the non-Freidlin to the Freidlin regime as $c$ increases.
  For example, the speed in any direction orthogonal to $\mathbf{e}$ is $2\sqrt{r_{2}}$ at $c=0$ but $2\sqrt{r_{1}}$ at the cylindrical transition.
  For $-2\sqrt{r_{2}}<c\leqslant-2\sqrt{r_{1}}$, the expansion set is a ball truncated by the moving half-space, as in the geometry found by Yi et al.
  \cite{yi2025}, but here both limiting habitats remain favorable.
  The complete classification also includes spherical and planar front pieces; see Figure \ref{fig:front-regimes-parameterized}.
  For fixed $N$ and $\mathbf e$, the expansion set depends only on $r_{1}$, $r_{2}$, and $c$, not on the detailed monotone growth profile or the choice of nontrivial compactly supported initial data.

 \textbf{(d) Nonlocal speed selection and distinct trajectory mechanisms.}
  The higher-growth habitat can enhance directional spreading speeds in the lower-growth region beyond the homogeneous value $2\sqrt{r_{1}}$.
  We give a variational interpretation of this nonlocal speed selection through two types of optimal trajectories, distinguished by whether they leave the evolving front.

  On the paraboloidal part of the front, away from its junctions, the optimal trajectories remain in $\{W>0\}$ at every intermediate time, outside the invaded region; see Theorem \ref{thm:main-front-trajectories}.
  These trajectories describe how exponentially small densities influenced by the higher-growth habitat determine directional spreading speeds.
  This provides a multidimensional counterpart of the optimal trajectory interpretation in \cite[Definition~3.6(ii)]{lam2027}.
  In particular, for $2\sqrt{r_{2}-r_{1}}<c<2\sqrt{r_{2}}$, this mechanism occurs in oblique directions while neither axial front exhibits it.
  In the corresponding one-dimensional problem with compactly supported initial data, nonlocal speed selection requires $c>2\sqrt{r_{2}}$ \cite{lam2022a}.
  This nonlocal speed selection mechanism in oblique directions complements those studied in \cite{holzer2014,girardin2019,liu2021a,liu2020}.
  It occurs for compactly supported initial data, unlike the additional mechanism in \cite{tao2026} for exponentially decaying initial data in decreasing habitats.

  On the conical part, away from its junctions, the optimal trajectories change direction while remaining on the propagation front at every intermediate time.
  Their running cost is zero almost everywhere; see Theorem \ref{thm:main-front-trajectories} and Remark \ref{remark:cone-front-nonlocal-selection}.
  This extends the interpretation of the stationary two-dimensional example in \cite[Section 5]{lam2027} to shifting interfaces and higher dimensions.

  At the critical speed $c=2\sqrt{r_{2}-r_{1}}$, Freidlin's condition holds.
  As on the conical part, the optimal trajectories reaching the cylindrical front have zero running cost almost everywhere and remain on the propagation front at every intermediate time; see Theorem~\ref{thm:main-front-trajectories}.

\begin{figure}[htpb]
  \centering
  \includegraphics[width=0.75\textwidth]{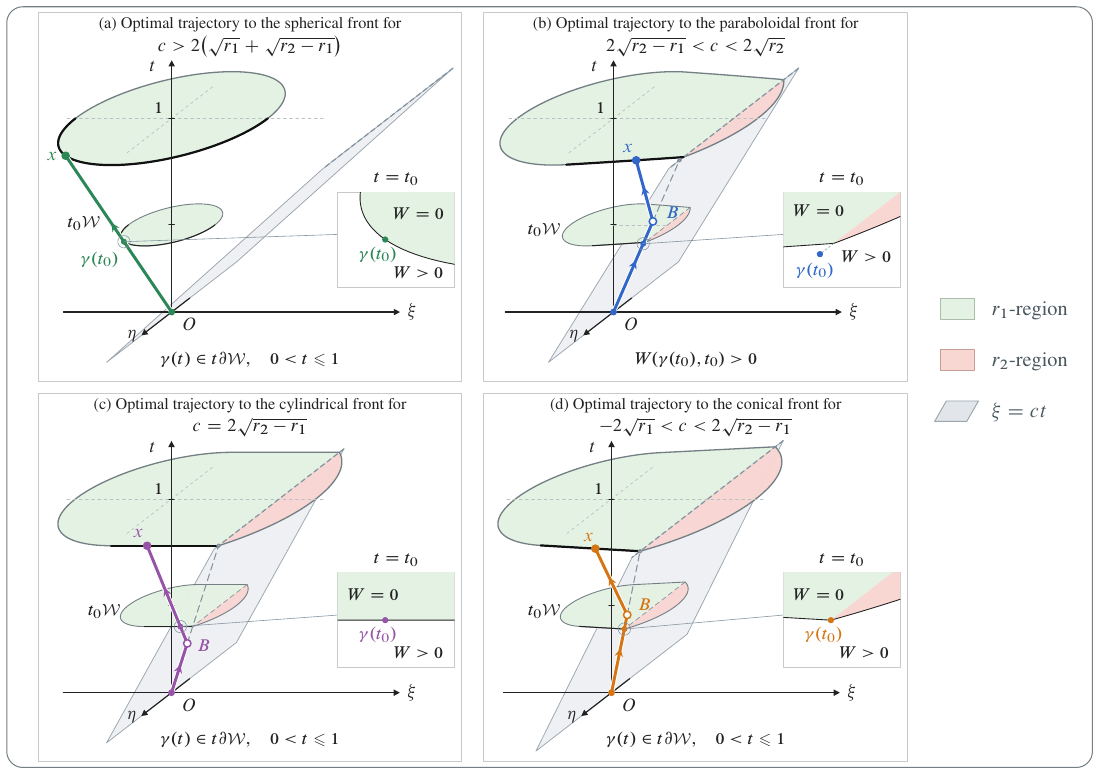}
  \caption{Optimal trajectories from the origin at time $0$ to $\mathbf{x}\in\partial\mathcal{W}$ at time $1$ in $\mathbb{R}^{2}$.
 The spatial coordinates are defined by $\mathbf{x}=\xi\mathbf{e}+\eta\mathbf{e}_{\perp}$ in the orthonormal basis $(\mathbf{e},\mathbf{e}_{\perp})$.
 Here $\gamma$ is parametrized forward in time; the corresponding backward path in the stopping-game formulation is $s\mapsto\gamma(1-s)$.
 Colored curves represent trajectories, black curves highlight the corresponding front pieces, and $B$ marks where a trajectory leaves the shifting interface.
 All panels show the section $t_0\mathcal{W}$ at the same intermediate time.
 The marked point $\gamma(t_0)$ lies outside this set in (b), where $W(\gamma(t_0),t_0)>0$, and on its boundary in (a), (c), and (d).
 The insets show these positions in greater detail.}
  \label{fig:front-trajectories}
\end{figure}

\textbf{(e) Optimal trajectories and their uniqueness.}
  In the geometric optics limit, we characterize the optimal trajectories from the initial source to the propagation front.
  On the spherical front pieces, the trajectories are straight.
  On the paraboloidal, cylindrical, and conical parts, away from their junctions with the spherical pieces, each optimal trajectory consists of two straight segments.
  In forward time, it first moves from the origin along the shifting interface and then enters the lower-growth region, where it continues to the prescribed front point.
  We determine the time and location at which it leaves the interface.
  Representative optimal trajectories are illustrated in Figure \ref{fig:front-trajectories}.
  In the Freidlin regimes, the unconstrained optimal control problem and the min-max formula have the same unique minimizing trajectory at each point of the propagation front.
  In the non-Freidlin regime, the min-max minimizer is unique at every front point except in the relative interior of the planar part, where infinitely many minimizing trajectories exist; see Theorem \ref{thm:main-front-trajectories}.
  These planar front points can be reached through the higher-growth region at speeds strictly below $2\sqrt{r_{2}}$.
  This allows distinct trajectories to attain the same zero min-max value; see Remark \ref{remark:planar-front-nonuniqueness}.

  At the critical speed $c=2\sqrt{r_{2}-r_{1}}$, the optimal trajectories reaching the cylindrical front first follow the shifting interface and then move only in directions orthogonal to $\mathbf e$.
  Their velocity component orthogonal to $\mathbf e$ remains constant throughout and has norm $2\sqrt{r_{1}}$.
  Their endpoints at time $1$ therefore lie at the same distance $2\sqrt{r_{1}}$ from the $\mathbf e$-axis, while their axial positions are determined by the time at which the trajectories leave the shifting interface.
  These endpoints form the cylindrical part of the front.

\textbf{Outline of the paper.}
In Section \ref{section:main-results} we state the main theorems.
In Section \ref{section:singular-hje}, we prove convergence to the viscosity solution of the variational inequality and establish self-similarity.
In Section \ref{section:stopping-game} we construct the stopping game representation and connect its zero set to spreading.
In Section \ref{section:reductions-wulff-shapes} we study the unconstrained action and prove the sharp Freidlin dichotomy.
The min-max and pointwise zero-set reductions are derived in Section \ref{section:stopping-game-reductions}.
In Section \ref{section:explicit-wulff-shape} we compute the Wulff shapes and directional speeds and classify the optimal trajectories at the front.
In Section \ref{section:discussion}, we discuss extensions to one dimension and favorable--unfavorable environments.
The comparison principle and the technical game estimates are collected in Appendices \ref{app:comparison-structure} and \ref{app:bounded-data}.


\section{Main results}
\label{section:main-results}

In this section, we state our main results.
We first introduce the hyperbolic scaling and the Hopf-Cole transformation:
\begin{align*}
 u^{\varepsilon}(\mathbf{x},t)=u\left(\frac{\mathbf{x}}{\varepsilon},\frac{t}{\varepsilon}\right), \quad w^{\varepsilon}(\mathbf{x},t)=-\varepsilon \ln u^{\varepsilon}(\mathbf{x},t), \quad (\mathbf{x},t)\in \mathbb{R}^{N}\times\mathbb{R}^{+},
\end{align*}
where $\varepsilon\in(0,1)$. Substituting them into \eqref{eq:rd-nd-main} yields
\begin{align}
\label{eq:w-eps-pde-main}
  \partial_{t} w^{\varepsilon}-\varepsilon \Delta w^{\varepsilon}+|D w^{\varepsilon}|^{2}+r\left(\frac{\mathbf{x}\cdot \mathbf{e}-ct}{\varepsilon}\right)
 =u^{\varepsilon}(\mathbf{x},t).
\end{align}
We define a macroscopic potential step function $R(\mathbf{x},t)$ by:
\begin{equation}
\label{eq:def-R-macro}
  R(\mathbf{x},t)=
  \begin{cases}
    r_{1} & \text{if } \mathbf{x}\cdot \mathbf{e}<ct,\\
    r_{2} & \text{if } \mathbf{x}\cdot \mathbf{e} \geqslant ct.
  \end{cases}
\end{equation}

\subsection{Geometric optics limit and stopping-game representation}

\begin{theorem}
\label{thm:wkb-direct-convergence}
Let $\overline{w}(\mathbf{x},t)$ and $\underline{w}(\mathbf{x},t)$ be the half-relaxed limits of $w^{\varepsilon}$, i.e.,
\begin{align*}
\underline{w}(\mathbf{x},t):=\liminf_{\substack{\varepsilon \rightarrow 0\\(\mathbf{y},s)\rightarrow(\mathbf{x},t)}}w^{\varepsilon}(\mathbf{y},s) \quad \text { and } \quad \overline{w}(\mathbf{x},t):=\limsup_{\substack{\varepsilon \rightarrow 0\\(\mathbf{y},s)\rightarrow(\mathbf{x},t)}} w^{\varepsilon}(\mathbf{y},s).
\end{align*}
Then we have
\begin{align}
\label{eq:unique-vis-sol}
  \overline w(\mathbf{x},t)=\underline w(\mathbf{x},t)=:w(\mathbf{x},t) \quad \text{for all } (\mathbf{x},t)\in\mathbb{R}^{N}\times\mathbb{R}^{+}.
\end{align}
Moreover, $ w^{\varepsilon}\rightarrow w$ locally uniformly in $\mathbb{R}^{N}\times\mathbb{R}^{+}$ and $w$ is the unique nonnegative continuous viscosity solution of the following variational inequality
\begin{align}
\label{eq:obstacle-direct}
  \min\left\{\partial_{t}w+|Dw|^{2}+R(\mathbf{x},t), \; w \right\}=0, \quad (\mathbf{x},t)\in \mathbb{R}^{N}\times\mathbb{R}^{+},
\end{align}
in the class of solutions for which there exists $C>0$ such that
\begin{align}
\label{eq:singular-asymptotic-precise}
  \frac{|\mathbf{x}|^{2}}{4t}-Ct \leqslant w(\mathbf{x},t) \leqslant \frac{|\mathbf{x}|^{2}}{4t}+Ct, \quad (\mathbf{x},t)\in\mathbb{R}^{N}\times\mathbb{R}^{+}.
\end{align}
In particular, $w$ satisfies the following singular limiting behavior:
\begin{align}\label{eq:singular-initial-data}
  \lim_{t\downarrow0}w(\mathbf{x},t)=
  \begin{cases}
    0& \text{for }\mathbf{x}=\mathbf{0},\\
    +\infty &\text{for }\mathbf{x}\in \mathbb{R}^{N}\backslash\{\mathbf{0}\}.
  \end{cases}
\end{align}
\end{theorem}

\begin{remark}
\label{remark:relax-initial-condition}
The singular initial condition \eqref{eq:singular-initial-data} in Theorem \ref{thm:wkb-direct-convergence} is understood in the relaxed boundary sense and does not require the existence of a joint limit at $(\mathbf{0},0)$.
Indeed, \eqref{eq:singular-asymptotic-precise} implies
\begin{align*}
  \overline w_{0}(\mathbf{x}) &:= \limsup_{\substack{(\mathbf{y},s)\to(\mathbf{x},0)\\s>0}} w(\mathbf{y},s)=+\infty \quad \text{for } \mathbf{x}\in\mathbb{R}^{N},\\
  \underline w_{0}(\mathbf{x}) &:= \liminf_{\substack{(\mathbf{y},s)\to(\mathbf{x},0)\\s>0}} w(\mathbf{y},s) =
  \begin{cases}
    0 & \text{for } \mathbf{x}=\mathbf{0},\\
    +\infty & \text{for } \mathbf{x}\neq\mathbf{0}.
  \end{cases}
\end{align*}
Thus the joint limit is $+\infty$ at every $\mathbf{x}_{0}\neq\mathbf{0}$, whereas at $\mathbf{x}_{0}=\mathbf{0}$ it does not exist, although the pointwise limit $w(\mathbf{0},t)\to0$ holds as $t\downarrow0$.
\end{remark}

We next represent the geometric optics limit $w$ in Theorem \ref{thm:wkb-direct-convergence} by a zero-sum stopping game.
For $(\mathbf{x},t)\in\mathbb{R}^{N}\times\mathbb{R}^{+}$, define
\begin{align*}
  \mathbb{X}
 =H^{1}_{\textrm{loc}}([0,\infty);\mathbb{R}^{N}), \quad \mathbb{X}^{\mathbf{x},t}
 =\left\{\gamma\in\mathbb{X}\mid \gamma(0)=\mathbf{x},\, \gamma(t)=\mathbf{0} \right\}.
\end{align*}
A map $\vartheta:\mathbb{X}\to[0,\infty]$ is called a \textit{stopping time} if, for all $\gamma,\widehat\gamma\in\mathbb{X}$ and all $s\geqslant 0$,
\begin{align*}
  \gamma(\tau)=\widehat\gamma(\tau)\ \text{ for } \tau\in[0,s], \quad \vartheta[\gamma]\leqslant s \quad\Longrightarrow
  \quad \vartheta[\gamma]=\vartheta[\widehat\gamma].
\end{align*}
We denote by $\Theta$ the set of all such stopping times on $\mathbb{X}$, and define the value of the zero-sum stopping game by
\begin{align}
\label{eq:stopping-game-value}
  W(\mathbf{x},t)=\sup_{\vartheta\in\Theta}\inf_{\gamma\in\mathbb{X}^{\mathbf{x},t}} \mathcal{P}(\mathbf{x},t;\vartheta,\gamma),
\end{align}
where the payoff $\mathcal{P}$ is defined by
\begin{align}
\label{eq:payoff-and-lagrangian}
  \mathcal{P}(\mathbf{x},t;\vartheta,\gamma)=\int_{0}^{t\wedge\vartheta[\gamma]} L(\dot\gamma(s),\gamma(s),t-s)\,ds,
  \quad \text{with} \quad L(\mathbf{v},\mathbf{x},\tau)=\frac{1}{4}|\mathbf{v}|^{2}-R(\mathbf{x},\tau),
\end{align}

\begin{theorem}
\label{thm:stopping-game-singular-obstacle-solution}
The value function $W$ defined by \eqref{eq:stopping-game-value} is the unique nonnegative continuous viscosity solution of the variational inequality \eqref{eq:obstacle-direct} satisfying the growth control condition \eqref{eq:singular-asymptotic-precise}.
Consequently, $W=w$ in $\mathbb{R}^{N}\times\mathbb{R}^{+}$, where $w$ is the geometric optics limit given by \eqref{eq:unique-vis-sol}.
\end{theorem}

The stopping value also admits the following equivalent pathwise min-max representation.

\begin{proposition}[Min-max representation of the stopping value]
\label{prop:stopping-value-minmax}
For every $(\mathbf{x},t)\in\mathbb{R}^{N}\times\mathbb{R}^{+}$,
\begin{align*}
  W(\mathbf{x},t)=\min_{\gamma\in\mathbb{X}^{\mathbf{x},t}} \max_{\sigma\in[0,t]} \int_{0}^{\sigma} L(\dot\gamma(s),\gamma(s),t-s)\,ds.
\end{align*}
\end{proposition}

\begin{remark}[Relation to the conjecture in \cite{lam2025}]
\label{remark:minmax-lam-conjecture}
Theorem \ref{thm:stopping-game-singular-obstacle-solution} and Proposition \ref{prop:stopping-value-minmax} remain valid for $N=1$, with the same proofs.
In particular, neither the stopping-game representation nor the min-max formula requires $N\geqslant2$.
To compare with \cite[Section~2.5]{lam2025}, take $N=1$, $\mathbf{e}=1$. 
For the nondecreasing $r(\cdot)$ considered here, the generalized principal eigenvalue $\Lambda_{1}$ defined as
\begin{align*}
  \Lambda_{1}=\inf\left\{\Lambda\in\mathbb{R}\mid \exists\,\phi\in C^{2}_{\mathrm{loc}}(\mathbb{R}),\; \phi>0,\, \text{s.t.}\;\phi''+r(y)\phi\leqslant\Lambda\phi\;\text{in}\; \mathbb{R}\right\}
\end{align*}
 satisfies $\Lambda_1=r_{2}$.
Indeed, \cite[Proposition~4.2(a)]{lam2025} implies $\Lambda_{1}\geqslant r_{2}$, while choosing $\phi\equiv1$ in the definition gives $\Lambda_{1}\leqslant r_{2}$.
Hence the cost on the shifting interface coincides with that in \eqref{eq:payoff-and-lagrangian}.
Their formulation deals with $x\geqslant0$ and uses nonnegative trajectories.
For $x\geqslant0$, taking the positive part of an admissible trajectory cannot increase its action at any intermediate time, so restricting to nonnegative trajectories does not change the min-max value.
Together with $w=W$, this confirms the representation conjectured in \cite[equation (2.22)]{lam2025} for a single increasing transition between two favorable regions.
\end{remark}

The value function $W$ also admits the following self-similarity.

\begin{proposition}[Self-similarity of the stopping game]
\label{prop:self-similarity}
The common function $W=w$ satisfies
\begin{align}
\label{eq:self-similarity}
  W(k\mathbf{x},kt)=kW(\mathbf{x},t) \quad \text{for all} \quad (\mathbf{x},t)\in\mathbb{R}^{N}\times\mathbb{R}^{+} \; \text{and} \;k>0.
\end{align}
\end{proposition}

The self-similarity of the stopping value function determines its space--time zero set through the unit time scaling, which governs the large time propagation of $u$, as stated next.

\begin{theorem}\label{thm:wulff}
If we denote the Wulff shape by 
\begin{align*}
  \mathcal{W}=\left\{\mathbf{x}\in \mathbb{R}^{N}\mid W(\mathbf{x},1)=0\right\},
\end{align*}
then for all sufficiently small $\eta>0$, it holds that
\begin{equation*}
  \begin{cases}
    \underset{t\rightarrow\infty} \lim\,\underset{\mathrm{dist}\{\frac{\mathbf{x}}{t},\mathcal{W}\}>\eta}\sup\, u(\mathbf{x},t)=0, \\
    \underset{t\rightarrow\infty} \liminf\,\underset{\mathrm{dist}\{\frac{\mathbf{x}}{t},\overline{\mathcal{W}}^{c}\}>\eta}\inf\, u(\mathbf{x},t)>0.
  \end{cases}
\end{equation*}
\end{theorem}

\subsection{Freidlin reduction and convex hull characterization}
\label{subsection:freidlin-reduction}

We next determine when the stopping game reduces to an unconstrained optimal control and describe the complementary regime through a pathwise characterization of the zero set.
For $(\mathbf{x},t)\in\mathbb{R}^{N}\times\mathbb{R}^{+}$, set
\begin{align}
\label{eq:hilbert-B}
  \mathbb B^{\mathbf{x},t}=\left\{\gamma\in H^{1}([0,t];\mathbb{R}^{N}) \mid \gamma(0)=\mathbf{x},\, \gamma(t)=\mathbf{0} \right\}.
\end{align}
For $\gamma\in\mathbb B^{\mathbf{x},t}$, define the unconstrained optimal control value function as
\begin{align}\label{eq:optimal-control-J}
  J(\mathbf{x},t)=\inf_{\gamma\in\mathbb B^{\mathbf{x},t}} \int_{0}^{t} L(\dot\gamma(s),\gamma(s),t-s)\,ds
\end{align}
Restriction to $[0,t]$ maps $\mathbb{X}^{\mathbf{x},t}$ onto $\mathbb B^{\mathbf{x},t}$.
Thus $J$ is precisely the value obtained by requiring the stopping player to wait until the terminal time. Let
\begin{align*}
  \mathbb P=\left\{(\mathbf{x},t)\in\mathbb{R}^{N}\times\mathbb{R}^{+} \mid J(\mathbf{x},t)>0 \right\}.
\end{align*}

\begin{definition}[Freidlin's condition]
\label{def:backward-freidlin-condition}
We say that $J$ satisfies \textit{Freidlin's condition} if, for every $(\mathbf{x},t)\in\partial\mathbb P$,
\begin{align}
\label{eq:backward-freidlin-condition}
  J(\mathbf{x},t)=\inf_{\gamma\in\mathbb{X}^{\mathbf{x},t}} \left\{\int_{0}^{t} L(\dot\gamma(s),\gamma(s),t-s)\,ds \,\middle|\, (\gamma(s),t-s)\in\mathbb P \text{ for every }s\in(0,t) \right\}.
\end{align}
\end{definition}

The following theorem gives both the reduction principle and the sharp range of validity of Freidlin's condition.
From now on, we use $\alpha=\sqrt{r_{2}-r_{1}}$ throughout the paper.

\begin{theorem}
\label{thm:main-freidlin-dichotomy}
The following assertions hold.
\begin{enumerate}[label=(\roman*),font=\normalfont, leftmargin=*]
  \item If $J$ satisfies Freidlin's condition, then
  \begin{align}
  \label{eq:main-freidlin-reduction}
  W(\mathbf{x},t)=\max\{0,J(\mathbf{x},t)\} \quad \text{in } \mathbb{R}^{N}\times\mathbb{R}^{+}.
  \end{align}

  \item If $c\leqslant-2\sqrt{r_{2}}$ or $c\geqslant2\alpha$, then $J$ satisfies Freidlin's condition.
  Moreover, we have
  \begin{align*}
  \mathcal{W}=\left\{\mathbf{x}\in\mathbb{R}^{N} \mid J(\mathbf{x},1)\leqslant0 \right\}.
  \end{align*}

  \item If $-2\sqrt{r_{2}}<c<2\alpha$, then $J$ does not satisfy Freidlin's condition. Moreover, by setting
  \begin{align*}
  D_{1}
 =\left\{\mathbf{x}\in\mathbb{R}^{N} \mid |\mathbf{x}|\leqslant2\sqrt{r_{1}} \right\},
  \quad D_{2}=\left\{\mathbf{x}\in\mathbb{R}^{N} \mid |\mathbf{x}|\leqslant2\sqrt{r_{2}}, \quad \mathbf{x}\cdot\mathbf{e}\geqslant c \right\},
  \end{align*}
  we have
  \begin{align*}
    \mathcal{W}=\mathrm{co}\left(D_{1}\cup D_{2}\right)\subsetneq\left\{\mathbf{x}\in\mathbb{R}^{N}\mid J(\mathbf{x},1)\leqslant0\right\}.
\end{align*}
  Thus the zero set of $\max\{0,J(\cdot,1)\}$ strictly contains the Wulff shape.
\end{enumerate}
\end{theorem}

\begin{remark}
\label{remark:scope-unconstrained-reduction}
Under the present assumptions with $N\geqslant2$, Theorem \ref{thm:main-freidlin-dichotomy} shows that Freidlin's condition is necessary and sufficient for the unconstrained reduction.
For $N=1$, however, the reduction remains valid at $c=-2\sqrt{r_{2}}$ despite the failure of Freidlin's condition; see Subsection \ref{subsec:discussion-1d}.
\end{remark}

\subsection{Explicit Wulff shapes and directional spreading speeds}
\label{subsection:explicit-wulff-shapes}

We now make the characterizations in Theorem \ref{thm:main-freidlin-dichotomy} explicit.
Throughout the rest of the paper, set $\lambda=\frac{c}{2}-\alpha$ and
\begin{align}
\label{eq:orthogonal-decompo}
  \xi=\mathbf{x}\cdot\mathbf{e},\quad {\boldsymbol\zeta}=\mathbf{x}-\xi\mathbf{e}, \quad \rho=|\mathbf{x}-\xi\mathbf{e}|.
\end{align}
For $t\in\mathbb{R}^{+}$ and $q\in\mathbb{R}$, define
\begin{align*}
  \Sigma_{1}^{q}(t) &= \left\{\mathbf{x}\in\mathbb{R}^{N} \mid |\mathbf{x}|\leqslant2\sqrt{r_{1}}t, \; \xi\leqslant qt \right\},\\
  \Sigma_{2}^{q}(t) &= \left\{\mathbf{x}\in\mathbb{R}^{N} \mid |\mathbf{x}|\leqslant2\sqrt{r_{2}}t, \; \xi\geqslant qt \right\},\\
  \Sigma_{\mathrm{m}}^{\lambda,q}(t) &= \left\{\mathbf{x}\in\mathbb{R}^{N} \mid 2\lambda t\leqslant\xi\leqslant qt, \, \rho^{2}\leqslant 4t\bigl[(r_{1}+\lambda^{2})t-\lambda\xi\bigr] \right\},
  \quad \text{if } \lambda\geqslant0,\; q\geqslant2\lambda.
\end{align*}

We first state the formulas in the Freidlin regimes.

\begin{theorem}[Explicit Wulff shapes in the Freidlin regimes]
\label{thm:main-wulff-freidlin}
The following assertions hold.
\begin{enumerate}[label=(\roman*),font=\normalfont, leftmargin=*]
  \item If $c\geqslant2(\alpha+\sqrt{r_{1}})$, then
  \begin{align*}
  \mathcal{W}= \overline{B_{2\sqrt{r_{1}}}(\mathbf{0})}.
  \end{align*}

  \item If $2\sqrt{r_{2}}<c<2(\alpha+\sqrt{r_{1}})$, set $c_\mathrm{nlp}=\lambda+\frac{r_{1}}{\lambda}$. Then $2\lambda<c_\mathrm{nlp}<c$ and
  \begin{align*}
  \mathcal{W}=\Sigma_{1}^{2\lambda}(1) \cup \Sigma_{\mathrm{m}}^{\lambda,c_\mathrm{nlp}}(1).
  \end{align*}

  \item  If $2\alpha\leqslant c\leqslant2\sqrt{r_{2}}$, then
\begin{align*}
  \mathcal{W}=\Sigma_{1}^{2\lambda}(1) \cup \Sigma_{\mathrm{m}}^{\lambda,c}(1) \cup \Sigma_{2}^{c}(1).
\end{align*}

  \item If $c\leqslant-2\sqrt{r_{2}}$, then
  \begin{align*}
  \mathcal{W}=\overline{B_{2\sqrt{r_{2}}}(\mathbf{0})}.
  \end{align*}
\end{enumerate}
In particular, if $c=2\alpha$, then $\lambda=0$ and the intermediate part of the Wulff shape is a solid cylinder of radius $2\sqrt{r_{1}}$ extending from $\xi=0$ to $\xi=2\alpha$, i.e., $\Sigma_{\mathrm{m}}^{0,2\alpha}(1)=\left\{\mathbf{x}\in\mathbb{R}^{N} |\, 0\leqslant\xi\leqslant2\alpha,\rho\leqslant2\sqrt{r_{1}} \right\}$.
\end{theorem}

We next make the convex hull characterization in Theorem \ref{thm:main-freidlin-dichotomy}-(iii) explicit.

\begin{theorem}[Explicit Wulff shapes in the non-Freidlin regime]
\label{thm:main-wulff-nonfreidlin}
Assume that $-2\sqrt{r_{2}}<c<2\alpha$.
Then the following assertions hold.
\begin{enumerate}[label=(\roman*),font=\normalfont, leftmargin=*]
  \item If $-2\sqrt{r_{2}}<c\leqslant-2\sqrt{r_{1}}$, then
  \begin{align*}
  \mathcal{W}=\Sigma_{2}^{c}(1).
  \end{align*}

  \item If $-2\sqrt{r_{1}}<c<2\alpha$, define
  \begin{align}
  \label{eq:contact-point}
  \xi_{*}=\frac{\sqrt{r_{1}}}{r_{2}} \left(c\sqrt{r_{1}}-\alpha\sqrt{4r_{2}-c^{2}} \right),
  \quad \rho_{*}=\frac{\sqrt{r_{1}}}{r_{2}} \left(\sqrt{r_{1}}\sqrt{4r_{2}-c^{2}}+\alpha c \right).
  \end{align}
  Then $-2\sqrt{r_{1}}<\xi_{*}<c$ and $\rho_{*}>0$, and
  \begin{align*}
  \mathcal{W}=\Sigma_{1}^{\xi_{*}}(1) \cup \left\{\mathbf{x}\in\mathbb{R}^{N} \mid \xi_{*}\leqslant\xi\leqslant c, \; \xi_{*}\xi+\rho_{*}\rho\leqslant4r_{1} \right\} \cup \Sigma_{2}^{c}(1).
\end{align*}
\end{enumerate}
\end{theorem}

Theorems \ref{thm:main-wulff-freidlin} and \ref{thm:main-wulff-nonfreidlin} show that $\mathcal{W}$ is a compact convex set containing $\mathbf{0}$ in its interior in every speed regime.
Hence each ray from the origin meets $\partial\mathcal{W}$ at a unique point, which leads to the following directional characterization.

\begin{corollary}[Characterization of the directional spreading speed]
\label{cor:wulff-directional-spreading}
For $\mathbf{n}\in\mathbb S^{N-1}$, define
\begin{align*}
  c_\mathrm{dir}(\mathbf{n})=\sup \left\{s\geqslant0 \mid s\mathbf{n}\in\mathcal{W} \right\}.
\end{align*}
Then $c_\mathrm{dir}(\mathbf{n})\in[2\sqrt{r_{1}},2\sqrt{r_{2}}]$ and $c_\mathrm{dir}(\mathbf{n})\mathbf{n}\in\partial\mathcal{W}$.
Moreover, for any nonempty compact set $K\subset\mathbb{R}^{N}$, we have
\begin{equation*}
  \begin{cases}
    \underset{t\to\infty}\liminf \underset{\mathbf{x}\in K}\inf\, u(\mathbf{x}+qt\mathbf{n},t) >0, & \text{if }q\in[0,c_\mathrm{dir}(\mathbf{n})),\\
    \underset{t\to\infty}\lim\, \underset{\mathbf{x}\in K}\sup\, u(\mathbf{x}+qt\mathbf{n},t) =0,&  \text{if }q\in(c_\mathrm{dir}(\mathbf{n}),\infty).
  \end{cases}
\end{equation*}
\end{corollary}

The preceding corollary reduces the computation of $c_\mathrm{dir}(\mathbf{n})$ to locating the intersection of the ray $\{s\mathbf{n}\mid s\geqslant0\}$ with the explicit boundaries above.
This yields the following formulas.

\begin{corollary}[Explicit formulas of the directional spreading speeds]
\label{cor:main-directional-speed}
For $\mathbf{n}\in\mathbb S^{N-1}$, set $\beta=\mathbf{n}\cdot\mathbf{e}$.
Then $\beta\in[-1,1]$ and the explicit directional spreading speeds are given as follows.
\begin{alignat*}{2}
  \textup{(i)}\quad &\text{if } c\geqslant2(\alpha+\sqrt{r_{1}}) \text{, then} \qquad &&c_\mathrm{dir}(\mathbf{n}) =2\sqrt{r_{1}}, \\
  \textup{(ii)}\quad &\text{if }2\sqrt{r_{2}}<c<2(\alpha+\sqrt{r_{1}})\text{, then}\qquad &&c_\mathrm{dir}(\mathbf{n}) =
  \begin{cases}
    2\sqrt{r_{1}} &\text{if }\beta\leqslant\frac{\lambda}{\sqrt{r_{1}}},\\
    c_{\lambda}(\beta) &\text{if }\beta>\frac{\lambda}{\sqrt{r_{1}}},
  \end{cases}\\
  \textup{(iii)}\quad &\text{if }2\alpha\leqslant c\leqslant2\sqrt{r_{2}}\text{, then}\qquad &&c_\mathrm{dir}(\mathbf{n}) =
  \begin{cases}
    2\sqrt{r_{1}}, &\text{if } \beta\leqslant\frac{\lambda}{\sqrt{r_{1}}},\\
    c_{\lambda}(\beta), &\text{if } \frac{\lambda}{\sqrt{r_{1}}} <\beta< \frac{c}{2\sqrt{r_{2}}},\\
    2\sqrt{r_{2}}, &\text{if } \beta\geqslant\frac{c}{2\sqrt{r_{2}}},
  \end{cases}\\
  \textup{(iv)}\quad &\text{if }-2\sqrt{r_{1}}<c<2\alpha\text{, then}\qquad &&c_\mathrm{dir}(\mathbf{n}) =
  \begin{cases}
    2\sqrt{r_{1}}, &\text{if }\beta\leqslant \frac{\xi_{*}}{2\sqrt{r_{1}}},\\[2pt]
    c_{*}(\beta) &\text{if } \frac{\xi_{*}}{2\sqrt{r_{1}}} <\beta< \frac{c}{2\sqrt{r_{2}}},\\[2pt]
    2\sqrt{r_{2}}, &\text{if } \beta\geqslant \frac{c}{2\sqrt{r_{2}}},
  \end{cases}\\
  \textup{(v)}\quad &\text{if }-2\sqrt{r_{2}}<c\leqslant -2\sqrt{r_{1}}\text{, then}\qquad &&c_\mathrm{dir}(\mathbf{n}) =
  \begin{cases}
    \frac{c}{\beta}, &\text{if } \beta<\frac{c}{2\sqrt{r_{2}}},\\
    2\sqrt{r_{2}}, &\text{if } \beta\geqslant \frac{c}{2\sqrt{r_{2}}},
  \end{cases}\\
  \textup{(vi)}\quad &\text{if }c\leqslant -2\sqrt{r_{2}}\text{, then}\qquad &&c_\mathrm{dir}(\mathbf{n})=2\sqrt{r_{2}},
\end{alignat*}
where $\xi_{*}$ in (iv) is given by \eqref{eq:contact-point} and $c_{\lambda}(\beta)$ in (ii)-(iii) and $c_{*}(\beta)$ in (iv) are defined as
\begin{align*}
  c_{\lambda}(\beta)=\frac{2(r_{1}+\lambda^{2})} {\sqrt{\lambda^{2}+r_{1}(1-\beta^{2})} +\lambda\beta},
  \quad c_{*}(\beta)=\frac{4r_{1}} {\xi_{*}\beta+ \rho_{*}\sqrt{1-\beta^{2}}}.
\end{align*}
\end{corollary}

\subsection{Optimal trajectories at the propagation front}
\label{subsection:front-trajectories}

We now characterize the optimal trajectories joining $\mathbf{x}\in\partial\mathcal{W}$ at time $1$ to the origin at time $0$. Uniqueness is understood on $[0,1]$.

\begin{theorem}[Optimal trajectories at the propagation front]
\label{thm:main-front-trajectories}
Fix $c\in\mathbb{R}$ and $\mathbf{x}\in\partial\mathcal{W}$.
\begin{enumerate}[label=(\roman*),font=\normalfont, leftmargin=*]
  \item \textbf{Freidlin regimes.}
  If $c\leqslant-2\sqrt{r_{2}}$ or $c\geqslant2\alpha$, then the unconstrained optimal control problem defining $J(\mathbf{x},1)$ and the min-max formula \eqref{eq:stopping-value-inf-action} have the same unique minimizing trajectory $\gamma_{*}$ on $[0,1]$.

  \item \textbf{Non-Freidlin regime.}
  If $-2\sqrt{r_{2}}<c<2\alpha$, then the min-max formula in Proposition \ref{prop:stopping-value-minmax} admits a minimizing trajectory at $(\mathbf{x},1)$.
  Precisely, there are infinitely many minimizing trajectories when
  \begin{align*}
    -2\sqrt{r_{2}}<c\leqslant-2\sqrt{r_{1}}, \qquad \xi=c, \qquad \rho<\sqrt{4r_{2}-c^{2}};
  \end{align*}
  that is, in the relative interior of the planar part of the front.
  At every other point of $\partial\mathcal{W}$, including all spherical parts and their junctions with the other front pieces, the minimizing trajectory is unique.
\end{enumerate}

The unique trajectories in \textup{(i)--(ii)} have the following explicit forms.
On the spherical parts of $\partial\mathcal{W}$, they are given by
\begin{align}
  \label{eq:front-straight-optimal-path}
  \gamma_{*}(s)=(1-s)\mathbf{x}, \quad s\in[0,1].
\end{align}
On the paraboloidal and cylindrical parts in the Freidlin regimes with $2\lambda<\xi<c$, set $\sigma_*=\frac{c-\xi}{2\alpha}$.
Then
\begin{align}
  \label{eq:freidlin-front-optimal-path}
  \gamma_{*}(s)=
  \begin{cases}
    (\xi-2\lambda s)\mathbf{e}+(1-s){\boldsymbol\zeta}, & s\in[0,\sigma_*],\\
    (1-s)(c\mathbf{e}+{\boldsymbol\zeta}), & s\in[\sigma_*,1].
  \end{cases}
\end{align}
On the conical part in the non-Freidlin regime, use $\xi_*$ and $\rho_{*}$ from \eqref{eq:contact-point} and set
\begin{align}
  \label{eq:cone-front-path-parameters}
  \mathbf{n}=\frac{{\boldsymbol\zeta}}{\rho}, \qquad
  P=c\mathbf{e}+\sqrt{4r_{2}-c^{2}}\,\mathbf{n}, \qquad
  Q=\xi_*\mathbf{e}+\rho_{*}\mathbf{n}, \qquad
  \tau_*=\frac{c-\xi}{c-\xi_*}.
\end{align}
Then
\begin{align}
  \label{eq:cone-front-optimal-path}
  \gamma_{*}(s)=
  \begin{cases}
    \mathbf{x}-sQ, & s\in[0,\tau_*],\\
    (1-s)P, & s\in[\tau_*,1].
  \end{cases}
\end{align}
Both two-segment formulas reduce to \eqref{eq:front-straight-optimal-path} at their junctions with the spherical parts.
The unique optimal trajectories have the following geometric characterization.
On the paraboloidal part in the Freidlin regime with $c>2\alpha$ and $2\lambda<\xi<c$, we have
\begin{align}
W(\gamma_{*}(s),1-s)>0, \qquad s\in(0,1),
\end{align}
so the optimal trajectory remains in the uninvaded region at every intermediate time.
At all other points where the minimizing trajectory is unique, it remains on the propagation front:
\begin{align}
\label{eq:main-optimal-path-on-front}
\gamma_{*}(s)\in (1-s)\partial\mathcal W,
\qquad
W(\gamma_{*}(s),1-s)=0, \qquad s\in[0,1).
\end{align}
Moreover, every unique minimizing trajectory in the non-Freidlin regime has zero running cost almost everywhere.
\end{theorem}

\begin{remark}[Conical fronts and nonlocal speed selection]
\label{remark:cone-front-nonlocal-selection}
For a point strictly between the two junctions of the conical front, the forward trajectory $s\mapsto\gamma_{*}(1-s)$ follows the shifting interface for time $1-\tau_*$ at a speed $2\sqrt{r_{2}}$ and then enters the lower-growth region at a speed $2\sqrt{r_{1}}$.
Since $|Q|^{2}=Q\cdot P=4r_{1}$ and $|P-Q|^{2}=4\alpha^{2}$, its endpoint satisfies
\begin{align*}
  |\mathbf{x}|^{2}=4r_{1}+4\alpha^{2}(1-\tau_*)^{2}>4r_{1}.
\end{align*}
Thus the directional spreading speed exceeds $2\sqrt{r_{1}}$, although the straight forward trajectory $s\mapsto s\mathbf{x}$ lies in the lower-growth region for $s\in(0,1)$.
This extends the interpretation of nonlocal speed selection in \cite[Section 5]{lam2027} from a stationary two-dimensional environment to shifting interfaces in any dimension $N\geqslant2$.
The conical trajectories remain on the propagation front, whereas the optimal trajectories reaching the paraboloidal part away from its junctions satisfy $W>0$ at all intermediate times, as stated in Theorem~\ref{thm:main-front-trajectories}.
\end{remark}

\begin{remark}
\label{remark:planar-front-nonuniqueness}
The nonuniqueness on the planar part reflects the strict inequality $|\mathbf{x}|<2\sqrt{r_{2}}$.
These points lie on the shifting interface but strictly inside the ball of radius $2\sqrt{r_{2}}$.
We can therefore perturb the straight trajectory into the higher-growth region while preserving its endpoints and keeping its speed below $2\sqrt{r_{2}}$.
For the trajectories constructed in the proof of Proposition \ref{prop:nonfreidlin-front-trajectories}, the action is negative at every positive backward time, so its maximum remains $0$, attained at time $0$.
These trajectories have the same min-max value, although their total actions differ.
At the junction with the spherical part, the equality $|\mathbf{x}|=2\sqrt{r_{2}}$ removes this freedom, and the Cauchy-Schwarz equality condition forces the unique straight trajectory.
\end{remark}

\section{Geometric optics limit and the associated variational inequality}
\label{section:singular-hje}

In this section, we derive the geometric optics limit of \eqref{eq:rd-nd-main}.
We first establish local estimates and identify the half-relaxed limits as viscosity sub- and supersolutions of the limiting variational inequality.
We then prove comparison in the singular growth class, obtain local uniform convergence and uniqueness, and finally establish the self-similarity of the limit.

For equations with discontinuous Hamiltonians, we use the notion of a viscosity solution in the discontinuous setting \cite{achdou2013}. We also refer to \cite{barles1987} for the introduction of the half-relaxed limit method.

\subsection{Preliminary estimates and half-relaxed limits}

\begin{lemma}
\label{lem:local-bound-nd}
Let $u(\mathbf{x},t)$ be the solution to \eqref{eq:rd-nd-main}.
For each compact subset $Q \subset \mathbb{R}^{N} \times \mathbb{R}^{+}$, there exists a constant $C(Q)>0$ independent of $\varepsilon$ such that $w^{\varepsilon}(\mathbf{x},t) \leqslant C(Q)$ for all $(\mathbf{x},t) \in Q$ and sufficiently small $\varepsilon>0$.
\end{lemma}

\begin{proof}
By the standard parabolic comparison principle, the solution satisfies $0 \leqslant u(\mathbf{x},t) \leqslant M$, where $M=\max\{r_{2}, \|g\|_{L^{\infty}}\}$.
Applying the Hopf-Cole transformation $w^{\varepsilon}=-\varepsilon\ln u^{\varepsilon}$, we obtain the global lower bound $w^{\varepsilon} \geqslant -\varepsilon\ln M$ in $\mathbb{R}^{N}\times\mathbb{R}^{+}$.
To finish the proof, it suffices to establish the local upper bound.

Recall that the environmental capacity satisfies $r(\cdot) \geqslant r_{1}$ globally.
Thus, the nonlinear reaction term admits a uniform linear lower bound:
\begin{align*}
  u[r(\mathbf{x}\cdot\mathbf{e}-ct)-u] \geqslant u(r_{1}-M):=-C_{0}u,
\end{align*}
where $C_{0} =M-r_{1} \geqslant 0$ is a constant independent of $t$, $\mathbf{x}$, and $\varepsilon$.
This indicates that $u(\mathbf{x},t)$ is a global supersolution to the linear heat equation:
\begin{align}
\label{eq:linear-heat-equation}
  \partial_{t}v-\Delta v =-C_{0}v \quad \text{in } \mathbb{R}^{N} \times \mathbb{R}^{+}.
\end{align}
Since the initial data $g \geqslant 0$ is non-trivial and continuous, there exists a ball $B_{\delta}(\mathbf{x}_{0}) \subset \mathbb{R}^{N}$ and a constant $\eta>0$ such that $g(\mathbf{x}) \geqslant \eta$ on $B_{\delta}(\mathbf{x}_{0})$.
By the comparison principle with the fundamental solution of \eqref{eq:linear-heat-equation}, we have
\begin{align*}
  u(\mathbf{x},t) \geqslant \frac{e^{-C_{0} t}}{(4\pi t)^{N/2}} \int_{B_{\delta}(\mathbf{x}_{0})} e^{-\frac{|\mathbf{x}-\mathbf{y}|^{2}}{4t}} \eta \, d\mathbf{y}.
\end{align*}
For any $\mathbf{y} \in B_{\delta}(\mathbf{x}_{0})$, since $|\mathbf{x}-\mathbf{y}| \leqslant |\mathbf{x}-\mathbf{x}_{0}|+|\mathbf{x}_{0}-\mathbf{y}| \leqslant |\mathbf{x}-\mathbf{x}_{0}|+\delta$,
\begin{align*}
  u(\mathbf{x},t) \geqslant \frac{\eta e^{-C_{0} t}}{(4\pi t)^{N/2}} \int_{B_{\delta}(\mathbf{x}_{0})} \exp\left(-\frac{(|\mathbf{x}-\mathbf{x}_{0}|+\delta)^{2}}{4t}\right) d\mathbf{y}
  =\frac{\eta |B_{\delta}(\mathbf{x}_{0})|}{(4\pi t)^{N/2}} \exp\left(-C_{0} t-\frac{(|\mathbf{x}-\mathbf{x}_{0}|+\delta)^{2}}{4t} \right).
\end{align*}
Rescaling to $u^{\varepsilon}(\mathbf{x},t)=u(\mathbf{x}/\varepsilon, t/\varepsilon)$ and applying $w^{\varepsilon}=-\varepsilon \ln u^{\varepsilon}$, we obtain
\begin{equation}
\label{eq:upper-bound-wep-optimal}
  \begin{aligned}
    w^{\varepsilon}(\mathbf{x},t) &\leqslant C_{0} t+\frac{\varepsilon^{2}}{4t} \left(\frac{|\mathbf{x}|}{\varepsilon}+|\mathbf{x}_{0}|+\delta \right)^{2}-\varepsilon \ln \left(\frac{\eta |B_{\delta}(\mathbf{x}_{0})|\varepsilon^{N/2}}{(4\pi t)^{N/2}}\right) \\
    &= C_{0} t+\frac{1}{4t} \left(|\mathbf{x}|+\varepsilon(|\mathbf{x}_{0}|+\delta)\right)^{2}-\varepsilon \ln \left(\frac{\eta |B_{\delta}(\mathbf{x}_{0})|\varepsilon^{N/2}}{(4\pi t)^{N/2}}\right).
  \end{aligned}
\end{equation}

Now, let $Q \subset \mathbb{R}^{N}\times \mathbb{R}^{+}$ be an arbitrary compact set.
There exist some constants $0<\tau<1$ and $A>0$ such that $Q \subset\overline{B_{A}(\mathbf{0})}  \times [\tau, 1/\tau]$.
As $\varepsilon\to0$, the spatial term converges uniformly on $Q$:
\begin{align*}
  \frac{1}{4t} \big( |\mathbf{x}|+\varepsilon(|\mathbf{x}_{0}|+\delta) \big)^{2} \longrightarrow \frac{|\mathbf{x}|^{2}}{4t} \leqslant \frac{A^{2}}{4\tau}.
\end{align*}
Meanwhile, the logarithmic term tends to zero uniformly for $t\in[\tau,1/\tau]$ as $\varepsilon\to0$.
Therefore, we can find a sufficiently small $\varepsilon_{0}>0$ such that for all $0<\varepsilon \leqslant \varepsilon_{0}$,
\begin{align*}
  w^{\varepsilon}(\mathbf{x},t) \leqslant \frac{4C_{0}+A^{2}}{4\tau}+1 =: C(Q) \quad \text{for all } (\mathbf{x},t) \in Q.
\end{align*}
The proof is complete.
\end{proof}

By Lemma \ref{lem:local-bound-nd}, the following half-relaxed limits are well-defined in $\mathbb{R}^{N}\times\mathbb{R}^{+}$:
\begin{align*}
  \overline{w}(\mathbf{x},t)=\limsup_{\substack{\varepsilon \to 0\\  (\mathbf{y},s)\to(\mathbf{x},t) }} w^{\varepsilon}(\mathbf{y},s),
  \quad \underline{w}(\mathbf{x},t)=\liminf_{\substack{ \varepsilon \to 0 \\ (\mathbf{y},s)\to(\mathbf{x},t)}} w^{\varepsilon}(\mathbf{y},s).
\end{align*}
The heat-kernel estimate also yields the following lower bound.

\begin{lemma}\label{lem:initial-lower-bound}
For any $t\in \mathbb{R}^{+}$, the lower half-relaxed limit of the macroscopic action $\underline{w}(\mathbf{x}, t)$ satisfies the following lower bound:
\begin{align}\label{eq:optimal-sub-estimate}
  \underline{w}(\mathbf{x}, t) \geqslant \frac{|\mathbf{x}|^{2}}{4t}-r_{2} t.
\end{align}
\end{lemma}

\begin{proof}
By the maximum principle, $u^{\varepsilon} \geqslant 0$ globally.
Since $r(\cdot) \leqslant r_{2}$, we know that $u^{\varepsilon}$ is a classical subsolution to the following linear parabolic equation:
\begin{align*}
  \partial_{t} v^{\varepsilon}-\varepsilon \Delta v^{\varepsilon}=\frac{r_{2}}{\varepsilon} v^{\varepsilon}
  \quad \text{in } \mathbb{R}^{N} \times \mathbb{R}^{+}.
\end{align*}
Applying the parabolic comparison principle, $u^{\varepsilon}$ is globally bounded from above by the fundamental solution of this linear equation initialized with $g(\frac{\mathbf{x}}{\varepsilon})$:
\begin{align*}
  u^{\varepsilon}(\mathbf{x}, t) \leqslant \frac{e^{r_{2} t / \varepsilon}}{(4\pi \varepsilon t)^{N/2}} \int_{\mathbb{R}^{N}} e^{-\frac{|\mathbf{x}-\mathbf{y}|^{2}}{4\varepsilon t}} g \left(\frac{\mathbf{y}}{\varepsilon}\right) \, d\mathbf{y}.
\end{align*}

Set the change of variables $\mathbf{z}=\frac{\mathbf{y}}{\varepsilon}$. Then
\begin{align}
\label{eq:integral-substituted}
  u^{\varepsilon}(\mathbf{x}, t) \leqslant \frac{\varepsilon^{N/2} e^{r_{2} t / \varepsilon}}{(4\pi t)^{N/2}} \int_{S} e^{-\frac{|\mathbf{x}-\varepsilon\mathbf{z}|^{2}}{4\varepsilon t}} g(\mathbf{z}) \, d\mathbf{z},
\end{align}
where $S=\mathrm{supp}(g)$ is the compact support of initial data $g \in C(\mathbb{R}^{N})$ with its Lebesgue measure satisfying $|S|>0$.
Let $R_{S}=\sup_{\mathbf{z} \in S} |\mathbf{z}|<\infty$.
For any $\mathbf{z} \in S$, it holds that
\begin{align*}
  |\mathbf{x}-\varepsilon \mathbf{z}| \geqslant \max\left(0, |\mathbf{x}|-\varepsilon R_{S}\right).
\end{align*}
Then we obtain
\begin{align*}
  \int_{S} e^{-\frac{|\mathbf{x}-\varepsilon\mathbf{z}|^{2}}{4\varepsilon t}} g(\mathbf{z}) \, d\mathbf{z}
  \leqslant  \int_{S} M \exp\left(-\frac{\max(0, |\mathbf{x}|-\varepsilon R_{S})^{2}}{4\varepsilon t} \right) d\mathbf{z}
  =M|S| \exp\left(-\frac{\max(0, |\mathbf{x}|-\varepsilon R_{S})^{2}}{4\varepsilon t} \right),
\end{align*}
where $M=\max\{r_{2}, \|g\|_{L^{\infty}}\}\geqslant \|g\|_{L^{\infty}}$.
Substituting this upper bound back into \eqref{eq:integral-substituted}, we have
\begin{align*}
  u^{\varepsilon}(\mathbf{x}, t) \leqslant \frac{\varepsilon^{N/2} M |S|}{(4\pi t)^{N/2}} \exp\left(\frac{r_{2} t}{\varepsilon}-\frac{\max(0, |\mathbf{x}|-\varepsilon R_{S})^{2}}{4\varepsilon t} \right).
\end{align*}
Applying the Hopf-Cole transformation $w^{\varepsilon}=-\varepsilon \ln u^{\varepsilon}$, we conclude that
\begin{align}
\label{eq:optimal-lower-bound}
  w^{\varepsilon}(\mathbf{x}, t) \geqslant -r_{2} t+\frac{\max(0, |\mathbf{x}|-\varepsilon R_{S})^{2}}{4t}-\varepsilon \ln \left(\frac{\varepsilon^{N/2} M |S|}{(4\pi t)^{N/2}} \right).
\end{align}
The right-hand side of \eqref{eq:optimal-lower-bound} converges locally uniformly in $\mathbb{R}^{N}\times\mathbb{R}^{+}$ to $\frac{|\mathbf{x}|^{2}}{4t}-r_{2}t$ as $\varepsilon\to0$.
Taking the lower half-relaxed limit in \eqref{eq:optimal-lower-bound} therefore yields \eqref{eq:optimal-sub-estimate}.
This completes the proof.
\end{proof}

\subsection{Uniqueness and self-similarity of the geometric optics limit}

We first record the corresponding comparison principle.

\begin{proposition}%
\label{prop:singular-comparison}
Let $U\in USC(\mathbb{R}^{N}\times\mathbb{R}^{+})$ be a viscosity subsolution of \eqref{eq:obstacle-direct}, and let $V\in LSC(\mathbb{R}^{N}\times\mathbb{R}^{+})$ be a viscosity supersolution of \eqref{eq:obstacle-direct}.
Assume that there exists $C>0$ such that
\begin{align}
\label{eq:singular-comparison-bounds}
  U(\mathbf{x},t) \leqslant \frac{|\mathbf{x}|^{2}}{4t}+Ct, \quad V(\mathbf{x},t) \geqslant \frac{|\mathbf{x}|^{2}}{4t}-Ct,
  \quad (\mathbf{x},t)\in\mathbb{R}^{N}\times\mathbb{R}^{+}.
\end{align}
Then
\begin{align*}
  U(\mathbf{x},t)\leqslant V(\mathbf{x},t) \quad \text{in }\mathbb{R}^{N}\times\mathbb{R}^{+}.
\end{align*}
Consequently, there is at most one nonnegative continuous viscosity solution $w$ of \eqref{eq:obstacle-direct} for which there exists $C_{w}>0$ such that
\begin{align}
\label{eq:growth-bound}
  \frac{|\mathbf{x}|^{2}}{4t}-C_{w}t \leqslant w(\mathbf{x},t) \leqslant \frac{|\mathbf{x}|^{2}}{4t}+C_{w}t, \quad (\mathbf{x},t)\in\mathbb{R}^{N}\times\mathbb{R}^{+}.
\end{align}
\end{proposition}

\begin{proof}
Fix an arbitrary $\tau\in\mathbb{R}^{+}$ and set $a_{\tau}=2\tau C$.
By \eqref{eq:singular-comparison-bounds},
\begin{align}
\label{eq:singular-order-at-tau}
  U(\mathbf{x},\tau) \leqslant V(\mathbf{x},\tau)+a_{\tau}, \quad \mathbf{x}\in\mathbb{R}^{N}.
\end{align}
Define $V^{\tau}(\mathbf{x},t)=V(\mathbf{x},t)+a_{\tau}$.
We first check that $V^{\tau}$ is still a viscosity supersolution of \eqref{eq:obstacle-direct} in $\mathbb{R}^{N}\times(\tau,\infty)$.
Let $\phi\in C^{1}$ touch $V^{\tau}$ from below at $(\mathbf{x}_{0},t_{0})$ with $t_{0}>\tau$.
Then $\phi-a_{\tau}$ touches $V$ from below at the same point.
Since $V$ is a viscosity supersolution,
\begin{align*}
  \min\left\{\partial_{t}\phi(\mathbf{x}_{0},t_{0})+|D\phi(\mathbf{x}_{0},t_{0})|^{2} +R^{*}(\mathbf{x}_{0},t_{0}), V(\mathbf{x}_{0},t_{0}) \right\}\geqslant0.
\end{align*}
Hence
\begin{align*}
  \partial_{t}\phi(\mathbf{x}_{0},t_{0})+|D\phi(\mathbf{x}_{0},t_{0})|^{2} +R^{*}(\mathbf{x}_{0},t_{0})\geqslant0
  \quad \text{and} \quad V^{\tau}(\mathbf{x}_{0},t_{0})\geqslant0,
\end{align*}
and therefore $V^{\tau}$ is a viscosity supersolution.

Fix an arbitrary $T>\tau$.
By Corollary \ref{coro:comparison-principle} and \eqref{eq:singular-order-at-tau}, we have
\begin{align*}
  U(\mathbf{x},t) \leqslant V^{\tau}(\mathbf{x},t)=V(\mathbf{x},t)+2C\tau, \quad (\mathbf{x},t)\in\mathbb{R}^{N}\times[\tau,T].
\end{align*}
Since $T>\tau$ is arbitrary, it follows that
\begin{align*}
  U(\mathbf{x},t) \leqslant V(\mathbf{x},t)+2C\tau, \quad (\mathbf{x},t)\in\mathbb{R}^{N}\times[\tau,\infty).
\end{align*}
For fixed $(\mathbf{x},t)\in\mathbb{R}^{N}\times\mathbb{R}^{+}$, let $\tau\downarrow0$ with $\tau<t$. We obtain
\begin{align*}
  U(\mathbf{x},t)\leqslant V(\mathbf{x},t), \quad (\mathbf{x},t)\in\mathbb{R}^{N}\times\mathbb{R}^{+}.
\end{align*}

Finally, let $w_{1}$ and $w_{2}$ be two nonnegative continuous viscosity solutions satisfying \eqref{eq:growth-bound} with constants $C_{1}$ and $C_{2}$, respectively.
Set $ C=\max\{C_{1},C_{2}\}$.
Applying the first part of the proposition with $U=w_{1}$ and $V=w_{2}$ gives $w_{1}\leqslant w_{2}$.
Interchanging $w_{1}$ and $w_{2}$, we have $w_{2}\leqslant w_{1}$. Hence $w_{1}=w_{2}$. The proof is complete.
\end{proof}

We next apply this comparison result to prove Theorem \ref{thm:wkb-direct-convergence}.

\begin{proof}[Proof of Theorem \ref{thm:wkb-direct-convergence}]
\textbf{Step 1.} We show the viscosity supersolution and subsolution  properties.

By the comparison principle, $u^{\varepsilon} \leqslant M$ with $M=\max\{r_{2}, \|g\|_{L^{\infty}}\}$, which implies $\overline w\geqslant \underline{w} \geqslant 0$.
For any test function $\phi\in C^{2}$, suppose that $\underline w-\phi$ attains its strict local minimum at $(\mathbf{x}_{0},t_{0})\in \mathbb{R}^{N}\times\mathbb{R}^{+}$.
By the definition of half-relaxed limits, up to extracting a subsequence, there exists a sequence $\{(\mathbf{x}_{\varepsilon},t_{\varepsilon})\}_{\varepsilon}$ such that $w^{\varepsilon}-\phi$ attains its local minimum at $(\mathbf{x}_{\varepsilon},t_{\varepsilon})$ and
\begin{align*}
  w^{\varepsilon}(\mathbf{x}_{\varepsilon},t_{\varepsilon}) \rightarrow \underline w(\mathbf{x}_{0},t_{0})
  \quad \text{and} \quad (\mathbf{x}_{\varepsilon},t_{\varepsilon})\rightarrow (\mathbf{x}_{0},t_{0})
\end{align*}
as $\varepsilon\rightarrow0$.
At the point $(\mathbf{x}_{\varepsilon},t_{\varepsilon})$, the maximum principle implies that
\begin{align}
\label{eq:p-sup-test}
  \partial_{t}(w^{\varepsilon}-\phi)=0, \quad D(w^{\varepsilon}-\phi)=0, \quad \Delta(w^{\varepsilon}-\phi) \geqslant 0.
\end{align}
From \eqref{eq:w-eps-pde-main}, by the nonnegativity of $u^{\varepsilon}$ and using \eqref{eq:p-sup-test}, we obtain
\begin{align*}
  \partial_{t}\phi(\mathbf{x}_{\varepsilon},t_{\varepsilon}) +|D\phi(\mathbf{x}_{\varepsilon},t_{\varepsilon})|^{2}+r\left(\frac{\mathbf{x}_{\varepsilon}\cdot \mathbf{e}-ct_{\varepsilon}}{\varepsilon}\right)
  \geqslant \varepsilon \Delta \phi(\mathbf{x}_{\varepsilon},t_{\varepsilon}).
\end{align*}
Since the limits and bounds of $r$ yield
\begin{align*}
  \limsup_{\varepsilon\to0}
  r\left(\frac{\mathbf{x}_{\varepsilon}\cdot\mathbf{e}-ct_{\varepsilon}}{\varepsilon}\right)\leqslant R^{*}(\mathbf{x}_{0},t_{0}),
\end{align*}
we then pass to the limsup as $\varepsilon\rightarrow 0$ to conclude that
\begin{align*}
  \partial_{t}\phi(\mathbf{x}_{0},t_{0})+|D\phi(\mathbf{x}_{0},t_{0})|^{2}+R^{*}(\mathbf{x}_{0},t_{0}) \geqslant 0.
\end{align*}
Thus, $\underline{w}$ is a viscosity supersolution of the variational inequality \eqref{eq:obstacle-direct}.

For the subsolution property, for any $\phi\in C^{2}$, suppose that $\overline w-\phi$ attains its strict local maximum at $(\mathbf{x}_{0},t_{0})\in\{\overline w>0\}$.
The upper half-relaxed limit argument gives, along a subsequence, local maximum points $(\mathbf{x}_{\varepsilon},t_{\varepsilon})$ of $w^{\varepsilon}-\phi$ such that
\begin{align*}
  w^{\varepsilon}(\mathbf{x}_{\varepsilon},t_{\varepsilon})\to\overline w(\mathbf{x}_{0},t_{0})>0
  \quad \text{and} \quad
   (\mathbf{x}_{\varepsilon},t_{\varepsilon}) \to(\mathbf{x}_{0},t_{0})
\end{align*}
as $\varepsilon\rightarrow0$.
At these points, it follows from \eqref{eq:w-eps-pde-main} that
\begin{align*}
  \partial_t\phi(\mathbf{x}_{\varepsilon},t_{\varepsilon})+|D\phi(\mathbf{x}_{\varepsilon},t_{\varepsilon})|^{2}+r\left(\frac{\mathbf{x}_{\varepsilon}\cdot\mathbf{e}
  -ct_{\varepsilon}}{\varepsilon}\right)
  \leqslant
  u^{\varepsilon}(\mathbf{x}_{\varepsilon},t_{\varepsilon})
  +\varepsilon\Delta\phi(\mathbf{x}_{\varepsilon},t_{\varepsilon}).
\end{align*}
The test condition implies $u^{\varepsilon}(\mathbf{x}_{\varepsilon},t_{\varepsilon})\to0$ in this case.
Meanwhile, we have
\begin{align*}
  \liminf_{\varepsilon\to0}
  r\left(\frac{\mathbf{x}_{\varepsilon}\cdot\mathbf{e}
  -ct_{\varepsilon}}{\varepsilon}\right)
  \geqslant R_{*}(\mathbf{x}_{0},t_{0}).
\end{align*}
Taking the liminf therefore gives
\begin{align*}
  \partial_t\phi(\mathbf{x}_{0},t_{0})
  +|D\phi(\mathbf{x}_{0},t_{0})|^{2}
  +R_{*}(\mathbf{x}_{0},t_{0})\leqslant0.
\end{align*}
Thus $\overline w$ is a viscosity subsolution of \eqref{eq:obstacle-direct}.

\textbf{Step 2.} We verify the growth control condition and then apply Proposition \ref{prop:singular-comparison}.

Using \eqref{eq:upper-bound-wep-optimal} derived from the fundamental solution, for all $t\in\mathbb{R}^{+}$, we have
\begin{align}
\label{eq:upper-w-tau}
  \overline{w}(\mathbf{x},t) \leqslant \frac{|\mathbf{x}|^{2}}{4t}+C_{0}t
\end{align}
for some $C_{0}\geqslant 0$.
Meanwhile, Lemma \ref{lem:initial-lower-bound} implies
\begin{align}
\label{eq:lower-w-tau}
  \underline{w}(\mathbf{x},t) \geqslant \frac{|\mathbf{x}|^{2}}{4t}-r_{2}t.
\end{align}
Combining these bounds with Step 1 and taking $C=\max\{C_{0},r_{2}\}$, we see that $\overline w$ and $\underline w$ satisfy the assumptions of Proposition \ref{prop:singular-comparison}. Then it follows that
\begin{align*}
  \overline w(\mathbf{x},t)\leqslant \underline w(\mathbf{x},t) \quad \text{in }\mathbb{R}^{N}\times\mathbb{R}^{+}.
\end{align*}
The reverse inequality $\underline w\leqslant\overline w$ holds by the definition of the half-relaxed limits. Therefore
\begin{align*}
  \overline w(\mathbf{x},t)=\underline w(\mathbf{x},t)=:w(\mathbf{x},t) \quad \text{in }\mathbb{R}^{N}\times\mathbb{R}^{+}.
\end{align*}
Since $\overline w$ is u.s.c. and $\underline w$ is l.s.c., $w$ is continuous.
In particular, $w^{\varepsilon}\rightarrow w$ locally uniformly in $\mathbb{R}^{N}\times\mathbb{R}^{+}$.
Moreover, $w$ is both a viscosity subsolution and a viscosity supersolution of \eqref{eq:obstacle-direct}.
The bounds \eqref{eq:upper-w-tau} and \eqref{eq:lower-w-tau} also show that $w$ satisfies \eqref{eq:singular-asymptotic-precise}.
Finally, the uniqueness assertion follows directly from Proposition \ref{prop:singular-comparison}.

To show the singular initial condition, taking $\mathbf{x}=\mathbf{0}$ in \eqref{eq:singular-asymptotic-precise} and recalling that $w\geqslant0$, we obtain $0\leqslant  w(\mathbf{0},t)\leqslant Ct$, and hence $w(\mathbf{0},t)\to0$ as $t\downarrow0$.
If $\mathbf{x}\neq\mathbf{0}$, then
\begin{align*}
  w(\mathbf{x},t) \geqslant \frac{|\mathbf{x}|^{2}}{4t}-Ct \rightarrow+\infty \quad\text{as }t\downarrow0.
\end{align*}
Therefore $w$ satisfies the singular initial condition stated in the theorem. We complete the proof.
\end{proof}

\begin{proof}[Proof of Proposition \ref{prop:self-similarity}]
Fix $k>0$ and define
\begin{align*}
  w_{k}(\mathbf{x},t)=\frac{1}{k} w(k\mathbf{x},kt).
\end{align*}
The shifting interface $\{ (\mathbf{x},t) \in \mathbb{R}^{N} \times \mathbb{R}^{+} \mid \mathbf{x}\cdot \mathbf{e}=ct \}$ is invariant under $(\mathbf{x},t)\mapsto(k\mathbf{x},kt)$, and hence $R(k\mathbf{x},kt)=R(\mathbf{x},t)$; the same identity holds for the semicontinuous envelopes of $R$.
Pulling back test functions by this scaling shows that $w_{k}$ is again a nonnegative continuous viscosity solution of \eqref{eq:obstacle-direct}.
Moreover, \eqref{eq:singular-asymptotic-precise} is invariant under the same scaling, so $w_{k}$ belongs to the uniqueness class in Theorem \ref{thm:wkb-direct-convergence}.
Therefore $w_{k}=w$, which is equivalent to \eqref{eq:self-similarity}.
\end{proof}

\section{A zero-sum stopping game representation for the variational inequality}
\label{section:stopping-game}

This section represents the solution of the limiting variational inequality as the value of a zero-sum stopping game.
We first establish the basic properties of the stopping value.
We then identify it as the unique viscosity solution through bounded-data approximation.
Finally, we derive a Wulff shape spreading result from the geometry of its zero set.

\subsection{The stopping game and basic properties of its value}

\begin{definition}
For fixed $(\mathbf{x},t)\in \mathbb{R}^{N}\times\mathbb{R}^{+}$, $\vartheta\in\Theta$, and $\epsilon\in(0,1)$, we say a trajectory $\gamma_{\epsilon}\in\mathbb{X}^{\mathbf{x}, t}$ is \textit{$\epsilon$-optimal} for the inner minimization of \eqref{eq:stopping-game-value} if
\begin{align*}
  \mathcal P(\mathbf{x},t;\vartheta,\gamma_{\epsilon}) \leqslant \inf_{\gamma\in\mathbb{X}^{\mathbf{x},t}} \mathcal P(\mathbf{x},t;\vartheta,\gamma)+\epsilon.
\end{align*}
\end{definition}

\begin{remark}\label{rem:epsilon-optimal-trajectories}
Whenever an $\epsilon$-optimal trajectory is used below, the corresponding inner infimum is a finite real number.
Indeed, the admissible class is nonempty, and a straight or piecewise linear admissible trajectory gives a finite upper bound.
Moreover, since $0\leqslant R\leqslant r_{2}$ and every terminal cost used below is nonnegative, each corresponding payoff is bounded below by $-r_{2}t$.
Therefore an $\epsilon$-optimal trajectory exists for every $\epsilon\in(0,1)$.
\end{remark}

\begin{proposition}
\label{prop:stopping-game-basic-properties} The stopping value function $W$ defined by \eqref{eq:stopping-game-value} is locally bounded, nonnegative, and continuous in $\mathbb{R}^{N}\times\mathbb{R}^{+}$.
Moreover, it satisfies \eqref{eq:singular-asymptotic-precise}.
\end{proposition}

\begin{proof}[Proof of Proposition \ref{prop:stopping-game-basic-properties}]
The nonnegativity is a direct consequence of the admissible stopping rule $\vartheta\equiv0$.
Indeed, for this rule the payoff is identically zero, and therefore
\begin{align}
\label{eq:W-nonnegative}
  W(\mathbf{x},t)\geqslant0.
\end{align}
On the other hand, fix an arbitrary stopping time $\vartheta\in\Theta$.
In the inner infimum, we may choose the zero-extended straight trajectory
\begin{align*}
  \gamma(s)=
  \begin{cases}
    (1-\frac{s}{t})\mathbf{x},&\text{if }s\in[0,t],\\
    \mathbf{0},&\text{if } s\in(t,\infty).
  \end{cases}
\end{align*}
Let $\tau_{\vartheta}=t\wedge\vartheta[\gamma]$.
Since $0\leqslant\tau_{\vartheta}\leqslant t$ and $R\geqslant r_{1}>0$, we have
\begin{align*}
  \inf_{\eta\in\mathbb{X}^{\mathbf{x},t}} \mathcal P(\mathbf{x},t;\vartheta,\eta)
  \leqslant \int_{0}^{\tau_{\vartheta}} \left[ \frac{1}{4}|\dot\gamma(s)|^{2}-R(\gamma(s),t-s) \right]ds
  \leqslant \tau_{\vartheta} \left(\frac{|\mathbf{x}|^{2}}{4t^{2}}-r_{1} \right)\leqslant \frac{|\mathbf{x}|^{2}}{4t}.
\end{align*}
Since $\vartheta\in\Theta$ is arbitrary, we conclude that
\begin{align}
\label{eq:stopping-value-upper-bound}
  W(\mathbf{x},t)\leqslant\frac{|\mathbf{x}|^{2}}{4t}.
\end{align}
This, together with \eqref{eq:W-nonnegative}, implies that $W$ is locally bounded in $\mathbb{R}^{N}\times\mathbb{R}^{+}$.

We next prove continuity away from $t=0$.
It suffices to work on a compact cylinder $K=\overline{B_{M}(\mathbf{0})}\times[\tau,T]$, with $M>0$ and $0<\tau<T<\infty$.
Let $(\mathbf{x},t),(\mathbf{y},\sigma)\in K$, and set $d=|\mathbf{x}-\mathbf{y}|+|t-\sigma|$.
We also fix
\begin{align}
\label{eq:kappa-choice}
  A_{K}=\frac{4}{\tau}, \quad B_{K}=|c|\left(1+\frac{6T}{\tau}\right), \quad \kappa=2B_{K}+1.
\end{align}
It suffices to prove the existence of a modulus of continuity $\omega_{K}$ with $\omega_{K}(d)\to0$ as $d\to0$ such that
\begin{align}
\label{eq:stopping-value-modulus}
  W(\mathbf{y}, \sigma) \leqslant W(\mathbf{x}, t)+\omega_{K}(d),
\end{align}
since interchanging the two points gives the reverse inequality.

Set $d_{K}=\min\left\{1, \frac{\tau^{2}}{16}\right\}$.
For the continuity estimate, it suffices to consider $d\in(0, d_{K})$. Write $\delta=\sqrt d$.
Then $d\leqslant\delta<\frac{\tau}{4}$.
Let $\Lambda=\frac{t}{\sigma-2\delta}$, and define $\Pi:[0,\infty]\to[0,t]$ by
\begin{align}
\label{eq:pullback-time}
  \Pi(s)=
  \begin{cases}
    0, &\text{if }s\in[0,\delta],\\
    \Lambda(s-\delta), &\text{if }s\in(\delta,\sigma-\delta),\\
    t, &\text{if }s\in[\sigma-\delta,\infty].
  \end{cases}
\end{align}
On the middle interval, we denote $v=\Pi(s)=\Lambda(s-\delta)$.

Fix $\tilde\vartheta\in\Theta$.
For every $\gamma\in\mathbb{X}$, define $\eta=\mathcal T\gamma\in\mathbb{X}$ by
\begin{align}
\label{eq:traj-transform}
  (\mathcal T\gamma)(s)=
  \begin{cases}
    \mathbf{y}+(\gamma(0)-\mathbf{x}) +\frac{s}{\delta}(\mathbf{x}+\kappa\delta\mathbf{e}-\mathbf{y}), &\text{if }s\in[0,\delta],\\
    \gamma(v)+\kappa\delta\mathbf{e},&\text{if }s\in(\delta,\sigma-\delta),\\
    \gamma(t)+\kappa(\sigma-s)\mathbf{e}, &\text{if }s\in[\sigma-\delta,\sigma],\\
    \gamma(t), &\text{if }s\in(\sigma,\infty).
  \end{cases}
\end{align}
Indeed, each branch in \eqref{eq:traj-transform} is locally $H^{1}$, and they agree at the three junction points $s=\delta$, $s=\sigma-\delta$, and $s=\sigma$.
Hence, $\mathcal T\gamma\in\mathbb{X}$.
In particular, if $\gamma\in\mathbb{X}^{\mathbf{x},t}$, then $\eta\in\mathbb{X}^{\mathbf{y},\sigma}$ is reduced to:
\begin{align}
\label{eq:traj-eta}
  \eta(s)=
  \begin{cases}
    \mathbf{y}+\frac{s}{\delta}(\mathbf{x}+\kappa\delta \mathbf{e}-\mathbf{y}), &\text{if } s\in[0,\delta], \\
    \gamma(v)+\kappa\delta\mathbf{e}, &\text{if } s\in(\delta,\sigma-\delta),\\
    \kappa(\sigma-s)\mathbf{e}, &\text{if } s\in[\sigma-\delta,\sigma],\\
    \mathbf{0}, &\text{if }s\in(\sigma,\infty).
  \end{cases}
\end{align}

Define
\begin{align*}
  \vartheta[\gamma]=\Pi(\tilde\vartheta[\eta]) \quad  \text{for}  \quad \gamma\in\mathbb{X}.
\end{align*}
We verify that $\vartheta$ is a stopping time on $\mathbb{X}$.
Let $\gamma_{1},\gamma_{2}\in\mathbb{X}$ agree on $[0,a]$ for some $a\geqslant0$, and let $\vartheta[\gamma_{1}]\leqslant a$.
Denote by $\eta_{i}=\mathcal{T}(\gamma_{i})$ the trajectory constructed from $\gamma_{i}$, $i=1,2$, and set $s_{1}=\tilde\vartheta[\eta_{1}]$.
Then $\Pi(s_{1})=\vartheta[\gamma_{1}]\leqslant a$.
If $s_{1}\geqslant\sigma-\delta$, then $\Pi(s_{1})=t\leqslant a$.
Hence $\gamma_{1}=\gamma_{2}$ on $[0,t]$, and therefore $\eta_{1}=\eta_{2}$ on $[0,\infty)$.
In particular, $\tilde\vartheta[\eta_{1}]=\tilde\vartheta[\eta_{2}]$.
Otherwise, $s_{1}<\sigma-\delta$.
On the first connecting interval, $\eta_{i}$ depends on $\gamma_{i}$ only through $\gamma_{i}(0)$.
Since $\gamma_{1}(0)=\gamma_{2}(0)$, we have
\begin{align*}
  \eta_{1}(s)=\eta_{2}(s) \quad\text{for }s\in[0,\min\{\delta,s_{1}\}].
\end{align*}
If $s_{1}\leqslant\delta$, the preceding equality already gives $\eta_{1}=\eta_{2}$ on $[0,s_{1}]$.
If $s_{1}>\delta$, it remains to consider $s\in(\delta,s_{1}]$.
Since $\Pi$ is nondecreasing, $\Pi(s)\leqslant\Pi(s_{1})\leqslant a$.
Thus, $\gamma_{1}(\Pi(s))=\gamma_{2}(\Pi(s))$, and consequently $\eta_{1}(s)=\eta_{2}(s)$ for $s\in(\delta,s_{1}]$. Therefore, in either case,
\begin{align*}
  \eta_{1}(s)=\eta_{2}(s) \quad\text{for all }s\in[0,s_{1}].
\end{align*}
In the case $s_1<\sigma-\delta$,  we have $s_1=\tilde\vartheta[\eta_1]<\infty$.
Since $\eta_{1}=\eta_{2}$ on $[0,s_{1}]$ and $\tilde\vartheta$ is a stopping time, it follows again that $\tilde\vartheta[\eta_{1}]=\tilde\vartheta[\eta_{2}]$.
Combining this with the former case, we conclude that $\vartheta[\gamma_{1}]=\vartheta[\gamma_{2}]$. Thus $\vartheta\in\Theta$.
Moreover, since $\Pi(\cdot)\in[0,t]$, we have $\vartheta[\gamma]\in [0,t]$ for every $\gamma\in\mathbb{X}$.

Since $(\mathbf{x},t),(\mathbf{y},\sigma)\in K$ and $\delta<\tau/4$, we have
\begin{align}
\label{eq:Lambda-bound}
  |\Lambda-1|=\frac{|t-\sigma+2\delta|}{\sigma-2\delta} \leqslant A_{K}(d+\delta).
\end{align}
Moreover, on the middle interval $s\in(\delta,\sigma-\delta)$, we have
\begin{align*}
  \left(\gamma(v)+\kappa\delta\mathbf{e}\right)\cdot\mathbf{e} -c(\sigma-s)-\left(\gamma(v)\cdot\mathbf{e}-c(t-v)\right)=\kappa\delta+c(t-v-\sigma+s)
  \geqslant \kappa\delta-B_{K}(d+\delta).
\end{align*}
Indeed, the last inequality follows from
\begin{align*}
  |t-v-\sigma+s| \leqslant |t-\sigma|+|\Lambda-1|s+\Lambda\delta \leqslant \left(1+\frac{6T}{\tau}\right)(d+\delta),
\end{align*}
where we used $\sigma\in[s, T]$ and $\Lambda\leqslant \frac{2T}{\tau}$.
Since $d\leqslant\delta$ by construction, the choice \eqref{eq:kappa-choice} yields
\begin{align*}
  \kappa\delta-B_{K}(d+\delta)\geqslant\delta>0.
\end{align*}
Note that $R(\mathbf{z},\tau)$ is nondecreasing with respect to $\mathbf{z}\cdot\mathbf{e}-c\tau$. The previous inequality implies
\begin{align}
\label{eq:middle-monotone}
  R\left(\gamma(v)+\kappa\delta\mathbf{e},\sigma-s\right) \geqslant R(\gamma(v),t-v) \quad\text{on the middle interval}.
\end{align}
By \eqref{eq:middle-monotone}, the term involving $-R$ along the transformed trajectory is not greater than the corresponding term along the original trajectory on the middle interval.

We now compare the full payoff along a near minimizer.
In the estimates below, $C_{K}$ denotes a positive constant depending only on $K,r_{2},c$ and $N$, whose value may change from line to line.
Let $\epsilon\in(0,1)$ and choose an $\epsilon$-optimal trajectory $\gamma_{\epsilon}$ for the inner minimization associated with the stopping time $\vartheta\in\Theta$.
The bound \eqref{eq:stopping-value-upper-bound} and the inequality $R\leqslant r_{2}$ imply
\begin{align}
\label{eq:energy-bound-for-near-minimizer}
  \int_{0}^{t\wedge\vartheta[\gamma_{\epsilon}]}|\dot\gamma_{\epsilon}(s)|^{2}\,ds \leqslant C_{K}.
\end{align}
Indeed,
\begin{align*}
  \frac{1}{4} \int_{0}^{t\wedge\vartheta[\gamma_{\epsilon}]} |\dot\gamma_{\epsilon}(s)|^{2}\,ds
  &= \mathcal P(\mathbf{x},t;\vartheta,\gamma_{\epsilon})+\int_{0}^{t\wedge\vartheta[\gamma_{\epsilon}]} R(\gamma_{\epsilon}(s),t-s)\,ds
  \leqslant W(\mathbf{x},t)+\epsilon+r_{2}T.
\end{align*}
Let $\eta_{\epsilon}\in\mathbb{X}^{\mathbf{y},\sigma}$ be the trajectory obtained from $\gamma_{\epsilon}$ by the construction \eqref{eq:traj-transform} (note that it is reduced to \eqref{eq:traj-eta}).
On the first connecting interval,
\begin{align*}
  \frac{1}{4}\int_{0}^{\delta} \left|\frac{\mathbf{x}+\kappa\delta\mathbf{e}-\mathbf{y}}{\delta}\right|^{2}ds \leqslant C_{K}\frac{d^{2}+\delta^{2}}{\delta}.
\end{align*}
On the branch $s\in[\sigma-\delta,\sigma]$, the integral of the quadratic term is bounded by $C_{K}\delta$.
The trajectory is constant after $s=\sigma$, which does not contribute to the payoff.
Since the quadratic integrand is nonnegative, the same bounds remain valid when either connecting interval is truncated by the stopping time. Let
\begin{align*}
  \tau_{\epsilon}=\tilde\vartheta[\eta_{\epsilon}], \qquad \widehat\tau_{\epsilon}=\sigma\wedge\tau_{\epsilon}.
\end{align*}
Thus, the payoff $\mathcal P(\mathbf{y},\sigma;\tilde\vartheta,\eta_{\epsilon})$ is integrated over $[0,\widehat\tau_{\epsilon}]$.
Define $a_{\epsilon}=\left(\widehat\tau_{\epsilon}\wedge(\sigma-\delta)\right)\vee\delta$.
If $\widehat\tau_{\epsilon}\leqslant\delta$, then $a_{\epsilon}=\delta$.
If $\delta<\widehat\tau_{\epsilon}<\sigma-\delta$, then $a_{\epsilon}=\widehat\tau_{\epsilon}$.
Finally, if $\widehat\tau_{\epsilon}\geqslant\sigma-\delta$, then $a_{\epsilon}=\sigma-\delta$.
Thus, up to a set of measure zero, $[\delta,a_{\epsilon}]$ coincides with $[0,\widehat\tau_{\epsilon}] \cap [\delta,\sigma-\delta]$.
Under the change of variables $v=\Pi(s)$, the interval $[\delta,a_{\epsilon}]$ is mapped to $[0,\Pi(a_{\epsilon})]$.
In particular, its upper endpoint satisfies
\begin{align}
\label{eq:middle-endpoint-identity}
  \Pi(a_{\epsilon}) =
  \begin{cases}
    0, &\text{if } \widehat\tau_{\epsilon}\in[0,\delta],\\
    \Pi(\widehat\tau_{\epsilon}), &\text{if } \widehat\tau_{\epsilon}\in(\delta,\sigma-\delta),\\
    t, &\text{if } \widehat\tau_{\epsilon}\geqslant\sigma-\delta,
  \end{cases}
  \quad  \Longrightarrow \Pi(a_{\epsilon})=\Pi(\widehat\tau_{\epsilon}) =\Pi(\tau_{\epsilon})=\vartheta[\gamma_{\epsilon}].
\end{align}
Here $\Pi(\widehat\tau_{\epsilon})=\Pi(\tau_{\epsilon})$ because either $\tau_{\epsilon}\leqslant\sigma$, in which case the two arguments coincide, or $\tau_{\epsilon}>\sigma$, in which case both arguments belong to the region where $\Pi\equiv t$.
On $[\delta,a_{\epsilon}]$, we have $v=\Pi(s)=\Lambda(s-\delta)$, and hence \eqref{eq:middle-monotone} implies
\begin{align*}
  -\int_{\delta}^{a_{\epsilon}} R(\eta_{\epsilon}(s),\sigma-s)\,ds
  =& -\frac{1}{\Lambda} \int_{0}^{\vartheta[\gamma_{\epsilon}]} R\left(\gamma_{\epsilon}(v)+\kappa\delta\mathbf{e}, \sigma-\delta-\frac{v}{\Lambda} \right)dv\\
  \leqslant& -\frac{1}{\Lambda} \int_{0}^{\vartheta[\gamma_{\epsilon}]} R(\gamma_{\epsilon}(v),t-v)\,dv\\
  \leqslant& -\int_{0}^{\vartheta[\gamma_{\epsilon}]} R(\gamma_{\epsilon}(v),t-v)\,dv+C_{K}(d+\delta),
\end{align*}
where the last inequality follows from $0< R\leqslant r_{2}$ and \eqref{eq:Lambda-bound} and the fact that $\Lambda$ is bounded away from zero on $K$.
On this same part, since $\dot\eta_{\epsilon}(s)=\Lambda\dot\gamma_{\epsilon}(v)$, it holds that
\begin{align*}
  \frac{1}{4} \int_{\delta}^{a_{\epsilon}} |\dot\eta_{\epsilon}(s)|^{2}ds
  \leqslant \frac{\Lambda}{4}\int_{0}^{\vartheta[\gamma_{\epsilon}]} |\dot\gamma_{\epsilon}(v)|^{2}dv.
\end{align*}
Using \eqref{eq:Lambda-bound}, \eqref{eq:energy-bound-for-near-minimizer}, and \eqref{eq:middle-endpoint-identity}, we obtain
\begin{align*}
  \frac{1}{4} \int_{\delta}^{a_{\epsilon}} |\dot\eta_{\epsilon}(s)|^{2}\,ds
  \leqslant \frac{1}{4} \int_{0}^{\vartheta[\gamma_{\epsilon}]} |\dot\gamma_{\epsilon}(v)|^{2}\,dv+C_{K}(d+\delta).
\end{align*}
Since $R\geqslant0$, omitting the term $-R$ on the two straight connecting branches gives an upper bound.
Moreover, by \eqref{eq:middle-endpoint-identity} and the fact that $\Pi$ takes values in $[0,t]$, we have $\vartheta[\gamma_{\epsilon}]\leqslant t$. Hence,
\begin{align*}
  \int_{0}^{\vartheta[\gamma_{\epsilon}]} L(\dot\gamma_{\epsilon}(v),\gamma_{\epsilon}(v),t-v)\,dv=\mathcal P(\mathbf{x},t;\vartheta,\gamma_{\epsilon}).
\end{align*}
Consequently, the estimates for two straight branches and the middle-interval yield
\begin{align}
\label{eq:path-perturbation-bound}
  \mathcal P(\mathbf{y},\sigma;\tilde\vartheta,\eta_{\epsilon})
  \leqslant \mathcal P(\mathbf{x},t;\vartheta,\gamma_{\epsilon})+C_{K}\left(d+\delta+\frac{d^{2}+\delta^{2}}{\delta} \right).
\end{align}
For a fixed $\tilde\vartheta\in\Theta$, the trajectory $\eta_{\epsilon}$ is admissible. Hence
\begin{align*}
  \inf_{\eta\in\mathbb{X}^{\mathbf{y},\sigma}} \mathcal P(\mathbf{y},\sigma;\tilde\vartheta,\eta) \leqslant \mathcal P(\mathbf{y},\sigma;\tilde\vartheta,\eta_{\epsilon}).
\end{align*}
Combining this with \eqref{eq:path-perturbation-bound} and the $\epsilon$-optimality of $\gamma_{\epsilon}$ yields
\begin{align*}
  \inf_{\eta\in\mathbb{X}^{\mathbf{y},\sigma}} \mathcal P(\mathbf{y},\sigma;\tilde\vartheta,\eta)
  \leqslant \inf_{\gamma\in\mathbb{X}^{\mathbf{x},t}} \mathcal P(\mathbf{x},t;\vartheta,\gamma)+C_{K}\left(d+\delta+\frac{d^{2}+\delta^{2}}{\delta}\right)+\epsilon.
\end{align*}
Since $\vartheta\in\Theta$, the last infimum is bounded above by $W(\mathbf{x},t)$.
Taking the supremum over $\tilde\vartheta\in\Theta$, we have
\begin{align*}
  W(\mathbf{y},\sigma) \leqslant W(\mathbf{x},t)+C_{K}\left(d+\delta+\frac{d^{2}+\delta^{2}}{\delta}\right)+\epsilon.
\end{align*}
Letting $\epsilon\downarrow0$ and using $\delta=\sqrt d$, we conclude
\begin{align*}
  W(\mathbf{y},\sigma) \leqslant W(\mathbf{x},t)+C_{K}\sqrt d \quad \text{for all }d\ll1.
\end{align*}
Then \eqref{eq:stopping-value-modulus} follows, and $W$ is continuous on $\mathbb{R}^{N}\times\mathbb{R}^ {+} $.

It remains to verify the growth control condition.
Choose the stopping time $\vartheta\equiv t$.
Then every $\gamma\in\mathbb{X}^{\mathbf{x},t}$ connects $\mathbf{x}$ to $\mathbf{0}$ in time $t$, and therefore, by the Cauchy-Schwarz inequality,
\begin{align*}
  |\mathbf{x}|^{2}=|\gamma(t)-\gamma(0)|^{2}= \left|\int_{0}^{t}\dot\gamma(s)\,ds\right|^{2} \leqslant t \left(\int_{0}^{t}|\dot\gamma(s)|^{2}\,ds\right).
\end{align*}
Since $R\leqslant r_{2}$, we have
\begin{align*}
  W(\mathbf{x},t) \geqslant \inf_{\gamma\in\mathbb{X}^{\mathbf{x},t}} \int_{0}^{t} \left[ \frac{1}{4}|\dot\gamma(s)|^{2}-r_{2} \right]ds
  \geqslant \frac{|\mathbf{x}|^{2}}{4t}-r_{2}t.
\end{align*}
Together with \eqref{eq:stopping-value-upper-bound}, this proves \eqref{eq:singular-asymptotic-precise} with $C=r_{2}$. The proof is complete.
\end{proof}

\subsection{Terminal-cost approximation and viscosity characterization}

We now identify the singular stopping value as the unique viscosity solution of the limiting variational inequality by approximating the hard terminal constraint with bounded terminal costs.

Let $w_{0}\in C_{b}(\mathbb{R}^{N})$ and $w_{0}\geqslant0$.
Consider the following variational inequality
\begin{align}
\label{eq:bounded-data-obstacle}
  \begin{cases}
    \min\left\{w_{t}+|Dw|^{2}+R(\mathbf{x},t),w\right\}=0, &(\mathbf{x},t)\in\mathbb{R}^{N}\times\mathbb{R}^{+},\\
    w(\mathbf{x},0)=w_{0}(\mathbf{x}), &\mathbf{x}\in\mathbb{R}^{N}.
  \end{cases}
\end{align}
For $\mathbf{x}\in\mathbb{R}^{N}$, set
\begin{align}
\label{eq:A-traj-set}
  \mathbb{A}^{\mathbf{x}}=\left\{\gamma\in\mathbb{X}\mid \gamma(0)=\mathbf{x}\right\}.
\end{align}
Define
\begin{align}
\label{eq:bounded-data-stopping-game}
  W_{0}(\mathbf{x},t)=\sup_{\vartheta\in\Theta} \inf_{\gamma\in\mathbb{A}^{\mathbf{x}}} \mathcal P_{0}(\mathbf{x},t;\vartheta,\gamma),
\end{align}
where
\begin{align*}
  \mathcal P_{0}(\mathbf{x},t;\vartheta,\gamma)=\int_{0}^{t\wedge\vartheta[\gamma]} L(\dot\gamma(s),\gamma(s),t-s)\,ds+\mathbbm{1}_{\{\vartheta[\gamma]\geqslant t\}}w_{0}(\gamma(t)).
\end{align*}

\begin{proposition}\label{prop:bounded-data-stopping-game}
The function $W_{0}$ defined by \eqref{eq:bounded-data-stopping-game} is the unique bounded continuous viscosity solution of \eqref{eq:bounded-data-obstacle}.
\end{proposition}

We first state the two auxiliary results that will be used in its proof.

\begin{lemma}
\label{lem:bounded-data-value-regularity}
The value function $W_{0}$ is continuous on $\mathbb{R}^{N}\times\mathbb{R}^{+}$ and satisfies
\begin{align*}
  0\leqslant W_{0}(\mathbf{x},t) \leqslant w_{0}(\mathbf{x}) \leqslant \|w_{0}\|_{L^{\infty}} \quad \text{for all } (\mathbf{x},t)\in\mathbb{R}^{N}\times\mathbb{R}^{+}.
\end{align*}
\end{lemma}

\begin{proposition}[Dynamic Programming Principle]
\label{prop:bounded-data-dpp}
For every $(\mathbf{x},t)\in\mathbb{R}^{N}\times\mathbb{R}^{+}$ and $h\in(0,t)$, $W_{0}$ satisfies
\begin{align}
\label{eq:bounded-data-dpp}
  W_{0}(\mathbf{x},t)=\sup_{\vartheta\in\Theta} \inf_{\gamma\in\mathbb{A}^{\mathbf{x}}} \left\{\int_{0}^{h\wedge\vartheta[\gamma]} L(\dot\gamma(s),\gamma(s),t-s)\,ds+\mathbbm{1}_{\{\vartheta[\gamma]>h\}} W_{0}(\gamma(h),t-h) \right\}.
\end{align}
\end{proposition}

Detailed proofs of the preceding three results are collected in Appendix \ref{app:bounded-data}.

We next use bounded terminal costs to approximate the fixed endpoint constraint $\gamma(t)=\mathbf{0}$ and obtain the viscosity characterization of the original stopping value.

\begin{proof}[Proof of Theorem \ref{thm:stopping-game-singular-obstacle-solution}]
We argue by approximation from bounded continuous initial data.

Let $\chi\in C([0,\infty))$ be nondecreasing, with $\chi=0$ on $[0,1]$ and $\chi=1$ on $[2,\infty)$.
For $\varrho>0$ and $M>0$, set $g_{\varrho,M}(\mathbf{z})=M\chi(\frac{|\mathbf{z}|}{\varrho})$.
Then $g_{\varrho,M}=0$ on $B_{\varrho}(\mathbf{0})$, $g_{\varrho,M}=M$ outside $B_{2\varrho}(\mathbf{0})$, and for each fixed $\mathbf{z}\ne\mathbf{0}$, $g_{\varrho,M}(\mathbf{z})\to M$ as $\varrho\downarrow0$, while $g_{\varrho,M}(\mathbf{0})=0$. Define
\begin{align*}
  W_{\varrho,M}(\mathbf{x},t)=\sup_{\vartheta\in\Theta} \inf_{\gamma\in\mathbb{A}^{\mathbf{x}}} \left\{\int_{0}^{t\wedge\vartheta[\gamma]} L(\dot\gamma(s),\gamma(s),t-s)\,ds+\mathbbm{1}_{\{\vartheta[\gamma]\geqslant t\}} g_{\varrho,M}(\gamma(t)) \right\},
\end{align*}
where $\mathbb{A}^{\mathbf{x}}$ is defined in \eqref{eq:A-traj-set}.
By Proposition \ref{prop:bounded-data-stopping-game}, $W_{\varrho,M}$ is the unique bounded continuous viscosity solution of
\begin{align*}
  \begin{cases}
    \min\{w_{t}+|Dw|^{2}+R(\mathbf{x},t),w\}=0, &(\mathbf{x},t)\in\mathbb{R}^{N}\times\mathbb{R}^{+},\\
    w(\mathbf{x},0)=g_{\varrho,M}(\mathbf{x}), &\mathbf{x}\in\mathbb{R}^{N}.
  \end{cases}
\end{align*}
We complete the proof in the following two steps.

\textbf{Step 1.} We construct a joint approximation which converges locally uniformly to $W$.
Since $\mathbb{X}^{\mathbf{x},t}\subset\mathbb{A}^{\mathbf{x}}$ and $g_{\varrho,M}(\gamma(t))=g_{\varrho,M}(\mathbf{0})=0$ for all $\gamma\in\mathbb{X}^{\mathbf{x},t}$, we have, for each $\vartheta\in\Theta$,
\begin{align*}
  \inf_{\gamma\in\mathbb{A}^{\mathbf{x}}} \left\{\int_{0}^{t\wedge\vartheta[\gamma]} L(\dot\gamma(s),\gamma(s),t-s)\,ds+\mathbbm{1}_{\{\vartheta[\gamma]\geqslant t\}} g_{\varrho,M}(\gamma(t)) \right\}
  \leqslant \inf_{\gamma\in\mathbb{X}^{\mathbf{x},t}} \int_{0}^{t\wedge\vartheta[\gamma]} L(\dot\gamma(s),\gamma(s),t-s)\,ds .
\end{align*}
Taking the supremum over $\vartheta\in\Theta$ implies
\begin{align}
\label{eq:bounded-value-below-hard-value}
  W_{\varrho,M}(\mathbf{x},t)\leqslant W(\mathbf{x},t).
\end{align}

Conversely, fix a compact cylinder $K=\overline{B_{A}(\mathbf{0})}\times[\tau,T]$ with $0<\tau<T$.
Assume $0<\varrho<\min\left\{1,\tau^{2}/4\right\}$ and let $\ell_{\varrho}=\sqrt\varrho$, so that $0<\ell_{\varrho}<\tau/2$.
We prove the reverse inequality uniformly on $K$.
Fix $(\mathbf{x},t)\in K$, $\vartheta\in\Theta$, and $\epsilon\in(0,1)$.
Use the following zero-extended straight trajectory:
\begin{align*}
  \bar\gamma(s) =
  \begin{cases}
    \left(1-\frac{s}{t}\right)\mathbf{x}, &s\in[0,t],\\
    \mathbf{0}, &s\in(t,\infty).
  \end{cases}
\end{align*}
Then $\bar\gamma\in\mathbb{X}^{\mathbf{x},t}$ and $g_{\varrho,M}(\bar\gamma(t))=0$.
Since $(\mathbf{x},t)\in K$ and $R\geqslant 0$,
\begin{align*}
  &\int_{0}^{t\wedge\vartheta[\bar\gamma]} L(\dot{\bar\gamma}(s),\bar\gamma(s),t-s)\,ds+\mathbbm{1}_{\{\vartheta[\bar\gamma]\geqslant t\}} g_{\varrho,M}(\bar\gamma(t))
  \leqslant \int_{0}^{t}\frac{|\mathbf{x}|^{2}}{4t^{2}}\,ds \leqslant C_{K,0}:=\frac{A^{2}}{4\tau}.
\end{align*}
Then choose an $\epsilon$-optimal trajectory $\gamma_{\epsilon}\in\mathbb{A}^{\mathbf{x}}$ satisfying
\begin{equation}
\label{eq:bounded-near-minimizer-bound}
  \begin{aligned}
    &\int_{0}^{t\wedge\vartheta[\gamma_{\epsilon}]} L(\dot\gamma_{\epsilon}(s),\gamma_{\epsilon}(s),t-s)\,ds+\mathbbm{1}_{\{\vartheta[\gamma_{\epsilon}]\geqslant t\}} g_{\varrho,M}(\gamma_{\epsilon}(t)) \\
    \leqslant& \inf_{\gamma\in\mathbb A^{\mathbf{x}}} \left\{\int_{0}^{t\wedge\vartheta[\gamma]} L(\dot\gamma(s),\gamma(s),t-s)\,ds+\mathbbm{1}_{\{\vartheta[\gamma]\geqslant t\}} g_{\varrho,M}(\gamma(t)) \right\} +\epsilon
    \leqslant C_{K,0}+\epsilon.
  \end{aligned}
\end{equation}
Assume $M>C_{K,0}+r_{2}T+1$.
To construct a trajectory in $\mathbb{X}^{\mathbf{x},t}$ from $\gamma_{\epsilon}\in\mathbb A^{\mathbf{x}}$, we consider whether $\vartheta[\gamma_{\epsilon}]<t$.
If so, the trajectory may be modified after the stopping time without changing the payoff; otherwise, the terminal penalty allows us to modify the trajectory near $t$ with only a small additional cost.

\textbf{Case 1.} Suppose $\vartheta[\gamma_{\epsilon}]<t$, set $\tau_{\epsilon}=\vartheta[\gamma_{\epsilon}]$, and define
\begin{align*}
  \eta_{\epsilon}(s) =
  \begin{cases}
    \gamma_{\epsilon}(s), &s\in[0,\tau_{\epsilon}),\\
    \frac{t-s}{t-\tau_{\epsilon}}\gamma_{\epsilon}(\tau_{\epsilon}), &s\in[\tau_{\epsilon},t),\\
    \mathbf{0}, &s\in[t,\infty).
  \end{cases}
\end{align*}
Then $\eta_{\epsilon}\in\mathbb{X}^{\mathbf{x},t}$.
Since the two trajectories agree up to $\tau_{\epsilon}$ and $\vartheta[\gamma_{\epsilon}]=\tau_{\epsilon}$, the stopping-time property implies $\vartheta[\eta_{\epsilon}]=\vartheta[\gamma_{\epsilon}]$. Consequently,
\begin{align}
\label{eq:converse-bound-1}
  \mathcal P(\mathbf{x},t;\vartheta,\eta_{\epsilon})=\int_{0}^{t\wedge\vartheta[\gamma_{\epsilon}]} L(\dot\gamma_{\epsilon}(s),\gamma_{\epsilon}(s),t-s)\,ds.
\end{align}

\textbf{Case 2.} Suppose $\vartheta[\gamma_{\epsilon}]\geqslant t$.
Since $L(\mathbf{v},\mathbf{x},s)\geqslant-r_{2}$, \eqref{eq:bounded-near-minimizer-bound} implies
\begin{align*}
  g_{\varrho,M}(\gamma_{\epsilon}(t)) \leqslant C_{K,0}+r_{2}t+\epsilon \leqslant C_{K,0}+r_{2}T+\epsilon <M.
\end{align*}
Since $g_{\varrho,M}=M$ outside $B_{2\varrho}(\mathbf{0})$, it follows that
\begin{align}
\label{eq:extend-endpoint-near-zero}
  |\gamma_{\epsilon}(t)|<2\varrho.
\end{align}
Set $\mathbf{z}_{\epsilon}=\gamma_{\epsilon}(t)$ and define
\begin{align*}
  \eta_{\epsilon}(s) =
  \begin{cases}
    \gamma_{\epsilon}(s), &s\in[0, t-\ell_{\varrho}),\\
    \gamma_{\epsilon}(s)-\frac{s-t+\ell_{\varrho}}{\ell_{\varrho}}\mathbf{z}_{\epsilon}, &s\in[t-\ell_{\varrho}, t),\\
    \mathbf{0}, &s\in[t,\infty).
  \end{cases}
\end{align*}
Then $\eta_{\epsilon}\in\mathbb{X}^{\mathbf{x},t}$.
Since $\eta_{\epsilon}=\gamma_{\epsilon}$ on $[0,t-\ell_{\varrho}]$ and $\vartheta[\gamma_{\epsilon}]\geqslant t$, we must have $\vartheta[\eta_{\epsilon}]>t-\ell_{\varrho}$.
Indeed, otherwise the stopping-time property, applied to $\eta_{\epsilon}$ and $\gamma_{\epsilon}$ on $[0,t-\ell_{\varrho}]$, would give $\vartheta[\gamma_{\epsilon}]=\vartheta[\eta_{\epsilon}]\leqslant t-\ell_{\varrho}$, a contradiction.
Set $\widehat\tau_{\epsilon}=t\wedge\vartheta[\eta_{\epsilon}]$. Then
\begin{align}
\label{eq:path-terminal-bound}
  t-\ell_{\varrho}<\widehat\tau_{\epsilon} \leqslant t.
\end{align}
Moreover, in the present case $\vartheta[\gamma_{\epsilon}]\geqslant t$, \eqref{eq:bounded-near-minimizer-bound}, $g_{\varrho,M}\geqslant0$ and $R\leqslant r_{2}$ provide the bound
\begin{align}
\label{eq:near-minimizer-energy-bound}
  \int_{0}^{t}|\dot\gamma_{\epsilon}(s)|^{2}\,ds \leqslant C_{K,1}:=4\left(C_{K,0}+r_{2}T+1\right).
\end{align}
On $(t-\ell_{\varrho},t)$, we have $\dot\eta_{\epsilon}(s)=\dot\gamma_{\epsilon}(s)-\frac{\mathbf{z}_{\epsilon}}{\ell_{\varrho}}$.
Since $\eta_{\epsilon}=\gamma_{\epsilon}$ on $[0,t-\ell_{\varrho}]$ and $\widehat\tau_{\epsilon}\leqslant t$, using \eqref{eq:extend-endpoint-near-zero}, \eqref{eq:near-minimizer-energy-bound}, and the Cauchy-Schwarz inequality, we have
\begin{equation}
\label{eq:bound-diff-kinetic}
  \begin{aligned}
    \frac{1}{4} \left(\int_{0}^{\widehat\tau_{\epsilon}}|\dot\eta_{\epsilon}(s)|^{2}\,ds-\int_{0}^{t}|\dot\gamma_{\epsilon}(s)|^{2}\,ds \right)
    \leqslant& \frac{1}{4}\int_{t-\ell_{\varrho}}^{t} \left| |\dot\eta_{\epsilon}(s)|^{2}-|\dot\gamma_{\epsilon}(s)|^{2} \right|ds\\
    \leqslant& \frac{|\mathbf{z}_{\epsilon}|}{2\ell_{\varrho}} \int_{t-\ell_{\varrho}}^{t}|\dot\gamma_{\epsilon}(s)|\,ds+\frac{|\mathbf{z}_{\epsilon}|^{2}}{4\ell_{\varrho}}\\
    \leqslant& \frac{\sqrt{C_{K,1}}|\mathbf{z}_{\epsilon}|}{2\sqrt{\ell_{\varrho}}}+\frac{|\mathbf{z}_{\epsilon}|^{2}}{4\ell_{\varrho}}\\
    \leqslant & \left(\sqrt{C_{K,1}}+1\right)\sqrt\varrho.
  \end{aligned}
\end{equation}
By \eqref{eq:path-terminal-bound}, the difference between the two potential terms is supported in $(t-\ell_{\varrho},t)$. Precisely,
\begin{align*}
  &\int_{0}^{\widehat\tau_{\epsilon}} R(\eta_{\epsilon}(s),t-s)\,ds-\int_{0}^{t} R(\gamma_{\epsilon}(s),t-s)\,ds
 =\int_{t-\ell_{\varrho}}^{\widehat\tau_{\epsilon}} R(\eta_{\epsilon}(s),t-s)\,ds-\int_{t-\ell_{\varrho}}^{t} R(\gamma_{\epsilon}(s),t-s)\,ds .
\end{align*}
Using $0\leqslant R\leqslant r_{2}$, we obtain
\begin{equation}
\label{eq:bound-diff-potential}
  \begin{aligned}
    \left| \int_{0}^{\widehat\tau_{\epsilon}} R(\eta_{\epsilon}(s),t-s)\,ds-\int_{0}^{t} R(\gamma_{\epsilon}(s),t-s)\,ds \right| \leqslant 2r_{2}\ell_{\varrho}
    \leqslant 2r_{2}\sqrt\varrho .
  \end{aligned}
\end{equation}
Combining \eqref{eq:bound-diff-kinetic} and \eqref{eq:bound-diff-potential} and using $g_{\varrho,M}(\gamma_{\epsilon}(t))\geqslant0$, we obtain
\begin{align}
\label{eq:converse-bound-2}
  \mathcal P(\mathbf{x},t;\vartheta,\eta_{\epsilon})
  \leqslant \int_{0}^{t} L(\dot\gamma_{\epsilon}(s),\gamma_{\epsilon}(s),t-s)\,ds+g_{\varrho,M}(\gamma_{\epsilon}(t))+C_{K,2}\sqrt\varrho,
\end{align}
where $C_{K,2}=\sqrt{C_{K,1}}+1+2r_{2}$.

To sum up, \eqref{eq:converse-bound-1} and \eqref{eq:converse-bound-2} imply that
\begin{align*}
  \mathcal P(\mathbf{x},t;\vartheta,\eta_{\epsilon})
  \leqslant \int_{0}^{t\wedge\vartheta[\gamma_{\epsilon}]} L(\dot\gamma_{\epsilon}(s),\gamma_{\epsilon}(s),t-s)\,ds+\mathbbm{1}_{\{\vartheta[\gamma_{\epsilon}]\geqslant t\}} g_{\varrho,M}(\gamma_{\epsilon}(t))+C_{K,2}\sqrt\varrho.
\end{align*}
Since $\eta_{\epsilon}\in \mathbb{X}^{\mathbf{x},t}$, taking the infimum over admissible trajectories, using the $\epsilon$-optimality of $\gamma_{\epsilon}$, and then letting $\epsilon\downarrow0$, we get
\begin{align*}
  \inf_{\eta\in\mathbb{X}^{\mathbf{x},t}} \mathcal P(\mathbf{x},t;\vartheta,\eta)
  \leqslant \inf_{\gamma\in \mathbb{A}^{\mathbf{x}}} \left\{\int_{0}^{t\wedge\vartheta[\gamma]} L(\dot\gamma(s),\gamma(s),t-s)\,ds+\mathbbm{1}_{\{\vartheta[\gamma]\geqslant t\}} g_{\varrho,M}(\gamma(t)) \right\}+C_{K,2}\sqrt\varrho.
\end{align*}
Taking the supremum over $\vartheta\in\Theta$ yields
\begin{align}
\label{eq:hard-value-above-bounded-value}
  W(\mathbf{x},t) \leqslant W_{\varrho,M}(\mathbf{x},t)+C_{K,2}\sqrt\varrho \quad \text{on }K,
\end{align}
provided $M>C_{K,0}+r_{2}T+1$.

We now choose an approximating sequence.
For $n\geqslant1$, set $\varrho_{n}=\frac{1}{n^{2}}$ and $M_{n}=n$.
Choose an arbitrary compact cylinder $K=\overline{B_{A}(\mathbf{0})}\times[\tau,T]$ with $A>0$ and $0<\tau<T$.
Since $\varrho_{n}\downarrow0$ and $M_{n}\uparrow\infty$, for all sufficiently large $n$ we have
\begin{align*}
  \varrho_{n}<\min\left\{1,\frac{\tau^{2}}{4} \right\}, \qquad M_{n}>\frac{A^{2}}{4\tau}+r_{2}T+1.
\end{align*}
Therefore, \eqref{eq:bounded-value-below-hard-value} and \eqref{eq:hard-value-above-bounded-value}, applied with $(\varrho,M)=(\varrho_{n},M_{n})$, yield
\begin{align*}
  0 \leqslant W(\mathbf{x},t)-W_{\varrho_{n},M_{n}}(\mathbf{x},t) \leqslant C_{K,2}\sqrt{\varrho_{n}}=\frac{C_{K,2}}{n}, \quad (\mathbf{x},t)\in K,
\end{align*}
for all sufficiently large $n$.
Since $C_{K,2}$ is independent of $n$, it follows that as $n\rightarrow\infty$,
\begin{align*}
  W_{\varrho_{n},M_{n}} \longrightarrow W \quad\text{uniformly on }K.
\end{align*}
Hence, by the arbitrariness of $K$, as $n\rightarrow\infty$,
\begin{align*}
  W_{\varrho_{n},M_{n}} \longrightarrow W \quad\text{locally uniformly in } \mathbb{R}^{N}\times\mathbb{R}^{+}.
\end{align*}

\textbf{Step 2.} We pass to the limit and prove uniqueness.
Each $W_{\varrho_{n},M_{n}}$ is a continuous viscosity solution of \eqref{eq:obstacle-direct} in $\mathbb{R}^{N}\times\mathbb{R}^{+}$.
Since $W_{\varrho_{n},M_{n}}$ converges locally uniformly to $W$ as $n\rightarrow\infty$, stability properties for viscosity solutions, together with Proposition \ref{prop:stopping-game-basic-properties}, show that $W$ is also a continuous viscosity solution of \eqref{eq:obstacle-direct} in $\mathbb{R}^{N}\times\mathbb{R}^{+}$.
It follows from Proposition \ref{prop:stopping-game-basic-properties} that $W$ is nonnegative and satisfies the growth control condition \eqref{eq:singular-asymptotic-precise}.
Hence Proposition \ref{prop:singular-comparison} implies that $W$ is the unique nonnegative continuous viscosity solution of \eqref{eq:obstacle-direct} satisfying \eqref{eq:singular-asymptotic-precise}.
In particular, if $w$ is the function obtained in Theorem \ref{thm:wkb-direct-convergence}, then
\begin{align*}
  W=w \quad\text{in }\mathbb{R}^{N}\times\mathbb{R}^{+}.
\end{align*}
This completes the proof.
\end{proof}

\subsection{The Wulff shape and spreading behavior}
Having identified $W$ with the geometric optics limit $w$, we now derive the spreading result from the geometry of the zero set of $W$.

\begin{proof}[Proof of Theorem \ref{thm:wulff}]

By the proof of Proposition \ref{prop:stopping-game-basic-properties}, $0\leqslant W(\mathbf{y},1) \leqslant \tfrac{|\mathbf{y}|^{2}}{4}$.
Hence $\mathbf{0}\in\mathcal{W}$.
Moreover, the lower bound in the same proposition yields $W(\mathbf{y},1) \geqslant \tfrac{|\mathbf{y}|^{2}}{4}-r_{2}$, and therefore
\begin{align}
\label{eq:wulff-contained-ball}
  \mathcal{W} \subset \overline{B_{2\sqrt{r_{2}}}(\mathbf{0})}.
\end{align}
Since $W(\cdot,1)$ is continuous, $\mathcal{W}$ is compact.

We next show that $\mathcal{W}$ has nonempty interior.
Fix $\mathbf{y}\in\overline{B_{2\sqrt{r_{1}}}(\mathbf{0})}$ and $\vartheta\in\Theta$. Consider
\begin{align*}
  \bar\gamma_\mathbf{y}(s) =
  \begin{cases}
    (1-s)\mathbf{y}, &s\in[0,1],\\
    \mathbf{0}, &s\in(1,\infty).
  \end{cases}
\end{align*}
Then $\bar\gamma_\mathbf{y}\in\mathbb{X}^{\mathbf{y},1}$.
Setting $\tau_{\vartheta}=1\wedge\vartheta[\bar\gamma_\mathbf{y}]$, and using $R\geqslant r_{1}$, we obtain
\begin{align*}
  \inf_{\gamma\in\mathbb{X}^{\mathbf{y},1}} \mathcal P(\mathbf{y},1;\vartheta,\gamma) \leqslant \mathcal P(\mathbf{y},1;\vartheta,\bar\gamma_\mathbf{y})
 =\int_{0}^{\tau_{\vartheta}} \left[ \frac{|\mathbf{y}|^{2}}{4}-R(\bar\gamma_\mathbf{y}(s),1-s) \right]ds
  \leqslant \tau_{\vartheta} \left(\frac{|\mathbf{y}|^{2}}{4}-r_{1} \right) \leqslant0.
\end{align*}
Taking the supremum over $\vartheta\in\Theta$ and recalling that $W\geqslant0$, we conclude that
\begin{align}
\label{eq:wulff-contains-ball}
  \overline{B_{2\sqrt{r_{1}}}(\mathbf{0})} \subset \mathcal{W}.
\end{align}
In particular, $\mathcal{W}$ has nonempty interior and $\textrm{dist}(\mathbf{0},\overline{\mathcal{W}}^{c})\geqslant 2\sqrt{r_{1}}$.
We complete the rest of the proof in the following two steps.

\textbf{Step 1.} We show the extinction property. We claim that
\begin{align}
\label{eq:extinction}
  \underset{\varepsilon\rightarrow 0}\lim\, u^{\varepsilon}=0 \quad \text{locally uniformly in }\{W>0\}.
\end{align}
Indeed, by Theorem \ref{thm:stopping-game-singular-obstacle-solution}, $W$ is the unique nonnegative continuous viscosity solution of \eqref{eq:obstacle-direct} with the growth control condition \eqref{eq:singular-asymptotic-precise}.
Then for every compact subset $K\subset\subset \{W>0\}$, there exists a constant $\kappa>0$ such that $W>\kappa$ in $K$.
Since $w^{\varepsilon}$ converges to $W$ locally uniformly in $\mathbb{R}^{N}\times\mathbb{R}^{+}$, we get that for $\varepsilon$ sufficiently small, $w^{\varepsilon}>\frac{\kappa}{2}$ in $K$. Therefore, we have
\begin{align*}
  u^{\varepsilon}=e^{-\frac{w^{\varepsilon}}{\varepsilon}} <e^{-\frac{\kappa}{2\varepsilon}} \quad \text{in } K.
\end{align*}
Since $u^{\varepsilon}$ is nonnegative, the claim holds.

Fix any $\eta>0$ and set
\begin{align*}
  K_{\eta}= \left\{\mathbf{y}\in\mathbb{R}^{N} \,\middle|\, \eta\leqslant\textrm{dist}(\mathbf{y},\mathcal{W}) \leqslant2\sqrt{r_{2}}+\eta \right\}.
\end{align*}
By \eqref{eq:wulff-contained-ball}, $K_{\eta}$ is compact.
Since every point of $K_{\eta}$ has positive distance from $\mathcal{W}$, we have $K_{\eta}\subset\subset\left\{\mathbf{y}\in\mathbb{R}^{N}\mid W(\mathbf{y},1)>0\right\}$.
In particular, $K_{\eta}\times\{1\}\subset\subset\{W>0\}$.
Then \eqref{eq:extinction}, restricted to $t=1$, yields
\begin{align}
\label{eq:extin-1}
  \underset{t\rightarrow\infty}\lim \underset{\textrm{dist}\{\frac{\mathbf{x}}{t},\mathcal{W}\}\in(\eta,2\sqrt{r_{2}}+\eta)} \sup\, u(\mathbf{x},t)
  =\underset{\varepsilon\rightarrow0} \lim \underset{\textrm{dist}\{\mathbf{y},\mathcal{W}\}\in(\eta,2\sqrt{r_{2}}+\eta)} \sup\,u^{\varepsilon}(\mathbf{y},1)=0.
\end{align}
Besides, since $\mathbf{0}\in\mathcal{W}$, $\textrm{dist}\left(\frac{\mathbf{x}}{t},\mathcal{W}\right)\leqslant\frac{|\mathbf{x}|}{t}$. Therefore,
\begin{align}
 \label{eq:set-unclusion}
  \left\{\mathbf{x}\mid \textrm{dist}\left(\frac{\mathbf{x}}{t},\mathcal{W}\right) \geqslant2\sqrt{r_{2}}+\eta \right\} \subset \left\{\mathbf{x}\mid|\mathbf{x}|\geqslant(2\sqrt{r_{2}}+\eta)t\right\}.
\end{align}
Let $S=\mathrm{supp}(g)$, $R_{S}=\sup_{\mathbf{z}\in S}|\mathbf{z}|$, $M=\max\{r_{2},\|g\|_{L^{\infty}}\}$, and $a=2\sqrt{r_{2}}+\eta$.
The heat kernel comparison used in the proof of Lemma \ref{lem:initial-lower-bound} gives, whenever $at\geqslant R_{S}$,
\begin{align*}
  \sup_{|\mathbf{x}|\geqslant at}u(\mathbf{x},t)
  \leqslant \frac{M|S|}{(4\pi t)^{N/2}}
  \exp\left(r_{2}t-\frac{(at-R_{S})^{2}}{4t}\right)
  \leqslant \frac{M|S|}{(4\pi t)^{N/2}}
  \exp\left[-\left(\frac{a^{2}}{4}-r_{2}\right)t
  +\frac{aR_{S}}{2}\right].
\end{align*}
Since $\frac{a^{2}}{4}-r_{2}>0$, the right-hand side tends to zero as $t\to\infty$.
Together with \eqref{eq:set-unclusion}, this yields
\begin{align}
\label{eq:extin-2}
  \underset{t\rightarrow\infty}\lim\,
  \underset{\textrm{dist}\{\frac{\mathbf{x}}{t},\mathcal{W}\}
  \geqslant2\sqrt{r_{2}}+\eta}\sup\,u(\mathbf{x},t)=0.
\end{align}
Combining \eqref{eq:extin-1} and \eqref{eq:extin-2}, we obtain the extinction assertion.

\textbf{Step 2.} We show the persistence property. We claim that
\begin{align}
\label{eq:persistence}
  \underset{\varepsilon\rightarrow 0}\liminf\, u^{\varepsilon}>0 \quad \text{locally uniformly in }\textrm{Int}\{W=0\}.
\end{align}
Let $S\subset\subset\textrm{Int}\{W=0\}$ be arbitrary. Choose $\delta>0$ such that
\begin{align*}
  S_{2\delta}=\left\{(\mathbf{x},t)\in\mathbb{R}^{N}\times\mathbb{R}^{+} \mid \textrm{dist}\{(\mathbf{x},t),S\}\leqslant2\delta\right\}
\end{align*}
is a compact subset of $\textrm{Int}\{W=0\}$.
Hence $W=0$ on $S_{2\delta}$.

Fix any $(\mathbf{x}_{0},t_{0})\in S$ and define
\begin{align*}
  \phi(\mathbf{x},t)=|\mathbf{x}-\mathbf{x}_{0}|^{2}+|t-t_{0}|^{2}.
\end{align*}
Since $W\equiv0$ on $S_{2\delta}$ and $w^\varepsilon\to W$ uniformly there, set
\begin{align*}
  \omega(\varepsilon)=\sup_{(\mathbf{x},t)\in S_{2\delta}} |w^{\varepsilon}(\mathbf{x},t)|.
\end{align*}
Then $\omega(\varepsilon)\to0$ as $\varepsilon\rightarrow0$.
For each $(\mathbf{x}_{0},t_{0})\in S$, let $(\mathbf{x}_{\varepsilon},t_{\varepsilon})$ be a maximum point of $w^{\varepsilon}-\phi$ on $\overline{B_{\delta}((\mathbf{x}_{0},t_{0}))}$.
Then we have $|\mathbf{x}_{\varepsilon}-\mathbf{x}_{0}|^{2} +|t_{\varepsilon}-t_{0}|^{2} \leqslant2\omega(\varepsilon)$.
Hence, for all sufficiently small $\varepsilon$, the maximum point $(\mathbf{x}_{\varepsilon},t_{\varepsilon})$ lies in $B_{\delta}((\mathbf{x}_{0},t_{0}))$, and $(\mathbf{x}_{\varepsilon},t_{\varepsilon})\rightarrow (\mathbf{x}_{0},t_{0})$ as $\varepsilon \rightarrow 0$, uniformly with respect to $(\mathbf{x}_{0},t_{0})\in S$.
At $(\mathbf{x}_{\varepsilon},t_{\varepsilon})$, $w^{\varepsilon}-\phi$ has a local maximum.
Hence, by \eqref{eq:w-eps-pde-main},
\begin{align*}
  u^{\varepsilon}(\mathbf{x}_{\varepsilon},t_{\varepsilon}) -r\left(\frac{\mathbf{x}_{\varepsilon}\cdot \mathbf{e}-ct_{\varepsilon}}{\varepsilon}\right)
  \geqslant \partial_{t}\phi(\mathbf{x}_{\varepsilon},t_{\varepsilon}) -\varepsilon\Delta \phi(\mathbf{x}_{\varepsilon},t_{\varepsilon}) +|D \phi(\mathbf{x}_{\varepsilon},t_{\varepsilon})|^{2}.
\end{align*}
Since $r\geqslant r_{1}$, it follows that for small $\varepsilon$, $u^{\varepsilon}(\mathbf{x}_{\varepsilon},t_{\varepsilon})\geqslant r_{1}-o(1)$, where the $o(1)$ term converges to $0$ as $\varepsilon\rightarrow0$, uniformly with respect to $(\mathbf{x}_{0},t_{0})\in S$.
Furthermore, by the Hopf-Cole transformation, we have $u^{\varepsilon}(\mathbf{x}_{0},t_{0})\geqslant u^{\varepsilon}(\mathbf{x}_{\varepsilon},t_{\varepsilon})$ for sufficiently small $\varepsilon$, and therefore,
\begin{align*}
  u^{\varepsilon}(\mathbf{x}_{0},t_{0}) \geqslant u^{\varepsilon}(\mathbf{x}_{\varepsilon},t_{\varepsilon}) \geqslant r_{1}-o(1) \quad \text{uniformly on } S.
\end{align*}
Hence
\begin{align*}
  \liminf_{\varepsilon\rightarrow0} \inf_{S} u^{\varepsilon} \geqslant r_{1}>0.
\end{align*}
This proves \eqref{eq:persistence}.

For any $\eta>0$, we set
\begin{align*}
  \bar K_{\eta}=\left\{\mathbf{y} \,\middle|\, \textrm{dist}(\mathbf{y},\overline{\mathcal{W}}^{c}) \geqslant\frac{\eta}{2} \right\}.
\end{align*}
Since $\mathcal{W}$ is compact, $\bar K_{\eta}$ is a compact subset of the spatial interior of $\mathcal{W}$; namely,
\begin{align*}
  \bar K_{\eta} \subset\subset \textrm{Int}\mathcal{W}=\textrm{Int} \left\{\mathbf{y}\in\mathbb{R}^{N} \mid W(\mathbf{y},1)=0 \right\}.
\end{align*}
Moreover, by the self-similarity, $\bar K_{\eta}\times\{1\}\subset\subset\textrm{Int}\left\{(\mathbf{y},s)\in\mathbb{R}^{N}\times\mathbb{R}^{+} \mid W(\mathbf{y},s)=0\right\}$.
Therefore, restricting \eqref{eq:persistence} to $t=1$, we obtain
\begin{align*}
  \liminf_{\varepsilon\to0} \inf_{\mathbf{y}\in\bar K_{\eta}} u^{\varepsilon}(\mathbf{y},1) \geqslant r_{1}>0.
\end{align*}
Since $\left\{\mathbf{y}\mid\textrm{dist}(\mathbf{y},\overline{\mathcal{W}}^{c})>\eta\right\} \subset \bar K_{\eta}$, we have
\begin{align*}
  \underset{t\rightarrow\infty}\liminf\, \underset{\mathrm{dist}\{\frac{\mathbf{x}}{t},\overline{\mathcal{W}}^{c}\}>\eta} \inf\, u(\mathbf{x},t)
  =\underset{\varepsilon\rightarrow0} \liminf \underset{\textrm{dist}\{\mathbf{y},\overline {\mathcal{W}}^{c}\}>\eta} \inf\, u^{\varepsilon}(\mathbf{y},1)
  \geqslant \underset{\varepsilon\rightarrow0} \liminf \inf_{\mathbf{y}\in\bar K_{\eta}}u^{\varepsilon}(\mathbf{y},1) >0.
\end{align*}
This completes the proof.
\end{proof}

\section{The unconstrained optimal control problem and Freidlin's condition}
\label{section:reductions-wulff-shapes}

We study the unconstrained optimal control value $J$ and determine exactly when it satisfies Freidlin's condition.
We first derive an explicit formula for $J$ and characterize its minimizing trajectories.
We then establish criteria for Freidlin's condition and prove the sharp dichotomy between its validity and failure regimes.
In Section \ref{section:stopping-game-reductions}, these results lead to two different reductions of the stopping game: an optimal control reduction in the Freidlin regimes and a pathwise characterization of the zero set in the complementary regime.

\subsection{The unconstrained optimal control problem}

Recall that $\mathbb{B}^{\mathbf{x},t}$ is given by \eqref{eq:hilbert-B}.
For $\gamma\in\mathbb B^{\mathbf{x},t}$, it admits a zero extension in $\mathbb{X}^{\mathbf{x},t}$.
For $s\in[0,t)$, define its time-shifted tail
\begin{align*}
  (\mathcal{Q}_{s}\gamma)(\tau)=\gamma(s+\tau), \quad \tau\in[0,t-s].
\end{align*}
Then $\mathcal{Q}_{s}\gamma\in\mathbb B^{\gamma(s),t-s}$.
Besides, we define an action functional associated to $\gamma\in\mathbb B^{\mathbf{x},t}$ as
\begin{align*}
  \mathcal I_{\mathrm{b}}(\mathbf{x},t;\gamma)=\int_{0}^{t} \left[ \frac{1}{4}|\dot\gamma(s)|^{2}-R(\gamma(s),t-s) \right]ds .
\end{align*}
Thus, the associated unconstrained value function is
\begin{align}
\label{eq:backward-unconstrained-value}
  J(\mathbf{x},t)=\inf_{\gamma\in\mathbb B^{\mathbf{x},t}} \mathcal I_{\mathrm{b}}(\mathbf{x},t;\gamma), \quad (\mathbf{x},t)\in \mathbb{R}^{N}\times\mathbb{R}^{+}.
\end{align}

\begin{proposition}
\label{prop:backward-sharp-value-function}
Recall the definitions $\alpha=\sqrt{r_{2}-r_{1}}$, $\lambda=\frac{c}{2}-\alpha$ and the orthogonal decomposition \eqref{eq:orthogonal-decompo}.
Then $J$ defined by \eqref{eq:backward-unconstrained-value} is explicitly solved as
\begin{align}
\label{eq:backward-J-decomposition}
  J(\mathbf{x},t)=\frac{|{\boldsymbol\zeta}|^{2}}{4t}+V(\xi,t), \quad (\mathbf{x},t)\in \mathbb{R}^{N}\times\mathbb{R}^{+},
\end{align}
where
\begin{align}
\label{eq:backward-V-piecewise}
  V(\xi,t)=
  \begin{cases}
    \frac{\xi^{2}}{4t}-r_{2}t, &\text{if } \xi\in[ct,\infty),\\
    \left(\frac {c}{2}-\alpha\right)\xi -\left(\frac{c^{2}}{4}+r_{2}-c\alpha\right)t, &\text{if } \xi\in(2\lambda t,ct),\\
    \frac{\xi^{2}}{4t}-r_{1}t, &\text{if } \xi\in(-\infty,2\lambda t].
  \end{cases}
\end{align}
Moreover, the minimum is uniquely attained by the following trajectory:
\begin{align*}
  \gamma_{*}(s)=Y_{*}(s)\mathbf{e}+ \left(1-\frac{s}{t}\right){\boldsymbol\zeta}, \quad s\in[0,t],
\end{align*}
where
\begin{align*}
  Y_{*}(s)=
  \begin{cases}
    \frac{\xi}{t}(t-s), & \text{if } \xi\in[ct,\infty),\\
    c(t-s), & \text{if } \xi\in(2\lambda t,ct) \text{ and } s\in[t-\tau_{*},t],\\
    c\tau_{*}+2\lambda(t-s-\tau_{*}), & \text{if } \xi\in(2\lambda t,ct) \text{ and } s\in[0,t-\tau_{*}),\\
    \frac{\xi}{t}(t-s), & \text{if } \xi\in(-\infty,2\lambda t],
  \end{cases}
\end{align*}
with
\begin{align}
\label{eq:tau-star}
  \tau_{*}=\frac{\xi-2\lambda t}{c-2\lambda}=\frac{\xi-2\lambda t}{2\alpha}.
\end{align}
\end{proposition}

\begin{proof}
Write $\gamma(s)=Y(s)\mathbf{e}+{\boldsymbol\eta}(s)$ with ${\boldsymbol\eta}(s)\cdot\mathbf{e}=0$.
The endpoint conditions are $Y(0)=\xi$, $Y(t)=0$, ${\boldsymbol\eta}(0)={\boldsymbol\zeta}$ and ${\boldsymbol\eta}(t)=\mathbf{0}$.
Since $R(\gamma(s),t-s)$ depends only on the moving coordinate $Y(s)-c(t-s)$ in the $\mathbf{e}$-direction, the component of $\gamma$ orthogonal to $\mathbf{e}$ decouples. By the Cauchy-Schwarz inequality,
\begin{align*}
  \int_{0}^{t}\frac{1}{4}|\dot{\boldsymbol\eta}(s)|^{2}\,ds \geqslant \frac{|{\boldsymbol\zeta}|^{2}}{4t},
\end{align*}
with equality for ${\boldsymbol\eta}(s)=\left(1-\frac{s}{t}\right){\boldsymbol\zeta}$.
Hence $J$ is reduced to \eqref{eq:backward-J-decomposition} with
\begin{align*}
  V(\xi,t)=\inf_{\substack{Y\in H^{1}([0,t];\mathbb{R})\\
  Y(0)=\xi,\ Y(t)=0}} \int_{0}^{t} \left[ \frac{1}{4}\dot Y(s)^{2}-r_{1}-\alpha^{2}\mathbbm{1}_{\{Y(s)\geqslant c(t-s)\}} \right]ds .
\end{align*}

Set $z(s)=Y(s)-c(t-s)$, $q=\xi-ct$.
Then $z(0)=q$, $z(t)=0$, and $\dot Y=\dot z-c$. Therefore
\begin{align}
\label{eq:backward-V-U}
  V(\xi,t)= \frac{c}{2}\xi-\frac{c^{2}}{4}t-r_{1}t+U(q,t),
\end{align}
where
\begin{align}
\label{eq:backward-U}
  U(q,t)=\inf_{\substack{z\in H^{1}([0,t];\mathbb{R})\\z(0)=q,\ z(t)=0}} \int_{0}^{t} \left[ \frac{1}{4}\dot z(s)^{2}-\alpha^{2}\mathbbm{1}_{\{z(s)\geqslant0\}} \right]ds .
\end{align}
We consider the following two cases to explicitly compute $U$.

\textbf{Case 1:} $q\geqslant0$.
For every admissible $z$, the Cauchy-Schwarz inequality implies
\begin{align*}
  \int_{0}^{t} \left[ \frac{1}{4}\dot z(s)^{2}-\alpha^{2}\mathbbm{1}_{\{z(s)\geqslant0\}} \right]ds
  \geqslant \frac{1}{4t} \left(\int_{0}^{t}\dot z(s)\,ds \right)^{2}-\alpha^{2}t=\frac{q^{2}}{4t}-\alpha^{2}t.
\end{align*}
The trajectory $z(s)=q(1-\frac{s}{t})$ stays in $\{z\geqslant0\}$ and attains equality. Therefore, $U(q,t)=\frac{q^{2}}{4t}-\alpha^{2}t$, which proves the first case in \eqref{eq:backward-V-piecewise}.

\textbf{Case 2:} $q<0$. For any admissible $z$, set
\begin{align*}
  E=\{s\in[0,t]\mid z(s)\geqslant0\}, \quad E^{c}=[0,t]\backslash E=\{s\in[0,t]\mid z(s)<0\}, \quad \tau_{z}=|E|.
\end{align*}
Since $z(0)=q<0$ and $z$ is continuous, $|E^{c}|=t-\tau_{z}>0$. Let $z^{-}=\max\{-z,0\}$.
Then $z^{-}\in H^{1}([0,t])$, $z^{-}(0)=-q$, $z^{-}(t)=0$, and $\dot z^{-}(s)=-\dot z(s)\mathbbm{1}_{E^{c}}(s)$ for a.e. $s\in[0,t]$.
Consequently,
\begin{align*}
  \int_{E^{c}} \dot z(s)\,ds=-\int_{0}^{t}\dot z^{-}(s)\,ds=z^{-}(0)-z^{-}(t)=-q .
\end{align*}
Furthermore, by the Cauchy-Schwarz inequality, we have
\begin{align*}
  \int_{E^{c}}\frac{1}{4}\dot z(s)^{2}\,ds \geqslant \frac{1}{4|E^{c}|}\left(\int_{E^{c}}\dot z(s)\,ds\right)^{2}=\frac{q^{2}}{4(t-\tau_{z})}.
\end{align*}
Since the kinetic energy on $E$ is nonnegative, for every admissible $z$, we obtain
\begin{align*}
  \int_{0}^{t} \left[ \frac{1}{4}\dot z(s)^{2}-\alpha^{2}\mathbbm{1}_{\{z(s)\geqslant0\}} \right]ds \geqslant \frac{q^{2}}{4(t-\tau_{z})}-\alpha^{2}\tau_{z}
  \geqslant \inf_{\tau\in[0,t)} \left\{\frac{q^{2}}{4(t-\tau)}-\alpha^{2}\tau \right\}.
\end{align*}
Taking the infimum over all admissible $z$, we have
\begin{align}
\label{eq:lower-inf-U}
  U(q,t) \geqslant \inf_{\tau\in[0,t)} \left\{\frac{q^{2}}{4(t-\tau)}-\alpha^{2}\tau \right\}.
\end{align}
Conversely, for each $\tau\in[0,t)$, define
\begin{align*}
  z_{\tau}(s)=
  \begin{cases}
    q\left(1-\frac{s}{t-\tau}\right), &\text{if } s\in[0,t-\tau],\\
    0, &\text{if } s\in(t-\tau,t].
  \end{cases}
\end{align*}
For every $\tau\in[0,t)$, the trajectory $z_{\tau}$ is admissible in \eqref{eq:backward-U}, and its total action is
\begin{align*}
  \int_{0}^{t} \left[ \frac{1}{4}\dot z_{\tau}(s)^{2}-\alpha^{2}\mathbbm{1}_{\{z_{\tau}(s)\geqslant0\}} \right]ds=\frac{q^{2}}{4(t-\tau)}-\alpha^{2}\tau .
\end{align*}
Taking the infimum over $\tau\in[0,t)$ and combining with \eqref{eq:lower-inf-U}, we obtain
\begin{align*}
  U(q,t)=\inf_{\tau\in[0,t)} \left\{\frac{q^{2}}{4(t-\tau)}-\alpha^{2}\tau \right\}.
\end{align*}
Let
\begin{align*}
  F(\tau)=\frac{q^{2}}{4(t-\tau)}-\alpha^{2}\tau, \quad \tau\in[0,t).
\end{align*}
If $q\in(-\infty,-2\alpha t]$, $F$ is nondecreasing on $[0,t)$, and hence $F$ attains its minimum at $\tau=0$.
If $q\in(-2\alpha t,0)$, then $F$ has the unique minimizer $\tau_{*}=t+\frac{q}{2\alpha}$. Therefore,
\begin{align*}
  U(q,t) =
  \begin{cases}
    \frac{q^{2}}{4t}, &\text{if }q\in(-\infty,-2\alpha t],\\[6pt]
    -\alpha q-\alpha^{2} t, &\text{if } q\in(-2\alpha t,0).
  \end{cases}
\end{align*}
Substitution into \eqref{eq:backward-V-U} yields the remaining two cases in \eqref{eq:backward-V-piecewise}. In the intermediate phase,
\begin{align*}
  \tau_{*}=t+\frac{q}{2\alpha}=\frac{\xi-2\lambda t}{c-2\lambda}, \qquad q=-2\alpha(t-\tau_{*}).
\end{align*}
Substituting the corresponding minimizer $z_{\tau_{*}}$ into $Y=z+c(t-s)$ gives the stated formula for $Y_{*}$.
In the other two phases, the straight minimizer $z(s)=q(1-\frac{s}{t})$ implies $Y_{*}(s)=\xi(1-\frac{s}{t})$.
Combining $Y_{*}$ with the minimizing orthogonal component yields $\gamma_{*}$.

The equality cases in the preceding estimates also give uniqueness of minimizing trajectory.
Indeed, equality in the orthogonal Cauchy-Schwarz estimate uniquely determines the orthogonal component.
In the homogeneous phases, equality in the corresponding Cauchy-Schwarz estimate determines the longitudinal component.
In the intermediate phase, strict convexity gives the unique minimizing time $\tau_{*}$, while equality in the lower-bound argument forces $z=z_{\tau_{*}}$.
Hence $\gamma_{*}$ is the unique minimizing trajectory.
\end{proof}

\begin{remark}
\label{rem:continuity-J}
It follows from Proposition \ref{prop:backward-sharp-value-function} that $J\in C(\mathbb{R}^{N}\times\mathbb{R}^{+})$.
Moreover, \eqref{eq:backward-J-decomposition} and \eqref{eq:backward-V-piecewise} imply that $J$ is strictly decreasing with respect to $t$.
\end{remark}

\begin{theorem}[The unconstrained HJE for the explicit value]
\label{thm:backward-J-viscosity}
The function $J$ defined by \eqref{eq:backward-unconstrained-value} is a continuous viscosity solution of
\begin{align}
\label{eq:backward-unconstrained-hje}
  J_{t}+|DJ|^{2}+R(\mathbf{x},t)=0 \quad\text{in }\mathbb{R}^{N}\times\mathbb{R}^{+}.
\end{align}
\end{theorem}

\begin{proof}
The continuity follows from Proposition \ref{prop:backward-sharp-value-function}.
On each of the three open regions $\{\xi>ct\}$, $\{2\lambda t<\xi<ct\}$, and $\{\xi<2\lambda t\}$, the corresponding formula in \eqref{eq:backward-J-decomposition}--\eqref{eq:backward-V-piecewise} is $C^{1}$, and direct differentiation gives \eqref{eq:backward-unconstrained-hje} pointwise, with $R=r_{2}$ in the first region and $R=r_{1}$ in the latter two.

It remains only to check the interfaces.
At $\xi=2\lambda t$, the two adjacent formulas are $C^{1}$-matched and both satisfy the same equation with $R=r_{1}$, so \eqref{eq:backward-unconstrained-hje} holds there in the classical, and hence viscosity, sense.
At $\xi=ct$, the one-sided slopes along $\xi$-direction are, respectively,
\begin{align*}
  \frac {c}{2}-\alpha \quad\text{and}\quad \frac {c}{2} .
\end{align*}
Thus the graph has an upward kink in the $\xi$-direction, and no $C^{1}$ test function can touch $J$ from above at such a point.
The subsolution condition holds automatically.
If $\phi$ touches $J$ from below at $(\mathbf{x},t)$ with $\xi=ct$, the touch condition implies, for some $\theta\in[0,1]$,
\begin{align*}
  D_{\boldsymbol\zeta}\phi=\frac{{\boldsymbol\zeta}}{2t}, \quad \phi_{\xi}=\frac {c}{2}-\alpha\theta\in\left[\frac{c}{2}-\alpha,\frac{c}{2}\right].
\end{align*}
Along the curve $s\rightarrow(\boldsymbol\zeta+c(t+s)\mathbf{e},t+s)$, $J-\phi$ attains its local minimum at $s=0$, so
\begin{align*}
  \phi_{t}+c\phi_{\xi}= \frac{d}{ds}J(\boldsymbol\zeta+c(t+s)\mathbf{e},t+s)|_{s=0}=-\frac{|\boldsymbol\zeta|^{2}}{4t^{2}} +\frac{c^{2}}{4}-r_{2}.
\end{align*}
Since $R^*(\mathbf{x},t)=r_{2}$ on the interface $\{\xi=ct\}$, we have
\begin{align*}
  \phi_{t}(\mathbf{x},t)+|D\phi|^{2}(\mathbf{x},t)+R^{*}(\mathbf{x},t)=\alpha^{2}\theta^{2} \geqslant0.
\end{align*}
Therefore, the supersolution condition holds at the discontinuous jumping interface.
\end{proof}

\subsection{Criteria and the sharp dichotomy for Freidlin's condition}

We now determine the sharp dichotomy for Freidlin's condition introduced in Definition \ref{def:backward-freidlin-condition}.
We first establish two useful criteria and then prove its sharp validity and failure regimes.

\begin{proposition}[A pointwise criterion for Freidlin's condition]
\label{prop:backward-freidlin-sufficient-condition}
Let $(\mathbf{x},t)\in\partial\mathbb P$.
Suppose that for every $\epsilon>0$ there exists $\gamma_{\epsilon}\in\mathbb B^{\mathbf{x},t}$ such that
\begin{align*}
  (\gamma_{\epsilon}(s),t-s)\in\mathbb P, \quad s\in(0,t),
\end{align*}
and
\begin{align*}
  \mathcal I_{\mathrm{b}}(\mathbf{x},t;\gamma_{\epsilon}) \leqslant J(\mathbf{x},t)+\epsilon.
\end{align*}
Then \eqref{eq:backward-freidlin-condition} holds at $(\mathbf{x},t)$.
\end{proposition}

\begin{proof}
The restricted infimum in \eqref{eq:backward-freidlin-condition} is always not smaller than the unrestricted value $J(\mathbf{x},t)$.
The assumed sequence implies the reverse inequality after letting $\epsilon\downarrow0$.
\end{proof}

\begin{corollary}[Perturbing a minimizing trajectory into the positive set]
\label{cor:backward-freidlin-perturbation}
Let $(\mathbf{x},t)\in\partial\mathbb P$ and $\gamma_{*}\in\mathbb B^{\mathbf{x},t}$ be a minimizing trajectory for $J(\mathbf{x},t)$.
Assume that this minimizer satisfies
\begin{align}
\label{eq:backward-relaxed-freidlin-condition}
  (\gamma_{*}(s),t-s)\in\overline{\mathbb P}, \quad  0\leqslant s<t,
\end{align}
where $\overline{\mathbb P}$ is the closure of the set $\mathbb P$.
Then \eqref{eq:backward-freidlin-condition} holds at $(\mathbf{x},t)$.
\end{corollary}

\begin{proof}
Let $\xi$ and ${\boldsymbol\zeta}$ be defined as in \eqref{eq:orthogonal-decompo}.
Then ${\boldsymbol\zeta}\cdot\mathbf{e}=0$. Set
\begin{align*}
  Y_{*}(s)=\gamma_{*}(s)\cdot\mathbf{e}, \qquad \mathbf{q}_{*}(s)=\gamma_{*}(s)-Y_{*}(s)\mathbf{e}.
\end{align*}
Since $\gamma_{*}\in\mathbb B^{\mathbf{x},t}$, we have $\mathbf{q}_{*}\in H^{1}([0,t];\mathbf{e}^{\perp})$, where $\mathbf{e}^{\perp}=\{\mathbf{z}\in\mathbb{R}^{N}|\mathbf{z}\cdot\mathbf{e}=0\}$.
In particular, $\mathbf{q}_{*}(0)={\boldsymbol\zeta}$ and $\mathbf{q}_{*}(t)=\mathbf{0}$.

For any ${\boldsymbol\eta}\in H^{1}([0,t];\mathbf{e}^{\perp})$ satisfying ${\boldsymbol\eta}(0)={\boldsymbol\zeta}$ and ${\boldsymbol\eta}(t)=\mathbf{0}$, define
\begin{align*}
  \tilde\gamma_{\boldsymbol\eta}(s)=Y_{*}(s)\mathbf{e}+{\boldsymbol\eta}(s), \quad s\in[0,t].
\end{align*}
The endpoint conditions give $\tilde\gamma_{\boldsymbol\eta}\in \mathbb B^{\mathbf{x},t}$.
Moreover, $\tilde\gamma_{\boldsymbol\eta}(s)\cdot\mathbf{e}=Y_{*}(s)=\gamma_{*}(s)\cdot\mathbf{e}$.
Since $R(\mathbf{y},\tau)$ depends on $\mathbf{y}$ only through $\mathbf{y}\cdot\mathbf{e}$, it follows that
\begin{align*}
  R(\tilde\gamma_{\boldsymbol\eta}(s),t-s)=R(\gamma_{*}(s),t-s).
\end{align*}
The orthogonal decomposition also shows that for almost every $s\in(0,t)$,
\begin{align*}
  |\dot\gamma_{*}(s)|^{2}=|\dot Y_{*}(s)|^{2}+|\dot{\mathbf{q}}_{*}(s)|^{2},
  \quad |\dot{\tilde\gamma}_{\boldsymbol\eta}(s)|^{2}=|\dot Y_{*}(s)|^{2}+|\dot{\boldsymbol\eta}(s)|^{2}.
\end{align*}
Therefore, the minimality of $\gamma_{*}$ yields
\begin{align*}
  0\leqslant \mathcal I_{\mathrm{b}} (\mathbf{x},t;\tilde\gamma_{\boldsymbol\eta})-\mathcal I_{\mathrm{b}} (\mathbf{x},t;\gamma_{*})
 =\frac{1}{4} \int_{0}^{t} \left(|\dot{\boldsymbol\eta}(s)|^{2}-|\dot{\mathbf{q}}_{*}(s)|^{2} \right)ds.
\end{align*}
Hence
\begin{align}
\label{eq:transverse-energy-minimality}
  \int_{0}^{t} |\dot{\mathbf{q}}_{*}(s)|^{2}ds \leqslant \int_{0}^{t} |\dot{\boldsymbol\eta}(s)|^{2}ds.
\end{align}
By the Cauchy-Schwarz inequality and \eqref{eq:transverse-energy-minimality} with ${\boldsymbol\eta}(s)=(1-\frac{s}{t}){\boldsymbol\zeta}$, we obtain
\begin{align*}
  \frac{|{\boldsymbol\zeta}|^{2}}{t} \leqslant \int_{0}^{t}|\dot{\mathbf{q}}_{*}(s)|^{2}ds \leqslant \frac{|{\boldsymbol\zeta}|^{2}}{t}.
\end{align*}
Thus equality holds, and the equality case in the Cauchy-Schwarz inequality yields
\begin{align}
\label{eq:backward-minimizer-orthogonal-component}
  \gamma_{*}(s)=Y_{*}(s)\mathbf{e}+\left(1-\frac{s}{t}\right){\boldsymbol\zeta}, \quad s\in[0,t].
\end{align}

If ${\boldsymbol\zeta}\neq\mathbf{0}$, set $\mathbf{m}=\frac{{\boldsymbol\zeta}}{|{\boldsymbol\zeta}|}$.
If ${\boldsymbol\zeta}=\mathbf{0}$, choose any unit vector $\mathbf{m}\perp\mathbf{e}$. Thus $|\mathbf{m}|=1$, $\mathbf{m}\cdot\mathbf{e}=0$ and ${\boldsymbol\zeta}\cdot\mathbf{m}=|{\boldsymbol\zeta}|$.
For $h>0$, define
\begin{align*}
  \gamma_{h}(s)=\gamma_{*}(s)+h s(t-s)\mathbf{m}, \quad  s\in[0,t].
\end{align*}
Then $\gamma_{h}\in\mathbb B^{\mathbf{x},t}$.
Since the perturbation $\gamma_{h}(s)-\gamma_{*}(s)=hs(t-s)\mathbf{m}$ is orthogonal to $\mathbf{e}$, we have
\begin{align}
\label{eq:backward-transverse-potential-invariance}
  R(\gamma_{h}(s),t-s)=R(\gamma_{*}(s),t-s), \quad  s\in[0,t).
\end{align}

Fix $s\in(0,t)$ and set $T=t-s$.
By \eqref{eq:backward-J-decomposition} and \eqref{eq:backward-minimizer-orthogonal-component},
\begin{align*}
  J(\gamma_{h}(s),T)-J(\gamma_{*}(s),T) &=\frac{1}{4T} \left(\left| \frac{T}{t}{\boldsymbol\zeta}+hsT\mathbf{m} \right|^{2}-\left| \frac{T}{t}{\boldsymbol\zeta} \right|^{2} \right)
 =\frac{hsT}{2t}|{\boldsymbol\zeta}|+\frac{h^{2}s^{2}T}{4} >0.
\end{align*}
Together with \eqref{eq:backward-relaxed-freidlin-condition}, this yields
\begin{align}
\label{eq:backward-perturbed-path}
  (\gamma_{h}(s),t-s)\in\mathbb P, \quad  s\in(0,t).
\end{align}

It remains to compare the actions.
From \eqref{eq:backward-minimizer-orthogonal-component},
\begin{align*}
  \dot\gamma_{*}(s)\cdot\mathbf{m}=-\frac{|{\boldsymbol\zeta}|}{t} \quad\text{for a.e. }s\in(0,t).
\end{align*}
Using \eqref{eq:backward-transverse-potential-invariance} and $\dot\gamma_{h}(s)=\dot\gamma_{*}(s)+h(t-2s)\mathbf{m}$, we obtain
\begin{align*}
  \mathcal I_{\mathrm{b}}(\mathbf{x},t;\gamma_{h})-\mathcal I_{\mathrm{b}}(\mathbf{x},t;\gamma_{*})
 =\frac {h}{2} \int_{0}^{t} (t-2s)\dot\gamma_{*}(s)\cdot\mathbf{m}\,ds+\frac{h^{2}}{4} \int_{0}^{t}(t-2s)^{2}\,ds=\frac{h^{2}t^{3}}{12}.
\end{align*}
Since $\mathcal I_{\mathrm{b}}(\mathbf{x},t;\gamma_{*})=J(\mathbf{x},t)$, for every $\epsilon>0$, choose $h\in(0,\sqrt{12\epsilon/t^{3}})$. Then
\begin{align*}
  \mathcal I_{\mathrm{b}}(\mathbf{x},t;\gamma_{h})=J(\mathbf{x},t)+\frac{h^{2}t^{3}}{12}<J(\mathbf{x},t)+\epsilon.
\end{align*}
Together with \eqref{eq:backward-perturbed-path}, Proposition \ref{prop:backward-freidlin-sufficient-condition} proves the assertion.
\end{proof}

\begin{theorem}[Freidlin's condition in the exterior speed regimes]
\label{thm:freidlin-sharp}
Assume that $c\leqslant-2\sqrt{r_{2}}$ or $c\geqslant 2\alpha=2\sqrt{r_{2}-r_{1}}$.
For every $(\mathbf{x},t)$ satisfying $J(\mathbf{x},t)=0$, let $\gamma_{*}\in\mathbb B^{\mathbf{x},t}$ be the minimizing backward trajectory described in Proposition \ref{prop:backward-sharp-value-function}. Then
\begin{align}
\label{eq:freidlin}
  J(\gamma_{*}(s),t-s)\geqslant0 \quad	\text{for all }s\in[0,t).
\end{align}
Consequently, $J$ satisfies Freidlin's condition.
\end{theorem}

\begin{proof}
Fix $(\mathbf{x},t)$ satisfying $J(\mathbf{x},t)=0$, denote $\xi=\mathbf{x}\cdot\mathbf{e}$ and ${\boldsymbol\zeta}=\mathbf{x}-\xi\mathbf{e}$, and set $F(s)=J(\gamma_{*}(s),t-s)$ for $s\in[0,t)$.
For each $s\in[0,t)$, Bellman's principle (see \cite[Chapter III]{bardi1997}) and minimality of $\gamma_{*}$ imply that its time-shifted tail is minimizing for $J(\gamma_{*}(s),t-s)$. Hence
\begin{align*}
  F(s)=\int_{s}^{t} \left[ \frac{1}{4}|\dot\gamma_{*}(\sigma)|^{2}-R(\gamma_{*}(\sigma),t-\sigma) \right]d\sigma .
\end{align*}

We first record a common observation for the two cases considered below.
If $\xi\geqslant ct$ or $\xi\leqslant 2\lambda t$, Proposition \ref{prop:backward-sharp-value-function} gives $\gamma_{*}(s)=(1-\frac{s}{t})\mathbf{x}$.
For $T=t-s$, we have $\gamma_{*}(s)\cdot\mathbf{e}=\frac{T}{t}\xi\geqslant cT$ in the first case and $\gamma_{*}(s)\cdot\mathbf{e}\leqslant 2\lambda T<cT$ in the second case.
Thus the straight trajectory remains in the same homogeneous branch for every $s\in[0,t)$, and the explicit formula for $J$ yields
\begin{align*}
  F(s)=\frac{T}{t}J(\mathbf{x},t)=0.
\end{align*}
Thus $F(s)=0$ in either homogeneous phase.

We now distinguish the two alternatives in the hypothesis.
First suppose that $c\leqslant-2\sqrt{r_{2}}$.
Since $R\leqslant r_{2}$, every $\gamma\in\mathbb B^{\mathbf{x},t}$ satisfies
\begin{align*}
  \mathcal I_{\mathrm{b}}(\mathbf{x},t;\gamma) \geqslant \frac{1}{4}\int_{0}^{t}|\dot\gamma(s)|^{2}ds-r_{2}t
  \geqslant \frac{1}{4t}\left|\int_{0}^{t}\dot\gamma(s)ds\right|^{2}-r_{2}t=\frac{|\mathbf{x}|^{2}}{4t}-r_{2}t,
\end{align*}
where the second inequality follows from Cauchy-Schwarz inequality and the endpoint conditions.
Taking the infimum and using $J(\mathbf{x},t)=0$, we obtain $|\mathbf{x}|\leqslant2\sqrt{r_{2}}t$ and hence
\begin{align*}
  \xi\geqslant-|\mathbf{x}|\geqslant-2\sqrt{r_{2}}t\geqslant ct.
\end{align*}
Thus $(\mathbf{x},t)$ belongs to the high-growth homogeneous phase, and the preceding calculation gives $F(s)=0$ for every $s\in[0,t)$.

Suppose next that $c\geqslant2\sqrt{r_{2}-r_{1}}$.
The preceding observation has already treated the two homogeneous phases, so it remains only to consider $2\lambda t<\xi<ct$. Set
\begin{align}
\label{eq:aux-para}
  \nu=\frac{\tau_{*}}{t}\in(0,1), \quad A=\frac{|{\boldsymbol\zeta}|^{2}}{4t^{2}}, \quad K_{1}=\frac{c^{2}}{4}-r_{2}, \quad \text{and} \quad K_{2}=\lambda^{2}-r_{1}.
\end{align}
where $\tau_{*}$ is given by \eqref{eq:tau-star}.
Substituting the minimizing trajectory from Proposition \ref{prop:backward-sharp-value-function} into the integral representation of $F$, we obtain
\begin{align*}
  F(s)=
  \begin{cases}
    (t-\tau_{*}-s)(A+K_{2})+\tau_{*}(A+K_{1}), &\text{if }s\in[0,t-\tau_{*}),\\
    (t-s)(A+K_{1}), &\text{if }s\in[t-\tau_{*},t).
  \end{cases}
\end{align*}
Since $F(0)=0$, we have $A+\nu K_{1}+(1-\nu)K_{2}=0$.
Together with $K_{1}-K_{2}=2\lambda\alpha$, this implies
\begin{align*}
  A+K_{1}=2\lambda\alpha(1-\nu) \quad \text{and} \quad A+K_{2}=-2\lambda\nu\alpha.
\end{align*}
Therefore,
\begin{align}
\label{eq:backward-front-value-profile}
  F(s) =
  \begin{cases}
    2\lambda\nu\alpha s, &s\in[0,(1-\nu)t),\\
    2\lambda\alpha(t-s)(1-\nu), &s\in[(1-\nu)t,t).
  \end{cases}
\end{align}
Since $\lambda\geqslant0$, \eqref{eq:backward-front-value-profile} implies $F(s)\geqslant0$.
In particular, if $\lambda>0$, the same formula and $\nu\in(0,1)$ give $J(\gamma_{*}(s),t-s)=F(s)>0$ for every $s\in(0,t)$.
This proves \eqref{eq:freidlin}.

We finally verify the relaxed condition \eqref{eq:backward-relaxed-freidlin-condition}.
Let $(\mathbf{x},t)\in\partial\mathbb P$ be arbitrary.
The continuity of $J$ and the definition of $\mathbb P$ imply $J(\mathbf{x},t)=0$, so \eqref{eq:freidlin} applies to the corresponding minimizer $\gamma_{*}$. Fix $s\in[0,t)$.
If $J(\gamma_{*}(s),t-s)>0$, then $(\gamma_{*}(s),t-s)\in\mathbb P$.
If $J(\gamma_{*}(s),t-s)=0$, the strict decrease of $J(\mathbf{y},\cdot)$ stated in Remark~\ref{rem:continuity-J} gives $J(\gamma_{*}(s),t-s-\delta)>0$ for every $\delta\in(0,t-s)$.
Letting $\delta\downarrow0$, we obtain
\begin{align*}
  (\gamma_{*}(s),t-s)\in\overline{\mathbb P}, \quad s\in[0,t).
\end{align*}
Since the boundary point is arbitrary, Corollary \ref{cor:backward-freidlin-perturbation} shows that \eqref{eq:backward-freidlin-condition} holds at every point of $\partial\mathbb P$.
Hence, $J$ satisfies Freidlin's condition.
\end{proof}

\begin{proposition}[Failure of Freidlin's condition in the intermediate regime]\label{prop:freidlin-failure-full-range}
If
\begin{align}
\label{eq:freidlin-failure-full-range}
  -2\sqrt{r_{2}} <c< 2\alpha,
\end{align}
then the value function $J$ does not satisfy Freidlin's condition \eqref{eq:backward-freidlin-condition}.
More precisely, for every $t\in\mathbb{R}^{+}$, there exists $(\mathbf{x},t)\in\partial\mathbb P$ such that
\begin{align*}
  0=J(\mathbf{x},t)<\inf_{\gamma\in\mathbb B^{\mathbf{x},t}} \left\{\mathcal I_{\mathrm{b}}(\mathbf{x},t;\gamma) \mid (\gamma(s),t-s)\in\mathbb P \text{ for }s\in(0,t) \right\}.
\end{align*}
\end{proposition}

\begin{proof}
Let $K_{1}$ and $K_{2}$ be defined as in \eqref{eq:aux-para}.
\eqref{eq:freidlin-failure-full-range} then implies $\lambda=\frac{c}{2}-\alpha<0$ and $K_{1}<0$.
Choose $\mu\in(0,1)$ sufficiently close to $1$ such that $\bar{K}=-\mu K_{1}-(1-\mu)K_{2}>0$.
Fix any $t\in\mathbb{R}^{+}$ and choose a vector ${\boldsymbol\zeta}\perp\mathbf{e}$ such that $|{\boldsymbol\zeta}|=2t\sqrt{\bar{K}}$. Define
\begin{align*}
  \xi=t\left[ \mu c+2\lambda(1-\mu) \right], \quad \mathbf{x}=\xi\mathbf{e}+{\boldsymbol\zeta}.
\end{align*}
Then $2\lambda t<\xi<ct$, and $\tau_{*}$ given by \eqref{eq:tau-star} satisfies $\tau_{*}=\mu t$.
Using \eqref{eq:backward-J-decomposition} and \eqref{eq:backward-V-piecewise}, we obtain
\begin{align*}
  J(\mathbf{x},t)=t\left[ \bar{K}+\mu K_{1}+(1-\mu)K_{2} \right]=0.
\end{align*}

We verify that $(\mathbf{x},t)\in\partial\mathbb P$. For $\epsilon\in(0,1)$, let $\mathbf{x}_{\epsilon}^{\pm}=\xi\mathbf{e}+(1\pm\epsilon){\boldsymbol\zeta}$.
Then we have
\begin{align*}
  J(\mathbf{x}_{\epsilon}^{\pm},t)=t\bar K\left[ (1\pm\epsilon)^{2}-1 \right].
\end{align*}
Thus $J(\mathbf{x}_{\epsilon}^{+},t)>0$ and $J(\mathbf{x}_{\epsilon}^{-},t)<0$.
Letting $\epsilon\downarrow0$ proves $(\mathbf{x},t)\in\partial\mathbb P$.

Let $\gamma_{*}$ be the minimizing trajectory given in Proposition \ref{prop:backward-sharp-value-function}, and define
\begin{align*}
  F(s)=J(\gamma_{*}(s),t-s), \qquad s\in[0,t).
\end{align*}
By Bellman's principle and optimality of $\gamma_{*}$, the tail $\mathcal{Q}_s\gamma_{*}$ is minimizing for $J(\gamma_{*}(s),t-s)$.
Along $\gamma_{*}$, the action density equals $\bar K+K_{2}$ on $[0,t-\tau_{*})$ and $\bar K+K_{1}$ on $[t-\tau_{*},t]$. Since $\tau_{*}=\mu t$,
\begin{align*}
  \bar{K}+K_{1}=2\lambda(1-\mu)\alpha, \quad \bar{K}+K_{2}=-2\lambda\mu\alpha.
\end{align*}
Together with $F(0)=J(\mathbf{x},t)=0$, we obtain
\begin{align*}
  F(s) =
  \begin{cases}
    2\lambda\mu\alpha s, &s\in[0,(1-\mu)t),\\
    2\lambda\alpha(t-s)(1-\mu), &s\in[(1-\mu)t, t).
  \end{cases}
\end{align*}
It follows from $\lambda<0$ and $\mu\in(0,1)$ that
\begin{align}
\label{eq:freidlin-failure-strict-negativity}
  J(\gamma_{*}(s),t-s)<0, \quad  s\in(0,t).
\end{align}

Suppose now, to the contrary, that the restricted infimum in \eqref{eq:backward-freidlin-condition} were zero.
Then there exists a sequence $\{\gamma_{n}\}_{n\geqslant 1}\subset\mathbb B^{\mathbf{x},t}$ such that $(\gamma_{n}(s),t-s)\in\mathbb P$ for $s\in(0,t)$ and
\begin{align*}
  0 \leqslant \mathcal I_{\mathrm{b}}(\mathbf{x},t;\gamma_{n}) \leqslant \frac{1}{n}.
\end{align*}
Since $R\leqslant r_{2}$, we have
\begin{align}
\label{eq:bound-kinetic}
  \frac{1}{4}\int_{0}^{t}|\dot\gamma_{n}(s)|^{2}ds \leqslant \mathcal I_{\mathrm{b}}(\mathbf{x},t;\gamma_{n})+r_{2}t \leqslant 1+r_{2}t.
\end{align}
Recalling that $\gamma_{n}(t)=\mathbf{0}$, the Poincar\'e inequality and \eqref{eq:bound-kinetic} yield
\begin{align*}
  \|\gamma_{n}\|_{H^{1}(0,t)}^{2} \leqslant (1+t^{2})\|\dot\gamma_{n}\|_{L^{2}(0,t)}^{2} \leqslant 4(1+t^{2})(1+r_{2}t).
\end{align*}
Hence $\{\gamma_{n}\}_{n}$ is bounded in $H^{1}([0,t];\mathbb{R}^{N})$.
After passing to a subsequence, we may assume that
\begin{align*}
  \gamma_{n}\rightharpoonup\overline\gamma \quad\text{in }H^{1}([0,t];\mathbb{R}^{N}), \quad \gamma_{n}\rightarrow\overline\gamma \quad\text{uniformly on }[0,t].
\end{align*}
The endpoint conditions pass to the limit, and so $\overline\gamma\in\mathbb B^{\mathbf{x},t}$.
By \cite[Proposition 3.5]{brezis2011}, we have
\begin{align*}
  \int_{0}^{t}|\dot{\overline\gamma}(s)|^{2}ds
  \leqslant\liminf_{n\to\infty}
  \int_{0}^{t}|\dot\gamma_n(s)|^{2}ds.
\end{align*}
Since $\gamma_n\to\overline\gamma$ uniformly and $R$ is u.s.c., for every $s\in(0,t)$ we have
\begin{align*}
  -r_{2}\leqslant-R(\overline\gamma(s),t-s)
  \leqslant\liminf_{n\to\infty}
  \left[-R(\gamma_n(s),t-s)\right].
\end{align*}
Fatou's lemma therefore yields
\begin{align*}
  -\int_{0}^{t}R(\overline\gamma(s),t-s)\,ds
  \leqslant\liminf_{n\to\infty}
  \int_{0}^{t}\left[-R(\gamma_n(s),t-s)\right]ds.
\end{align*}
Combining these two inequalities, we obtain
\begin{align}
\label{eq:lsc-inequality}
  \mathcal I_{\mathrm{b}}(\mathbf{x},t;\overline\gamma)
  \leqslant\liminf_{n\to\infty}
  \mathcal I_{\mathrm{b}}(\mathbf{x},t;\gamma_n)=0.
\end{align}
On the other hand, since $\overline\gamma\in\mathbb B^{\mathbf{x},t}$ and $J(\mathbf{x},t)=0$, one has $\mathcal I_{\mathrm{b}}(\mathbf{x},t;\overline\gamma) \geqslant J(\mathbf{x},t)=0$.
Thus $\overline\gamma$ is a minimizing trajectory.
By the uniqueness assertion in Proposition \ref{prop:backward-sharp-value-function}, we obtain $\overline\gamma=\gamma_{*}$.

For each fixed $s\in(0,t)$, the choice of $\gamma_n$ gives $J(\gamma_n(s),t-s)>0$.
Hence, by the continuity of $J$ and the uniform convergence of $\gamma_n$, we have
\begin{align*}
  J(\gamma_{*}(s),t-s)=\lim_{n\to\infty}J(\gamma_n(s),t-s) \geqslant0.
\end{align*}
This contradicts \eqref{eq:freidlin-failure-strict-negativity}.
Hence the restricted infimum in \eqref{eq:backward-freidlin-condition} cannot be zero.
Since it is no smaller than $J(\mathbf{x},t)=0$, it is strictly positive.
Therefore, Freidlin's condition fails at $(\mathbf{x},t)$.
\end{proof}

\section{Variational reductions of the stopping game}
\label{section:stopping-game-reductions}

This section develops the variational reductions underlying the Wulff shape analysis.
We first reduce the stopping value to the unconstrained value in the Freidlin regimes.
We then derive a min-max representation and a pathwise characterization of the zero set valid for every $c\in\mathbb{R}$.
In the non-Freidlin regime, we further obtain a pointwise characterization and a convex hull representation.

For each $t\in\mathbb{R}^{+}$, define the propagation region
\begin{align*}
  \mathcal{Z}(t)=\left\{\mathbf{x}\in\mathbb{R}^{N}\mid W(\mathbf{x},t)=0 \right\}.
\end{align*}
By Proposition \ref{prop:self-similarity},
\begin{align}
\label{eq:propagation-region-scaling}
  \mathcal{Z}(t)=t\mathcal{Z}(1), \quad \partial\mathcal{Z}(t)=t\partial\mathcal{Z}(1).
\end{align}
Accordingly,
\begin{align}
\label{eq:wulff-shape-definition}
  \mathcal{W}=\mathcal{Z}(1), \quad \partial \mathcal{W}= \partial\mathcal{Z}(1).
\end{align}

\subsection{Optimal control reduction in the Freidlin regimes}

\begin{theorem}[Reduction to the unconstrained value under Freidlin's condition]
\label{thm:backward-W-is-viscosity-solution}
Assume that $J$ defined by \eqref{eq:backward-unconstrained-value} satisfies Freidlin's condition \eqref{eq:backward-freidlin-condition}.
Then the stopping game value $W$ satisfies
\begin{align}
\label{eq:backward-W-max-J}
  W(\mathbf{x},t)=\max\{0,J(\mathbf{x},t)\} \quad\text{in }\mathbb{R}^{N}\times\mathbb{R}^{+}.
\end{align}
\end{theorem}

\begin{proof}
Set $w(\mathbf{x},t)=\max\{0,J(\mathbf{x},t)\}$.
Since $J$ is continuous in $\mathbb{R}^{N}\times\mathbb{R}^{+}$, so is $w$.
We first prove that $w$ is a viscosity solution of \eqref{eq:obstacle-direct}.

Let $\phi\in C^{1}$ touch $w$ from above at $(\mathbf{x}_{0},t_{0})\in \mathbb{R}^{N}\times\mathbb{R}^{+}$.
If $w(\mathbf{x}_{0},t_{0})=0$, the subsolution inequality follows from the obstacle term.
If $w(\mathbf{x}_{0},t_{0})>0$, then $w=J$ near $(\mathbf{x}_{0},t_{0})$, and the subsolution inequality follows from Theorem \ref{thm:backward-J-viscosity}.

We next verify the supersolution property.
Let $\phi\in C^{1}$ touch $w$ from below at $(\mathbf{x}_{0},t_{0})\in \mathbb{R}^{N}\times\mathbb{R}^{+}$.
Since $w\geqslant0$, it is enough to prove
\begin{align}
\label{eq:super-pde-goal}
  \phi_{t}(\mathbf{x}_{0},t_{0})+|D\phi(\mathbf{x}_{0},t_{0})|^{2}+R^{*}(\mathbf{x}_{0},t_{0}) \geqslant0.
\end{align}
If $(\mathbf{x}_{0},t_{0})\in\mathbb P$, then $w=J$ locally, and \eqref{eq:super-pde-goal} follows from Theorem \ref{thm:backward-J-viscosity}.
If $(\mathbf{x}_{0},t_{0})$ is an interior point of $\{w=0\}$, then $D\phi(\mathbf{x}_{0},t_{0})=0$ and $\phi_{t}(\mathbf{x}_{0},t_{0})=0$, so \eqref{eq:super-pde-goal} holds because $R\geqslant0$.

It remains to consider $(\mathbf{x}_{0},t_{0})\in\partial\mathbb P$.
Then $J(\mathbf{x}_{0},t_{0})=w(\mathbf{x}_{0},t_{0})=0$.
Suppose by contradiction that, for some $\eta>0$,
\begin{align*}
  \phi_{t}(\mathbf{x}_{0},t_{0})+|D\phi(\mathbf{x}_{0},t_{0})|^{2}+R^{*}(\mathbf{x}_{0},t_{0}) \leqslant -4\eta .
\end{align*}
Choose $\delta_{0}\ll1$ such that $\phi\leqslant w$ in $B_{\delta_{0}}(\mathbf{x}_{0}) \times(t_{0}-\delta_{0},t_{0}+\delta_{0})$ and, throughout this neighborhood,
\begin{align*}
  \phi_{t}+|D\phi|^{2}+R\leqslant-2\eta .
\end{align*}
For $h\in(0,\min\{t_{0},\delta_{0}\})$, by Freidlin's condition and the definition of the infimum, applied with approximation error $h^{2}$, there exists $\gamma_{h}\in\mathbb B^{\mathbf{x}_{0},t_{0}}$ such that
\begin{align}
\label{eq:boundary-full-action-small}
  (\gamma_{h}(s),t_{0}-s)\in\mathbb P, \quad s\in(0,t_{0}); \quad \text{and} \quad  \mathcal I_{\mathrm{b}}(\mathbf{x}_{0},t_{0};\gamma_{h})\leqslant h^{2} .
\end{align}
The time-shifted tail $\mathcal{Q}_{h}\gamma_{h}\in\mathbb B^{\gamma_{h}(h),t_{0}-h}$ is admissible for $J(\gamma_{h}(h),t_{0}-h)$. Therefore
\begin{align*}
  J(\gamma_{h}(h),t_{0}-h) \leqslant \int_{h}^{t_{0}} \left[ \frac{1}{4}|\dot\gamma_{h}(s)|^{2}-R(\gamma_{h}(s),t_{0}-s) \right]ds .
\end{align*}
Combining this with \eqref{eq:boundary-full-action-small}, and using $w=J$ at $(\gamma_{h}(h),t_{0}-h)$, we have
\begin{align}
\label{eq:boundary-first-segment-small}
  \int_{0}^{h} \left[ \frac{1}{4}|\dot\gamma_{h}(s)|^{2}-R(\gamma_{h}(s),t_{0}-s) \right]ds+w(\gamma_{h}(h),t_{0}-h) \leqslant h^{2} .
\end{align}
Since $J(\gamma_{h}(h),t_{0}-h)>0$ and $R\leqslant r_{2}$, \eqref{eq:boundary-first-segment-small} implies that
\begin{align}
\label{eq:energy-bound}
  \frac{1}{4}\int_{0}^{h}|\dot\gamma_{h}(s)|^{2}\,ds \leqslant h^{2}+r_{2}h .
\end{align}
We claim that $(\gamma_{h}(s),t_{0}-s)\in B_{\delta_{0}}(\mathbf{x}_{0})\times(t_{0}-\delta_{0},t_{0}+\delta_{0})$ for $0\leqslant s\leqslant h$, provided $h$ is sufficiently small.
Indeed, the Cauchy-Schwarz inequality and \eqref{eq:energy-bound} give, for every $s\in[0,h]$,
\begin{align*}
  |\gamma_{h}(s)-\mathbf{x}_{0}|=\left| \int_{0}^{s}\dot\gamma_{h}(\sigma)\,d\sigma \right| \leqslant 2h\sqrt{h+r_{2}}.
\end{align*}
Moreover, $|(t_{0}-s)-t_{0}|=s\leqslant h$.
Choose $h>0$ sufficiently small that $h<\delta_{0}$ and $2h\sqrt{h+r_{2}}<\delta_{0}$.
This proves the claim.

Since $\phi$ touches $w$ from below and $w(\mathbf{x}_{0},t_{0})=0$,
\begin{align*}
  \phi(\gamma_{h}(h),t_{0}-h)-\phi(\mathbf{x}_{0},t_{0}) \leqslant w(\gamma_{h}(h),t_{0}-h).
\end{align*}
Hence
\begin{align}
\label{eq:boundary-phi-action-upper}
  \int_{0}^{h} \left[ \frac{1}{4}|\dot\gamma_{h}(s)|^{2}-R(\gamma_{h}(s),t_{0}-s) \right]ds+\phi(\gamma_{h}(h),t_{0}-h)-\phi(\mathbf{x}_{0},t_{0}) \leqslant h^{2} .
\end{align}
On the other hand, by the fundamental theorem of calculus and Fenchel's inequality, one has
\begin{align*}
  &\int_{0}^{h} \left[ \frac{1}{4}|\dot\gamma_{h}(s)|^{2}-R(\gamma_{h}(s),t_{0}-s) \right]ds+\phi(\gamma_{h}(h),t_{0}-h)-\phi(\mathbf{x}_{0},t_{0})\\
  \geqslant& \int_{0}^{h} \left[ -|D\phi|^{2}-\phi_{t}-R \right](\gamma_{h}(s),t_{0}-s)\,ds \geqslant 2\eta h ,
\end{align*}
which contradicts \eqref{eq:boundary-phi-action-upper} for $h<2\eta$.
Hence $w$ is a nonnegative viscosity solution of \eqref{eq:obstacle-direct}.

Finally, since $r_{1}\leqslant R\leqslant r_{2}$, by using the Cauchy-Schwarz inequality for the lower bound and by testing with the straight trajectory for the upper bound, we have
\begin{align*}
  \frac{|\mathbf{x}|^{2}}{4t}-r_{2}t \leqslant J(\mathbf{x},t) \leqslant \frac{|\mathbf{x}|^{2}}{4t}-r_{1}t .
\end{align*}
Hence $w=\max\{0,J\}$ satisfies \eqref{eq:singular-asymptotic-precise} with $C=r_{2}$.
By Theorem \ref{thm:stopping-game-singular-obstacle-solution}, the stopping value $W$ is the unique nonnegative continuous viscosity solution of \eqref{eq:obstacle-direct} with this growth control condition.
Therefore $W=w$, which proves \eqref{eq:backward-W-max-J}.
\end{proof}

In the Freidlin regimes, Theorem \ref{thm:backward-W-is-viscosity-solution} identifies the propagation region with the zero sublevel set of $J$.
The corresponding boundary is identified in the next lemma, while the explicit formulas are derived later.

\begin{lemma}
\label{lem:freidlin-zero-level-characterization}
Assume that $c\leqslant-2\sqrt{r_{2}}$ or $c\geqslant2\alpha$.
Then, for every $t\in\mathbb{R}^{+}$,
\begin{align}
\label{eq:freidlin-zero-level-front}
  \mathcal{Z}(t)=\left\{\mathbf{x}\in\mathbb{R}^{N} \mid J(\mathbf{x},t)\leqslant0 \right\},
  \quad \partial\mathcal{Z}(t)=\left\{\mathbf{x}\in\mathbb{R}^{N} \mid J(\mathbf{x},t)=0 \right\}.
\end{align}
\end{lemma}

\begin{proof}
By Theorem \ref{thm:freidlin-sharp}, $J$ satisfies Freidlin's condition.
Theorem \ref{thm:backward-W-is-viscosity-solution} therefore concludes that $W=\max\{0,J\}$.
Moreover, Theorems \ref{thm:wkb-direct-convergence} and \ref{thm:stopping-game-singular-obstacle-solution} give $w=W$.
Consequently, $\mathcal{Z}(t)=\{J(\cdot,t)\leqslant 0\}$.

We next verify that $\partial\mathcal{Z}(t)=\{J(\cdot,t)=0\}$.
Since $J(\cdot,t)$ is continuous, $\partial\{J(\cdot,t)\leqslant0\}\subset\{J(\cdot,t)=0\}$.
Conversely, fix $\mathbf{x}\in\mathbb{R}^{N}$ such that $J(\mathbf{x},t)=0$.
Notice that $\mathbf{x}\neq\mathbf{0}$, since the explicit formula gives $J(\mathbf{0},t)<0$.
For $h\in(0,t)$, the strict decrease of $J$ with respect to $t$ implies that
\begin{align*}
  J(\mathbf{x},t-h)>J(\mathbf{x},t)=0.
\end{align*}
Set $k_{h}=\frac{t}{t-h}>1$.
Using $W=w$, the self-similarity of $w$, and $W=\max\{0,J\}$, we obtain
\begin{align*}
  W(k_{h}\mathbf{x},t)=w(k_{h}\mathbf{x},k_{h}(t-h))=k_{h}w(\mathbf{x},t-h)=k_{h}J(\mathbf{x},t-h) >0.
\end{align*}
Since $k_{h}\mathbf{x}\notin\mathcal{Z}(t)$ and $k_{h}\mathbf{x}\to\mathbf{x}\in\mathcal{Z}(t)$ as $h\downarrow0$, we have $\mathbf{x}\in\partial\mathcal{Z}(t)$.
This proves \eqref{eq:freidlin-zero-level-front}.
\end{proof}

\subsection{A general min-max formula and a pathwise characterization}
We next drop Freidlin's condition and derive an exact min-max representation of the stopping value, valid for every $c\in\mathbb{R}$.
As a consequence, $\mathcal Z(t)$ is characterized by the existence of an admissible trajectory along which every action is nonpositive.
Throughout this and the next subsection, every trajectory defined on $[0,t]$ with terminal value $\mathbf{0}$ is identified with its zero extension to $[0,\infty)$.

For each fixed $t\in \mathbb{R}^{+}$ and $\gamma\in\mathbb{X}$, define
\begin{align*}
  A_{\gamma}(\sigma)=\int_{0}^{\sigma} L(\dot\gamma(s),\gamma(s),t-s) \,ds, \quad \sigma\in[0,t].
\end{align*}
Since the integrand belongs to $L^{1}(0,t)$, $A_{\gamma}\in AC([0,t])$. Set
\begin{align}
\label{eq:St}
  \mathcal{S}(t)=\left\{\mathbf{x}\in\mathbb{R}^{N} \mid \exists\gamma \in\mathbb{X}^{\mathbf{x},t} \;\ \text{s.t.}\;\ A_\gamma(\sigma)\leqslant 0, \;\forall\sigma\in[0,t] \right\}.
\end{align}

We prove Proposition \ref{prop:stopping-value-minmax} in a form that also identifies the zero set of the stopping value.

\begin{proposition}
\label{prop:inf-action-characterization}
Fix $(\mathbf{x},t)\in\mathbb{R}^{N}\times\mathbb{R}^{+}$, let $\mathcal{S}(t)$ be given by \eqref{eq:St}. Then
\begin{align}
\label{eq:stopping-value-inf-action}
  W(\mathbf{x},t) =\min_{\gamma\in\mathbb{X}^{\mathbf{x},t}}\max_{\sigma\in[0,t]}A_{\gamma}(\sigma).
\end{align}
Consequently,
\begin{align}
\label{eq:zero-set-safe-action}
  \mathcal{Z}(t)=\mathcal{S}(t).
\end{align}
\end{proposition}

\begin{proof}
For $\gamma\in\mathbb{X}^{\mathbf{x},t}$, set $M(\gamma)=\max_{\sigma\in[0,t]}A_{\gamma}(\sigma)$. We first show that
\begin{align}
\label{eq:stopping-to-inf-M}
  W(\mathbf{x},t)=\inf_{\gamma\in\mathbb{X}^{\mathbf{x},t}}M(\gamma).
\end{align}
For every $\vartheta\in\Theta$ and every $\gamma\in\mathbb{X}^{\mathbf{x},t}$,
\begin{align*}
  \mathcal P(\mathbf{x},t;\vartheta,\gamma)=A_{\gamma}\left(t\wedge\vartheta[\gamma]\right) \leqslant M(\gamma).
\end{align*}
Taking first the infimum over $\gamma$ and then the supremum over $\vartheta$ implies
\begin{align*}
  W(\mathbf{x},t) \leqslant \inf_{\gamma\in\mathbb{X}^{\mathbf{x},t}}M(\gamma).
\end{align*}

For the converse inequality, let $V$ denote the right-hand side of \eqref{eq:stopping-to-inf-M}.
Since $A_{\gamma}(0)=0$, we have $M(\gamma)\geqslant0$ for every admissible $\gamma$, and hence $V\geqslant0$.
The straight admissible trajectory has finite action, so $V<\infty$. Thus, $V\in[0,\infty)$. Suppose first that $V>0$.
For $m\in(0,V)$, define $\vartheta_{m}:\mathbb{X}\to[0,\infty]$ by
\begin{align}
\label{eq:stopping-the-m}
  \vartheta_{m}[\gamma] =
  \begin{cases}
    \displaystyle \inf \left\{\sigma\in[0,t]\mid A_{\gamma}(\sigma)\geqslant m \right\}, &\text{if } \left\{\sigma\in[0,t]\mid A_{\gamma}(\sigma)\geqslant m \right\} \neq\emptyset,\\
    t, &\text{if } \left\{\sigma\in[0,t]\mid A_{\gamma}(\sigma)\geqslant m \right\} =\emptyset.
  \end{cases}
\end{align}
By construction, $\vartheta_m[\gamma]\in[0,t]$ for every $\gamma\in\mathbb{X}$.
We verify that $\vartheta_m\in\Theta$.
Let $\gamma,\widehat\gamma\in\mathbb{X}$ agree on $[0,a]$, and suppose that $\vartheta_m[\gamma]\leqslant a$. For $\eta\in\mathbb{X}$, set
\begin{align*}
  E_{\eta}=\left\{\sigma\in[0,t]\mid A_{\eta}(\sigma)\geqslant m \right\}.
\end{align*}
If $a<t$, then $\vartheta_m[\gamma]\leqslant a<t$ implies $E_{\gamma}\neq\emptyset$.
Since $A_{\gamma}$ is continuous, $E_{\gamma}$ is compact, and hence
\begin{align*}
  \vartheta_m[\gamma]=\min E_{\gamma}=\min\left(E_{\gamma}\cap[0,a]\right).
\end{align*}
Moreover, $A_{\gamma}=A_{\widehat\gamma}$ on $[0,a]$, so
\begin{align*}
  E_{\gamma}\cap[0,a]=E_{\widehat\gamma}\cap[0,a] \neq\emptyset.
\end{align*}
Thus, by the continuity of $A_{\widehat\gamma}$, $E_{\widehat\gamma}$ is nonempty and compact, with $\min E_{\widehat\gamma}\leqslant a$. Therefore,
\begin{align*}
  \vartheta_m[\widehat\gamma]=\min E_{\widehat\gamma}=\min\left(E_{\widehat\gamma}\cap[0,a]\right)
 =\min\left(E_{\gamma}\cap[0,a]\right)=\vartheta_m[\gamma].
\end{align*}
If $a\geqslant t$, then $A_{\gamma}=A_{\widehat\gamma}$ on $[0,t]$, and hence $E_{\gamma}=E_{\widehat\gamma}$.
Definition \eqref{eq:stopping-the-m} gives $\vartheta_m[\gamma]=\vartheta_m[\widehat\gamma]$.
Therefore, $\vartheta_m\in\Theta$.

For every $\gamma\in\mathbb{X}^{\mathbf{x},t}$ and $m\in(0,V)$, it holds that
\begin{align*}
  A_{\gamma}(0)=0<m<V\leqslant M(\gamma).
\end{align*}
By the continuity of $A_{\gamma}$ and the definition of $\vartheta_m$, we have
\begin{align*}
  \mathcal P(\mathbf{x},t;\vartheta_m,\gamma)=A_{\gamma}\left(t\wedge\vartheta_m[\gamma]\right)=A_{\gamma}\left(\vartheta_m[\gamma]\right)=m.
\end{align*}
Since $\vartheta_m\in\Theta$, it follows that
\begin{align*}
  W(\mathbf{x},t) \geqslant \inf_{\gamma\in\mathbb{X}^{\mathbf{x},t}} \mathcal P(\mathbf{x},t;\vartheta_m,\gamma)=m.
\end{align*}
Letting $m\uparrow V$ proves $W(\mathbf{x},t)\geqslant V$ when $V>0$.
If $V=0$, the admissible stopping time $\vartheta\equiv0$ gives $W(\mathbf{x},t)\geqslant0=V$.
Thus, we obtain \eqref{eq:stopping-to-inf-M}.

It remains to show that the infimum is attained.
Let $\{\gamma_n\}_{n}\subset\mathbb{X}^{\mathbf{x},t}$ be a minimizing sequence.
Since $M(\gamma_n)\to V$ and $M(\gamma_n)\geqslant A_{\gamma_n}(t)$, we have
\begin{align*}
  \frac{1}{4} \int_{0}^{t}|\dot\gamma_n(s)|^{2}ds \leqslant M(\gamma_n)+r_{2}t.
\end{align*}
Together with $\gamma_n(t)=\mathbf{0}$ and the Poincar\'e inequality, this shows that $\{\gamma_n\}_{n}$ is bounded in $H^{1}([0,t];\mathbb{R}^{N})$.
After passing to a subsequence, there exists $\gamma^{0}\in H^{1}([0,t];\mathbb{R}^{N})$ such that
\begin{align*}
  \gamma_n \rightharpoonup \gamma^{0} \quad\text{in }H^{1}([0,t];\mathbb{R}^{N}), \quad  \gamma_n \longrightarrow \gamma^{0} \quad\text{uniformly on }[0,t].
\end{align*}
The endpoint conditions pass to the limit, so $\gamma^{0}(0)=\mathbf{x}$ and $\gamma^{0}(t)=\mathbf{0}$.
Let $\gamma$ be the zero extension of $\gamma^{0}$ to $[0,\infty)$. Then $\gamma\in\mathbb{X}^{\mathbf{x},t}$.


For each fixed $\sigma\in(0,t]$, restriction to $[0,\sigma]$ preserves the weak convergence of the derivatives and the uniform convergence of the trajectories.
Applying the same weak lower-semicontinuity and Fatou's Lemma used to prove \eqref{eq:lsc-inequality}, now on $[0,\sigma]$, yields
\begin{align*}
  A_{\gamma}(\sigma) \leqslant \liminf_{n\to\infty}A_{\gamma_n}(\sigma).
\end{align*}
The same inequality is trivial for $\sigma=0$.
Choose $\sigma_{*}\in[0,t]$ such that $M(\gamma)=A_{\gamma}(\sigma_{*})$. Then
\begin{align*}
  M(\gamma)=A_{\gamma}(\sigma_{*}) \leqslant \liminf_{n\to\infty}A_{\gamma_n}(\sigma_{*}) \leqslant \liminf_{n\to\infty}M(\gamma_n)=V.
\end{align*}
Thus $\gamma$ attains the infimum in \eqref{eq:stopping-value-inf-action}.

Finally, $A_{\gamma}(0)=0$ implies $M(\gamma)\geqslant0$.
The minimum representation further implies that
\begin{align*}
  \mathbf{x}\in\mathcal{Z}(t) \iff \min_{\gamma}M(\gamma)=0 &\iff \exists\gamma \in\mathbb{X}^{\mathbf{x},t} \;\ \text{s.t.}\;\ M(\gamma)=0\\
  &\iff \exists\gamma \in\mathbb{X}^{\mathbf{x},t} \;\ \text{s.t.}\;\ A_\gamma(\sigma)\leqslant 0 \quad\forall\sigma\in[0,t]\\
  &\iff \mathbf{x}\in\mathcal S(t).
\end{align*}
This proves \eqref{eq:zero-set-safe-action}.
\end{proof}

\subsection{Zero set geometry in the non-Freidlin regime}

We now turn to the non-Freidlin regime.
Using the preceding min-max characterization, we derive an equivalent pointwise condition for the zero set geometry and then obtain its convex hull representation.

For $t\in\mathbb{R}^{+}$, define
\begin{align}
\label{eq:Kt}
  \mathcal K(t)=\left\{\mathbf{x}\in\mathbb{R}^{N} \,\middle|\, \exists\gamma \in\mathbb{X}^{\mathbf{x},t} \;\ \text{s.t.}\;\ |\dot\gamma(s)|^{2} \leqslant 4R(\gamma(s),t-s)
  \quad\text{for a.e. }s\in(0,t) \right\}.
\end{align}

\begin{theorem}[Zero set and convex hull characterization]
\label{thm:zero-set-slow-interface}
Assume that $-2\sqrt{r_{2}}<c<2\alpha$.
Then, for every $t\in \mathbb{R}^{+}$,
\begin{align}
\label{eq:zero-set-pointwise-reduction}
  \mathcal{Z}(t)=\mathcal K(t).
\end{align}
Moreover, setting
\begin{align*}
  D_{1}(t)
 =\left\{\mathbf{x}\in\mathbb{R}^{N}\mid |\mathbf{x}|\leqslant2\sqrt{r_{1}}t \right\},
  \quad D_{2}(t)=\left\{\mathbf{x}\in\mathbb{R}^{N}\mid |\mathbf{x}|\leqslant2\sqrt{r_{2}}t,\quad \mathbf{x}\cdot\mathbf{e}\geqslant ct \right\},
\end{align*}
we have
\begin{align}
  \mathcal{Z}(t)=\mathrm{co} \left(D_{1}(t)\cup D_{2}(t) \right).
\label{eq:zero-set-convex-hull}
\end{align}
\end{theorem}

\subsubsection{A pointwise characterization of the zero set}

\begin{proposition}
\label{prop:safe-action-pointwise-char}
Assume that $-2\sqrt{r_{2}}<c<2\alpha$.
Let $\mathcal{S}(\cdot)$ be given by \eqref{eq:St}. Then
\begin{align}
\label{eq:pathwise-to-pointwise}
  \mathcal{S}(t)=\mathcal K(t) \quad\text{for every }t\in \mathbb{R}^{+}.
\end{align}
\end{proposition}

The inclusion $\mathcal K(t)\subset\mathcal S(t)$ follows directly from the definitions.
To prove the reverse inclusion, we need the following sharp estimate at the first contact with the shifting interface.

\begin{lemma}
\label{lem:fixed-contact-optimization}
Assume that $-2\sqrt{r_{2}}<c<2\alpha$.
For fixed $(\mathbf{x},t)\in \mathbb{R}^{N}\times\mathbb{R}^{+}$, denote $\mathbf{x}=\xi\mathbf{e}+{\boldsymbol\zeta}$, where $\xi=\mathbf{x}\cdot\mathbf{e}$ and ${\boldsymbol\zeta}\perp\mathbf{e}$, and set $\rho=|{\boldsymbol\zeta}|$.
Assume that $\xi<ct$, and define
\begin{align*}
  \delta=ct-\xi, \qquad \mu_{c}=r_{2}-\frac{c^{2}}{4}.
\end{align*}
Let $\gamma\in\mathbb{X}^{\mathbf{x},t}$ satisfy $A_\gamma(\sigma)\leqslant0$ for all $\sigma\in[0,t]$, and define its first contact time by
\begin{align}\label{eq:contact-time}
  \tau:=\min\left\{s\in[0,t]\mid \gamma(s)\cdot\mathbf{e}=c(t-s) \right\}.
\end{align}
Assume that $\tau<t$.
Then $\tau\in(0,t)$ and, with $B(\tau)=4r_{1}\tau^{2}-(\delta-c\tau)^{2}$, we have $B(\tau)\geqslant0$ and
\begin{align}
\label{eq:fixed-contact-radius-bound}
  \rho \leqslant \rho_{\tau}(\xi,t):=\sqrt{B(\tau)}+2(t-\tau)\sqrt{\mu_{c}}.
\end{align}
Moreover, the bound in \eqref{eq:fixed-contact-radius-bound} is attained: for every $\tau\in(0,t)$ such that $B(\tau)\geqslant0$ and every unit vector $\mathbf{n}\perp\mathbf{e}$, the point $\mathbf{x}_{\tau} =\xi\mathbf{e}+\rho_{\tau}(\xi,t)\mathbf{n}$ admits a trajectory $\widehat\gamma_{\tau}\in\mathbb{X}^{\mathbf{x}_{\tau},t}$ whose first contact time is $\tau$ and which satisfies
\begin{align}
\label{eq:fixed-contact-pointwise-equality}
  \frac{1}{4} |\dot{\widehat\gamma}_{\tau}(s)|^{2}=R(\widehat\gamma_{\tau}(s),t-s) \quad\text{for a.e. }s\in(0,t).
\end{align}
\end{lemma}

\begin{proof}
Since $\gamma(t)=\mathbf{0}$ and $\gamma$ is continuous, the minimum defining $\tau$ is attained.
Moreover, $\xi<ct$ implies $\tau>0$.
Together with the assumption $\tau<t$, we have $\tau\in(0,t)$.

The parameter range implies $\mu_{c}>0$, while $\xi<ct$ gives $\delta>0$.
By the definition of $\tau$, $\gamma(\tau)\cdot\mathbf{e}=c(t-\tau)$.
Hence, denoting the $\mathbf{e}^{\perp}$-component of the contact point $\gamma(\tau)$ by $\mathbf{p}$, we write
\begin{align}\label{eq:gamm-deco}
  \gamma(\tau)=\mathbf{p}+c(t-\tau)\mathbf{e}, \qquad \mathbf{p}\cdot\mathbf{e}=0.
\end{align}
Let $\tilde{\mathbf{p}}$ be the point of $[\mathbf{0},{\boldsymbol\zeta}]$ closest to $\mathbf{p}$, namely, $|\tilde{\mathbf{p}}-\mathbf{p}| =\min_{\mathbf{q}\in[\mathbf{0},{\boldsymbol\zeta}]}|\mathbf{q}-\mathbf{p}|$.
Then $\tilde{\mathbf{p}}=\eta{\boldsymbol\zeta}$ for some $\eta\in[0,1]$, with the convention $\eta=0$ when ${\boldsymbol\zeta}=\mathbf{0}$, and
\begin{align*}
  |\tilde{\mathbf{p}}|\leqslant|\mathbf{p}|, \qquad |{\boldsymbol\zeta}-\tilde{\mathbf{p}}| \leqslant |{\boldsymbol\zeta}-\mathbf{p}|.
\end{align*}
Set $\rho_{1}:=|{\boldsymbol\zeta}-\tilde{\mathbf{p}}|=(1-\eta)\rho$ and $\rho_{2}:=|\tilde{\mathbf{p}}|=\eta\rho$.
We complete the proof in three steps.

\textbf{Step 1.} We derive the constraints on $\rho_{1}$ and $\rho_{2}$.
Since $\tau$ is the first contact time, $\gamma$ lies in the $r_{1}$-region on $[0,\tau)$.
Applying the Cauchy-Schwarz inequality on $[0,\tau]$ and $[\tau,t]$, using $R\leqslant r_{2}$ on the latter interval, and invoking \eqref{eq:gamm-deco}, we obtain
\begin{align*}
  \frac{\rho_{1}^{2}+(\delta-c\tau)^{2}}{4\tau} -r_{1}\tau \leqslant A_{\gamma}(\tau) \leqslant0
\end{align*}
and
\begin{align*}
  \frac{\rho_{1}^{2}+(\delta-c\tau)^{2}}{4\tau} -r_{1}\tau+\frac{\rho_{2}^{2}}{4(t-\tau)} -\mu_{c}(t-\tau) \leqslant A_{\gamma}(t) \leqslant0.
\end{align*}
Consequently, $B(\tau)\geqslant0$ and, with the notation $D(\tau)=\frac{B(\tau)}{4\tau}+\mu_{c}(t-\tau)$, we have
\begin{align*}
  \rho_{1}^{2} \leqslant B(\tau), \qquad \rho_{2}^{2} \leqslant 4(t-\tau) \left(D(\tau)-\frac{\rho_{1}^{2}}{4\tau} \right).
\end{align*}

\textbf{Step 2.} We prove \eqref{eq:fixed-contact-radius-bound}.
Noting that $\rho=\rho_{1}+\rho_{2}$, we have $\rho\leqslant G_{\tau}(\rho_{1})$, where
\begin{align*}
  G_{\tau}(q)=q+ 2\sqrt{ (t-\tau) \left(D(\tau)-\frac{q^{2}}{4\tau} \right) }, \quad q\in\left[0,\sqrt{B(\tau)}\right] .
\end{align*}
Since $D(\tau)-B(\tau)/(4\tau)=\mu_{c}(t-\tau)>0$, a direct differentiation shows that $G_{\tau}'(q)>0$ whenever $q\in[0,q_{0}(\tau))$, where $q_{0}(\tau):=2\tau\sqrt{\frac{D(\tau)}{t}}$.
Moreover, a direct calculation implies that
\begin{align*}
  B(\tau)-q_{0}(\tau)^{2}=\frac{t-\tau}{t}\Phi(\tau), \quad \Phi(\tau)=2c\delta\tau-\delta^{2} -4(r_{2}-r_{1})\tau^{2}.
\end{align*}
If $c\leqslant0$, then $\Phi(\tau)<0$.
If $c>0$, the assumption $c<2\alpha$ implies
\begin{align*}
  \Phi(\tau)=-4(r_{2}-r_{1}) \left(\tau-\frac{c\delta}{4(r_{2}-r_{1})} \right)^{2}+\delta^{2} \left(\frac{c^{2}}{4(r_{2}-r_{1})}-1 \right) <0.
\end{align*}
Thus $q_{0}(\tau)>\sqrt{B(\tau)}$, and $G_{\tau}'(q)>0$ on $[0,\sqrt{B(\tau)}]$, and hence $G_{\tau}$ is strictly increasing on this interval. Therefore,
\begin{align*}
  \rho \leqslant G_{\tau}\left(\sqrt{B(\tau)}\right)=\sqrt{B(\tau)}+2(t-\tau)\sqrt{\mu_{c}}=\rho_{\tau}(\xi,t).
\end{align*}

\textbf{Step 3.} We construct a trajectory attaining the upper bound.
Fix $\tau\in(0,t)$ with $B(\tau)\geqslant0$, and set $\rho_{1,\tau}=\sqrt{B(\tau)}$ and $\rho_{2,\tau}=2(t-\tau)\sqrt{\mu_{c}}$.
Recall that $\mathbf{x}_{\tau} =\xi\mathbf{e} +(\rho_{1,\tau}+\rho_{2,\tau})\mathbf{n}$. Define
\begin{align*}
  \mathbf{q}_{\tau}=\rho_{2,\tau}\mathbf{n}+c(t-\tau)\mathbf{e}.
\end{align*}
Then $\mathbf{q}_{\tau}\cdot\mathbf{e}=c(t-\tau)$, and so $\mathbf{q}_{\tau}$ lies on the shifting interface. We then define
\begin{align*}
  \widehat\gamma_{\tau}(s) =
  \begin{cases}
    \left(1-\frac{s}{\tau} \right) \mathbf{x}_{\tau}+\frac{s}{\tau}\mathbf{q}_{\tau}, &\text{if }s\in[0,\tau),\\
    \frac{t-s}{t-\tau}\rho_{2,\tau}\mathbf{n}+c(t-s)\mathbf{e}, &\text{if } s\in[\tau,t].
  \end{cases}
\end{align*}
The two expressions match continuously at $s=\tau$ and satisfy $\widehat\gamma_{\tau}(0)=\mathbf{x}_{\tau}$ and $\widehat\gamma_{\tau}(t)=\mathbf{0}$.
Hence $\widehat\gamma_{\tau}\in\mathbb{X}^{\mathbf{x}_{\tau},t}$.
For $0\leqslant s<\tau$, we have
\begin{align*}
  \widehat\gamma_{\tau}(s)\cdot\mathbf{e}-c(t-s) =-\delta\left(1-\frac{s}{\tau}\right)<0,
\end{align*}
and therefore, its first segment lies in the $r_{1}$-region.
Besides, by construction, the second segment lies on the shifting interface.
Thus $\tau$ is the first contact time of $\widehat\gamma_{\tau}$. Finally, we have
\begin{align*}
  \frac{1}{4} |\dot{\widehat\gamma}_{\tau}(s)|^{2} =
  \begin{cases}
    \frac{ \rho_{1,\tau}^{2}+(\delta-c\tau)^{2} }{4\tau^{2}}=r_{1}, &\text{for a.e. }s\in(0,\tau),\\
    \frac{\rho_{2,\tau}^{2}}{4(t-\tau)^{2}}+\frac{c^{2}}{4}=r_{2}, &\text{for a.e. }s\in(\tau,t).
  \end{cases}
\end{align*}
Thus \eqref{eq:fixed-contact-pointwise-equality} holds, and hence $A_{\widehat\gamma_{\tau}}(\sigma)=0$ for every $\sigma\in[0,t]$. Moreover,
\begin{align*}
  |\mathbf{x}_{\tau}-\xi\mathbf{e}|=\rho_{1,\tau}+\rho_{2,\tau}=\rho_{\tau}(\xi,t).
\end{align*}
Therefore, the upper bound is attained, and the proof is complete.
\end{proof}

\begin{proof}[Proof of Proposition \ref{prop:safe-action-pointwise-char}]
On one hand, let $\mathbf{x}\in\mathcal K(t)$ and choose $\gamma\in\mathbb{X}^{\mathbf{x},t}$ satisfying the pointwise inequality in \eqref{eq:Kt}.
Integrating this inequality over $[0,\sigma]$ for each $\sigma\in[0,t]$ shows that $A_{\gamma}(\sigma)\leqslant0$.
Hence $\mathbf{x}\in\mathcal S(t)$, and therefore $\mathcal K(t)\subset\mathcal S(t)$.

On the other hand, let $\mathbf{x}\in\mathcal S(t)$ and choose $\gamma\in\mathbb{X}^{\mathbf{x},t}$ such that
\begin{align*}
  A_{\gamma}(\sigma)\leqslant0 \quad\text{for every }\sigma\in[0,t].
\end{align*}
As before, write $\mathbf{x}=\xi\mathbf{e}+{\boldsymbol\zeta}$, where $\xi=\mathbf{x}\cdot\mathbf{e}$, ${\boldsymbol\zeta}\perp\mathbf{e}$, and $\rho=|{\boldsymbol\zeta}|$.

\textbf{Case 1:} $\xi\geqslant ct$.
Since $R\leqslant r_{2}$, the Cauchy-Schwarz inequality implies
\begin{align*}
  0 \geqslant A_{\gamma}(t) \geqslant \frac{|\mathbf{x}|^{2}}{4t}-r_{2}t.
\end{align*}
Thus $|\mathbf{x}|\leqslant2\sqrt{r_{2}}t$. The straight trajectory $\overline\gamma(s)=\left(1-\frac{s}{t} \right)\mathbf{x}$, $s\in[0,t]$, satisfies
\begin{align*}
  \overline\gamma(s)\cdot\mathbf{e} \geqslant c(t-s) \quad \text{and} \quad |\dot{\overline\gamma}(s)|^{2} \leqslant 4r_{2}.
\end{align*}
Hence $\mathbf{x}\in\mathcal K(t)$.

\textbf{Case 2:} $\xi<ct$.
Define the first contact time of $\gamma$ as:
\begin{align*}
  \tau=\min \left\{s\in[0,t] \mid \gamma(s)\cdot\mathbf{e}=c(t-s) \right\}.
\end{align*}
Since $h(s)=\gamma(s)\cdot\mathbf{e}-c(t-s)$ is continuous, $h(0)<0$, and $h(t)=0$, the above minimum is well-defined and $\tau\in(0,t]$.
If $\tau=t$, by the continuity of $h(\cdot)$ and the definition of the first contact time $\tau$, $\gamma$ lies in the $r_{1}$-region for a.e. $s\in(0,t)$, and therefore
\begin{align*}
  0 \geqslant A_{\gamma}(t) \geqslant \frac{|\mathbf{x}|^{2}}{4t}-r_{1}t.
\end{align*}
It follows that $|\mathbf{x}|\leqslant2\sqrt{r_{1}}t$.
The straight trajectory from $\mathbf{x}$ to $\mathbf{0}$ remains in the $r_{1}$-region and has a speed at most $2\sqrt{r_{1}}$.
Hence $\mathbf{x}\in\mathcal K(t)$.

It remains to consider $\tau<t$.
Lemma \ref{lem:fixed-contact-optimization} first implies
\begin{align*}
  B(\tau)\geqslant0, \qquad \rho \leqslant \rho_{\tau}(\xi,t):=\sqrt{B(\tau)}+2(t-\tau)\sqrt{\mu_{c}}.
\end{align*}
Since $\tau<t$ and $\mu_{c}>0$, $\rho_{\tau}(\xi,t)>0$.
Choose $\mathbf{n}=\frac{{\boldsymbol\zeta}}{\rho}$ if $\rho>0$, and choose any unit vector $\mathbf{n}\perp\mathbf{e}$ if $\rho=0$. Set
\begin{align*}
  \mathbf{x}_{\tau}=\xi\mathbf{e}+\rho_{\tau}(\xi,t)\mathbf{n}.
\end{align*}
Since $B(\tau)\geqslant0$, Lemma \ref{lem:fixed-contact-optimization} further concludes the existence of a trajectory $\widehat\gamma_{\tau}\in \mathbb{X}^{\mathbf{x}_{\tau},t}$ satisfying
\begin{align}
\label{eq:sharp-gam-ta}
  \frac{1}{4} |\dot{\widehat\gamma}_{\tau}(s)|^{2}=R(\widehat\gamma_{\tau}(s),t-s) \quad\text{for a.e. }s\in(0,t).
\end{align}
Set $\theta=\frac{\rho}{\rho_{\tau}(\xi,t)}\in[0,1]$ and define
\begin{align*}
  \gamma_{\theta}(s)=\left(\widehat\gamma_{\tau}(s)\cdot\mathbf{e} \right)\mathbf{e}+\theta \left[ \widehat\gamma_{\tau}(s)-\left(\widehat\gamma_{\tau}(s)\cdot\mathbf{e} \right)\mathbf{e} \right].
\end{align*}
By the choices of $\mathbf{n}$ and $\theta$, it holds that $\theta\rho_{\tau}(\xi,t)\mathbf{n}=\rho\mathbf{n}={\boldsymbol\zeta}$, and hence, we have
\begin{align*}
  \gamma_{\theta}(0)=\xi\mathbf{e}+{\boldsymbol\zeta}=\mathbf{x}, \quad \text{and} \quad \gamma_{\theta}(t)=\mathbf{0}.
\end{align*}
Since the above transformation is linear, $\gamma_{\theta}\in H^{1}([0,t];\mathbb{R}^{N})$.
Together with the endpoint conditions and the zero-extension convention, this implies $\gamma_{\theta}\in\mathbb{X}^{\mathbf{x},t}$.

By the definition of $\gamma_{\theta}$,
\begin{align*}
  \gamma_{\theta}(s)\cdot\mathbf{e}=\widehat\gamma_{\tau}(s)\cdot\mathbf{e}, \quad \forall  s\in[0,t].
\end{align*}
Notice that $R(\mathbf{y},t-s)$ depends on $\mathbf{y}$ only through $\mathbf{y}\cdot\mathbf{e}$. Therefore,
\begin{align}
\label{eq:equal-potential-two-traj}
  R(\gamma_{\theta}(s),t-s)=R(\widehat\gamma_{\tau}(s),t-s), \quad \forall  s\in[0,t].
\end{align}
For a.e. $s\in(0,t)$, set
\begin{align*}
  \dot{\widehat\gamma}_{\tau,\perp}(s)=\dot{\widehat\gamma}_{\tau}(s)-\left(\dot{\widehat\gamma}_{\tau}(s)\cdot\mathbf{e} \right)\mathbf{e}.
\end{align*}
Differentiating the definition of $\gamma_{\theta}$ gives
\begin{align*}
  \dot\gamma_{\theta}(s)=\left(\dot{\widehat\gamma}_{\tau}(s)\cdot\mathbf{e} \right)\mathbf{e}+\theta \dot{\widehat\gamma}_{\tau,\perp}(s).
\end{align*}
Since the two terms on the right-hand side are orthogonal, we have
\begin{align*}
  |\dot\gamma_{\theta}(s)|^{2}=|\dot{\widehat\gamma}_{\tau}(s)|^{2}-(1-\theta^{2}) |\dot{\widehat\gamma}_{\tau,\perp}(s)|^{2}.
\end{align*}
Combining this with \eqref{eq:sharp-gam-ta} and \eqref{eq:equal-potential-two-traj}, we obtain
\begin{align*}
  \frac{1}{4}|\dot\gamma_{\theta}(s)|^{2}-R(\gamma_{\theta}(s),t-s)=-\frac{1-\theta^{2}}{4} |\dot{\widehat\gamma}_{\tau,\perp}(s)|^{2} \leqslant0
  \quad \text{for a.e. } s\in(0,t).
\end{align*}
Thus $\gamma_{\theta}$ verifies that $\mathbf{x}\in\mathcal K(t)$.
To sum up, in all cases, we have $\mathbf{x}\in \mathcal{K}(t)$, and hence $\mathcal{S}(t)\subset \mathcal{K}(t)$.
This completes the proof of \eqref{eq:pathwise-to-pointwise}.
\end{proof}

\subsubsection{Convex hull characterization}

By \eqref{eq:zero-set-safe-action} and Proposition \ref{prop:safe-action-pointwise-char}, $\mathcal Z(t)=\mathcal K(t)$.
It remains to identify the geometry of $\mathcal K(t)$.

\begin{proposition}[Convex hull geometry]
\label{prop:pointwise-convex-hull}
Assume that $-2\sqrt{r_{2}}<c<2\alpha$. Then
\begin{align}
\label{eq:pointwise-convex-hull}
  \mathcal K(t)=\mathrm{co} \left(D_{1}(t)\cup D_{2}(t) \right) \quad\text{for every }t\in \mathbb{R}^{+}.
\end{align}
\end{proposition}

\begin{proof}
Define the truncated interface by
\begin{align*}
  C(t) &= D_{2}(t) \cap \left\{\mathbf{y}\in\mathbb{R}^{N}\mid \mathbf{y}\cdot\mathbf{e}=ct \right\}
 =\left\{ct\mathbf{e}+{\boldsymbol\zeta}\mid {\boldsymbol\zeta}\cdot\mathbf{e}=0,\; |{\boldsymbol\zeta}| \leqslant t\sqrt{4r_{2}-c^{2}} \right\}.
\end{align*}

\textbf{Step 1.} We prove
\begin{align}
\label{eq:pointwise-cup-convex}
  \mathcal K(t)=D_{2}(t) \cup \mathrm{co} \left(D_{1}(t)\cup C(t) \right).
\end{align}

\textbf{(a).}
We first prove the ``$\subset$'' part.
Let $\mathbf{x}\in\mathcal K(t)$ and choose $\gamma\in\mathbb{X}^{\mathbf{x},t}$ satisfying the pointwise bound.
Since $R\leqslant r_{2}$, the Cauchy-Schwarz inequality yields $|\mathbf{x}|\leqslant2\sqrt{r_{2}}t$.

\textbf{Case (a1):} $\mathbf{x}\cdot\mathbf{e}\geqslant ct$.
Then it is trivial that $\mathbf{x}\in D_{2}(t)$.

\textbf{Case (a2):} $\mathbf{x}\cdot\mathbf{e}<ct$.
As in Case 2 of the proof of Proposition \ref{prop:safe-action-pointwise-char}, define the first contact time:
\begin{align*}
  \tau=\min \left\{s\in[0,t]\mid \gamma(s)\cdot\mathbf{e}=c(t-s) \right\}.
\end{align*}
The same continuity argument gives $\tau\in(0,t]$.

If $\tau=t$, then $\gamma$ remains in the $r_{1}$ region for a.e. $s\in(0,t)$.
Hence the pointwise bound implies $|\mathbf{x}|\leqslant2\sqrt{r_{1}}t$.
Therefore, $\mathbf{x}\in D_{1}(t)$.
It remains to consider $\tau\in(0,t)$. Set
\begin{align}
\label{eq:convex-para}
  \kappa=\frac{\tau}{t}, \qquad Q=\frac{t}{\tau} \left(\mathbf{x}-\gamma(\tau) \right), \qquad P=\frac{t}{t-\tau}\gamma(\tau).
\end{align}
Before time $\tau$, the definition of $\tau$ implies that the trajectory lies in the $r_{1}$ region. Thus
\begin{align*}
  |Q| \leqslant \frac{t}{\tau} \int_{0}^{\tau}|\dot\gamma(s)|ds \leqslant 2\sqrt{r_{1}}t,
\end{align*}
and hence $Q\in D_{1}(t)$. Moreover, we have
\begin{align*}
  P\cdot\mathbf{e}=ct, \qquad |P| \leqslant \frac{t}{t-\tau} \int_{\tau}^{t}|\dot\gamma(s)|ds \leqslant 2\sqrt{r_{2}}t.
\end{align*}
Therefore, $P\in C(t)$.
Notably, by \eqref{eq:convex-para}, we have $\mathbf{x}=\kappa Q+(1-\kappa)P$, and consequently, $\mathbf{x} \in \mathrm{co} \left(D_{1}(t)\cup C(t) \right)$.
Combining this with the preceding cases, we obtain
\begin{align}
\label{eq:convex-subset}
  \mathcal K(t) \subset D_{2}(t) \cup \mathrm{co} \left(D_{1}(t)\cup C(t) \right).
\end{align}

\textbf{(b).}
We then prove the ``$\supset$'' part.
If $\mathbf{x}\in D_{2}(t)$, then the straight trajectory
\begin{align*}
  \gamma(s)=\left(1-\frac{s}{t} \right)\mathbf{x}, \quad s\in[0,t],
\end{align*}
satisfies
\begin{align*}
  \gamma(s)\cdot\mathbf{e} \geqslant c(t-s), \qquad |\dot\gamma(s)|=\frac{|\mathbf{x}|}{t} \leqslant 2\sqrt{r_{2}}.
\end{align*}
Hence $D_{2}(t)\subset\mathcal K(t)$.
Next, if $\mathbf{x}\in\mathrm{co}(D_{1}(t)\cup C(t))$, since $D_{1}(t)$ and $C(t)$ are convex, it can be written as
\begin{align*}
  \mathbf{x}=\bar\kappa \bar Q+(1-\bar\kappa)\bar P \qquad \text{for some} \quad  \bar Q\in D_{1}(t), \quad \bar P\in C(t), \quad \text{and }\bar \kappa\in[0,1].
\end{align*}

\textbf{Case (b1):} $\mathbf{x}\cdot\mathbf{e}<ct$.
Then $\bar \kappa>0$, and we set $\tau=\bar \kappa t$. If $\bar\kappa<1$, define
\begin{align*}
  \gamma(s) =
  \begin{cases}
    \mathbf{x}-\frac{s}{t}\bar Q, &\text{if }s\in[0,\tau],\\[3pt]
    \frac{t-s}{t}\bar P, &\text{if }s\in[\tau, t].
  \end{cases}
\end{align*}
The two expressions agree at $s=\tau$.
Moreover, a direct computation yields
\begin{align*}
  \gamma(s)\cdot\mathbf{e}-c(t-s)=\left(\frac{s}{\tau}-1 \right) \left(ct-\mathbf{x}\cdot\mathbf{e} \right) <0 \quad \text{for }s\in[0,\tau).
\end{align*}
Thus the first segment lies in the $r_{1}$ region and $\frac{|\bar Q|}{t}\leqslant2\sqrt{r_{1}}$.
The second segment lies on the shifting interface and $\frac{|\bar P|}{t}\leqslant2\sqrt{r_{2}}$.
Therefore, $\gamma$ satisfies the pointwise bound in $\mathcal{K}(t)$.
If $\bar \kappa=1$, then $\mathbf{x}=\bar Q\in D_{1}(t)$.
Since $\mathbf{x}\cdot\mathbf{e}<ct$, the straight trajectory from $\mathbf{x}$ to $\mathbf{0}$ remains in the $r_{1}$ region for $s<t$ and has a speed at most $2\sqrt{r_{1}}$.
Hence $\mathbf{x}\in\mathcal K(t)$.

\textbf{Case (b2):} $\mathbf{x}\cdot\mathbf{e}\geqslant ct$. It follows from the fact
\begin{align*}
  \mathrm{co} \left(D_{1}(t)\cup C(t) \right) \subset \overline {B_{2\sqrt{r_{2}}t}(\mathbf{0})}
\end{align*}
that $\mathbf{x}\in D_{2}(t)\subset\mathcal K(t)$.
Combining the cases above with \eqref{eq:convex-subset} proves \eqref{eq:pointwise-cup-convex}.

\textbf{Step 2.} We eliminate the auxiliary truncated interface $C(t)$.
Since $C(t)\subset D_{2}(t)$, \eqref{eq:pointwise-cup-convex} gives
\begin{align*}
  \mathcal K(t) \subset \mathrm{co} \left(D_{1}(t)\cup D_{2}(t) \right).
\end{align*}
Conversely, for any $\mathbf{x}\in\mathrm{co}(D_{1}(t)\cup D_{2}(t))$, because $D_{1}(t)$ and $D_{2}(t)$ are convex, we can write
\begin{align*}
  \mathbf{x}=\tilde\kappa  \tilde Q+(1- \tilde\kappa) \tilde P, \qquad \text{for some } \tilde Q\in D_{1}(t), \quad \tilde P\in D_{2}(t),
  \quad \text{and }\tilde\kappa\in[0,1].
\end{align*}

\textbf{Case 1:} $\mathbf{x}\cdot\mathbf{e}\geqslant ct$.
Then $\mathbf{x}\in D_{2}(t)$ since $D_{1}(t)\cup D_{2}(t)\subset\overline {B_{2\sqrt{r_{2}}t}(\mathbf{0})}$.
Together with Step 1, we have $\mathbf{x}\in \mathcal{K}(t)$.

\textbf{Case 2:} $\mathbf{x}\cdot\mathbf{e}<ct$.
Then we must have $ \tilde Q\cdot\mathbf{e}<ct$.
Since $\tilde P\in D_{2}(t)$, we also have $\tilde P\cdot\mathbf{e}\geqslant ct$.
Hence, there exists a unique $P_{0}$ on the segment joining $\tilde Q$ to $\tilde P$ satisfying $P_{0}\cdot\mathbf{e}=ct$.
Since the segment is contained in $\overline {B_{2\sqrt{r_{2}}t}(\mathbf{0})}$, we know $P_{0}\in C(t)$.
Moreover, as $\mathbf{x}$ lies on the same segment and $\mathbf{x}\cdot\mathbf{e}<ct=P_0\cdot\mathbf{e}$, it follows that $\mathbf{x}$ lies on the segment joining $ \tilde Q$ to $P_{0}$. Therefore,
\begin{align*}
  \mathbf{x} \in \mathrm{co} \left(D_{1}(t)\cup C(t) \right) \subset \mathcal K(t).
\end{align*}
Combining these cases, we obtain \eqref{eq:pointwise-convex-hull}.
\end{proof}

\begin{proof}[Proof of Theorem \ref{thm:zero-set-slow-interface}]
Proposition \ref{prop:inf-action-characterization} together with Proposition \ref{prop:safe-action-pointwise-char} implies \eqref{eq:zero-set-pointwise-reduction}.
Finally, Proposition \ref{prop:pointwise-convex-hull} yields \eqref{eq:zero-set-convex-hull}.
\end{proof}

\begin{proof}[Proof of Theorem \ref{thm:main-freidlin-dichotomy}]
Assertion (i) is Theorem \ref{thm:backward-W-is-viscosity-solution}.
For (ii), Theorem \ref{thm:freidlin-sharp} gives Freidlin's condition, so (i) yields \eqref{eq:main-freidlin-reduction}.
Lemma \ref{lem:freidlin-zero-level-characterization}, together with \eqref{eq:wulff-shape-definition}, then gives the zero-set identity.

For assertion (iii), Proposition \ref{prop:freidlin-failure-full-range} gives the failure of Freidlin's condition in the intermediate regime.
Theorem \ref{thm:zero-set-slow-interface}, evaluated at $t=1$ and combined with \eqref{eq:wulff-shape-definition}, gives the convex hull formula.
Taking the constant stopping times $\vartheta\equiv0$ and $\vartheta\equiv1$ in \eqref{eq:stopping-game-value} yields
\begin{align*}
  W(\mathbf{x},1)\geqslant\max\{0,J(\mathbf{x},1)\},
\end{align*}
and hence $\mathcal{W}\subseteq\{\mathbf{x}\in\mathbb{R}^{N}\mid J(\mathbf{x},1)\leqslant0\}$. It remains to prove the strict inclusion.

Choose $\mathbf{x}_{0}$ and its unconstrained minimizing trajectory $\gamma_{*}$ as in the proof of Proposition \ref{prop:freidlin-failure-full-range}, with $t=1$.
Then $J(\mathbf{x}_{0},1)=0$.
If $W(\mathbf{x}_{0},1)=0$, Proposition \ref{prop:inf-action-characterization} provides $\gamma\in\mathbb{X}^{\mathbf{x}_{0},1}$ such that
\begin{align*}
  0=J(\mathbf{x}_{0},1)\leqslant A_{\gamma}(1)\leqslant\max_{\sigma\in[0,1]}A_{\gamma}(\sigma)=0.
\end{align*}
Thus $\gamma$ is the unique minimizer $\gamma_{*}$ for the optimal control problem \eqref{eq:backward-unconstrained-value} on $[0,1]$ by Proposition \ref{prop:backward-sharp-value-function}.
However, Bellman's principle along the minimizing trajectory $\gamma_{*}$ and \eqref{eq:freidlin-failure-strict-negativity} give
\begin{align*}
  A_{\gamma_{*}}(s)=J(\mathbf{x}_{0},1)-J(\gamma_{*}(s),1-s)>0\quad \text{for } s\in(0,1),
\end{align*}
contradicting $\max_{\sigma\in[0,1]}A_{\gamma}(\sigma)=0$.
Therefore $W(\mathbf{x}_{0},1)>0=J(\mathbf{x}_{0},1)$, which proves the strict inclusion and completes the proof of (iii).
\end{proof}

\section{Wulff shapes, optimal trajectories and directional spreading speeds}
\label{section:explicit-wulff-shape}

The reductions obtained in Section \ref{section:stopping-game-reductions} characterize the propagation region through either an unconstrained value function or a convex hull representation.
We now determine the Wulff shapes explicitly, characterize the optimal trajectories reaching the front, and derive the directional spreading speeds in all speed regimes.
We retain the notation introduced in Subsection \ref{subsection:explicit-wulff-shapes}.

\subsection{Wulff shapes and optimal trajectories in the Freidlin regimes}

\begin{proposition}
\label{prop:fast-receding-front}
Assume that $c\leqslant-2\sqrt{r_{2}}$.
Then, for every $t\in\mathbb{R}^{+}$,
\begin{align}
\label{eq:fast-receding-region}
  \mathcal{Z}(t)=\left\{\mathbf{x}\in\mathbb{R}^{N} \mid |\mathbf{x}|\leqslant2\sqrt{r_{2}}t \right\}.
\end{align}
Consequently,
\begin{align}
\label{eq:fast-receding-front}
  \partial\mathcal{Z}(t)=\left\{\mathbf{x}\in\mathbb{R}^{N} \mid |\mathbf{x}|=2\sqrt{r_{2}}t \right\}.
\end{align}
\end{proposition}

\begin{proof}
By Lemma \ref{lem:freidlin-zero-level-characterization}, it suffices to determine the set where $J(\cdot,t)\leqslant0$.

Since $R\leqslant r_{2}$, every $\gamma\in\mathbb B^{\mathbf{x},t}$ satisfies
\begin{align*}
  \mathcal I_{\mathrm{b}}(\mathbf{x},t;\gamma) \geqslant \frac{1}{4} \int_{0}^{t}|\dot\gamma(s)|^{2}ds-r_{2}t \geqslant \frac{|\mathbf{x}|^{2}}{4t}-r_{2}t.
\end{align*}
Taking the infimum over $\gamma$, we obtain
\begin{align}
\label{eq:negative-fast-global-lower}
  J(\mathbf{x},t) \geqslant \frac{|\mathbf{x}|^{2}}{4t}-r_{2}t.
\end{align}

Now suppose that $|\mathbf{x}|\leqslant2\sqrt{r_{2}}t$.
Then $\xi\geqslant-|\mathbf{x}|\geqslant-2\sqrt{r_{2}}t\geqslant ct$.
Thus $\mathbf{x}$ lies in the high-growth branch of Proposition \ref{prop:backward-sharp-value-function}, which gives directly
\begin{align*}
  J(\mathbf{x},t)=\frac{|\mathbf{x}|^{2}}{4t}-r_{2}t \leqslant0.
\end{align*}

Conversely, if $|\mathbf{x}|>2\sqrt{r_{2}}t$, then \eqref{eq:negative-fast-global-lower} gives $J(\mathbf{x},t)>0$.
Therefore, $\{J(\cdot,t)\leqslant0\}$ is precisely the ball in \eqref{eq:fast-receding-region}.
The result immediately follows from Lemma \ref{lem:freidlin-zero-level-characterization}.
\end{proof}

\begin{proposition}
\label{prop:directional-front-case1}
Assume that $2\alpha\leqslant c\leqslant2\sqrt{r_{2}}$.
Then, for every $t\in\mathbb{R}^{+}$,
\begin{align}
\label{eq:front-region-decomposition}
  \mathcal{Z}(t)=\Sigma_{1}^{2\lambda}(t) \cup \Sigma_{\mathrm{m}}^{\lambda,c}(t) \cup \Sigma_{2}^{c}(t).
\end{align}
Moreover, let $\Gamma_{1}(t)$, $\Gamma_{\mathrm{m}}(t)$, and $\Gamma_{2}(t)$ be given by
\begin{align*}
  \Gamma_{1}(t) &= \left\{\mathbf{x}\mid \xi^{2}+\rho^{2}=4r_{1}t^{2}, \, \xi\leqslant2\lambda t \right\},\\
  \Gamma_{\mathrm{m}}(t) &= \left\{\mathbf{x}\mid \rho^{2}=4t\left[ (r_{1}+\lambda^{2})t-\lambda\xi \right], \, 2\lambda t\leqslant\xi\leqslant ct \right\},\\
  \Gamma_{2}(t) &= \left\{\mathbf{x}\mid \xi^{2}+\rho^{2}=4r_{2}t^{2}, \, \xi\geqslant ct \right\}.
\end{align*}
Then
\begin{align}
\label{eq:front-decomposition}
  \partial\mathcal{Z}(t)=\Gamma_{1}(t)\cup\Gamma_{\mathrm{m}}(t)\cup\Gamma_{2}(t).
\end{align}
\end{proposition}

\begin{proof}
By Lemma \ref{lem:freidlin-zero-level-characterization}, it suffices to determine the zero sublevel set and the zero level set of $J$.

The assumptions on $c$ imply $0\leqslant\lambda\leqslant\sqrt{r_{2}}-\alpha<\sqrt{r_{1}}$.
Since $r_{1}+\lambda^{2}=c^{2}/4+r_{2}-c\alpha$, Proposition \ref{prop:backward-sharp-value-function} yields
\begin{align}
\label{eq:J-case-1}
  J(\mathbf{x},t)=\frac{\rho^{2}}{4t} +
  \begin{cases}
    \frac{\xi^{2}}{4t}-r_{1}t, &\text{if }\xi\leqslant2\lambda t,\\
    \lambda\xi-(r_{1}+\lambda^{2})t, &\text{if }2\lambda t<\xi<ct,\\
    \frac{\xi^{2}}{4t}-r_{2}t, &\text{if }\xi\geqslant ct.
  \end{cases}
\end{align}
Solving $J(\mathbf{x},t)\leqslant 0$ in the three intervals given by \eqref{eq:J-case-1} implies \eqref{eq:front-region-decomposition}.
At $\xi=2\lambda t$, both the first two expressions give $\rho^{2}=4(r_{1}-\lambda^{2})t^{2}$.
At $\xi=ct$, both the last two expressions give $\rho^{2}=(4r_{2}-c^{2})t^{2}$, where we used the fact that $ r_{1}+\lambda^{2}-c\lambda=r_{2}-\frac{c^{2}}{4}$.
These quantities are nonnegative under the stated assumptions, so the three pieces join at their common endpoints.
Combining this with Lemma \ref{lem:freidlin-zero-level-characterization}, we obtain \eqref{eq:front-decomposition}. The proof is complete.
\end{proof}

\begin{proposition}
\label{prop:intermediate-fast-front}
Assume that $2\sqrt{r_{2}}<c<2(\alpha+\sqrt{r_{1}})$.
Recall that $c_\mathrm{nlp}=\lambda+\frac{r_{1}}{\lambda}$.
Then $2\lambda<c_\mathrm{nlp}<c$.
Moreover, for every $t\in\mathbb{R}^{+}$,
\begin{align}
\label{eq:intermediate-fast-region}
  \mathcal Z(t)=\Sigma_{1}^{2\lambda}(t) \cup \Sigma_{\mathrm{m}}^{\lambda,c_\mathrm{nlp}}(t).
\end{align}
Let $\Gamma_{1}(t)$ and $\Gamma_{\mathrm{m}}(t)$ be given by
\begin{align}
  \Gamma_{1}(t) &= \left\{\mathbf{x}\mid \xi^{2}+\rho^{2}=4r_{1}t^{2}, \, \xi\leqslant2\lambda t \right\}, \nonumber\\
\label{eq:intermediate-fast-middle}
  \Gamma_{\mathrm{m}}(t) &= \left\{\mathbf{x}\mid \rho^{2}=4t\left[ (r_{1}+\lambda^{2})t-\lambda\xi \right], \, 2\lambda t\leqslant\xi\leqslant c_\mathrm{nlp}t \right\}.
\end{align}
Then
\begin{align}
\label{eq:intermediate-fast-front}
  \partial\mathcal Z(t)=\Gamma_{1}(t)\cup\Gamma_{\mathrm{m}}(t).
\end{align}
\end{proposition}

\begin{proof}
The assumptions imply $0<\sqrt{r_{2}}-\alpha<\lambda<\sqrt{r_{1}}$. Consequently,
\begin{align*}
  c_\mathrm{nlp}-2\lambda=\frac{r_{1}-\lambda^{2}}{\lambda} >0, \qquad c-c_\mathrm{nlp}=\frac{c^{2}/4-r_{2}}{\lambda} >0.
\end{align*}
Thus $2\lambda<c_\mathrm{nlp}<c$.

By Lemma \ref{lem:freidlin-zero-level-characterization}, it suffices to determine the zero sublevel set and the zero level set of $J$.
Since $r_{1}+\lambda^{2}=\frac{c^{2}}{4}+r_{2}-c\alpha$, \eqref{eq:backward-J-decomposition} and \eqref{eq:backward-V-piecewise} take the same forms as \eqref{eq:J-case-1}.

On the branch $\xi\leqslant2\lambda t$, the inequality $J(\mathbf{x},t)\leqslant0$ is equivalent to
\begin{align*}
  \xi^{2}+\rho^{2}\leqslant4r_{1}t^{2},
\end{align*}
which gives $\Sigma_{1}^{2\lambda}(t)$.
On the branch $2\lambda t\leqslant\xi\leqslant c_\mathrm{nlp}t$, it is equivalent to
\begin{align*}
  \rho^{2} \leqslant 4t\left[ (r_{1}+\lambda^{2})t-\lambda\xi \right].
\end{align*}
Since $\lambda>0$, the right-hand side is nonnegative if and only if $\xi\leqslant c_\mathrm{nlp}t$.
Together with $2\lambda<c_\mathrm{nlp}<c$, this gives $\Sigma_{\mathrm{m}}^{\lambda,c_\mathrm{nlp}}(t)$.

The high-growth branch contributes no point to the zero sublevel set.
Indeed, if $\xi\geqslant ct$, then
\begin{align*}
  J(\mathbf{x},t) \geqslant \left(\frac{c^{2}}{4}-r_{2} \right)t >0.
\end{align*}
Therefore, \eqref{eq:intermediate-fast-region} holds.

Replacing the two preceding sublevel inequalities with equalities gives $\Gamma_{1}(t)$ and $\Gamma_{\mathrm{m}}(t)$, respectively.
At $\xi=2\lambda t$, the two equations agree and give $\rho^{2}=4(r_{1}-\lambda^{2})t^{2}$, while at $\xi=c_\mathrm{nlp}t$ the intermediate component satisfies $\rho=0$.
Hence Lemma \ref{lem:freidlin-zero-level-characterization} yields \eqref{eq:intermediate-fast-front}.
\end{proof}

\begin{proposition}
\label{prop:extremely-fast-front}
Assume that $c\geqslant2(\alpha+\sqrt{r_{1}})$.
Then, for every $t\in\mathbb{R}^{+}$,
\begin{align}
\label{eq:extremely-fast-region}
  \mathcal{Z}(t)=\left\{\mathbf{x}\in\mathbb{R}^{N} \,\middle|\, |\mathbf{x}|\leqslant2\sqrt{r_{1}}t \right\}.
\end{align}
Consequently,
\begin{align}
\label{eq:extremely-fast-front}
  \partial\mathcal{Z}(t)=\left\{\mathbf{x}\in\mathbb{R}^{N} \,\middle|\, |\mathbf{x}|=2\sqrt{r_{1}}t \right\}.
\end{align}
\end{proposition}

\begin{proof}
Recall that $\lambda=\frac{c}{2}-\alpha$.
Then $\lambda\geqslant\sqrt{r_{1}}$. Moreover, $c>2\sqrt{r_{2}}$.
By Lemma \ref{lem:freidlin-zero-level-characterization}, it suffices to determine the set where $J(\cdot,t)\leqslant0$.

Suppose first that $|\mathbf{x}|\leqslant2\sqrt{r_{1}}t$.
Then $\xi\leqslant|\mathbf{x}| \leqslant2\sqrt{r_{1}}t \leqslant2\lambda t$.
Proposition \ref{prop:backward-sharp-value-function} therefore gives
\begin{align*}
  J(\mathbf{x},t)=\frac{|\mathbf{x}|^{2}}{4t}-r_{1}t \leqslant0.
\end{align*}

Conversely, suppose that $|\mathbf{x}|>2\sqrt{r_{1}}t$.
If $\xi\leqslant2\lambda t$, then
\begin{align*}
  J(\mathbf{x},t)=\frac{|\mathbf{x}|^{2}}{4t}-r_{1}t >0.
\end{align*}
If $2\lambda t<\xi<ct$, then
\begin{align*}
  J(\mathbf{x},t)=\frac{\rho^{2}}{4t} +\lambda\xi-(r_{1}+\lambda^{2})t>\frac{\rho^{2}}{4t} +(\lambda^{2}-r_{1})t \geqslant0.
\end{align*}
Finally, if $\xi\geqslant ct$, then
\begin{align*}
  J(\mathbf{x},t)=\frac{|\mathbf{x}|^{2}}{4t}-r_{2}t \geqslant \left(\frac{c^{2}}{4}-r_{2}\right)t >0.
\end{align*}
Thus $\{J(\cdot,t)\leqslant0\}$ is precisely the ball in \eqref{eq:extremely-fast-region}.
The boundary formula follows immediately.
\end{proof}

\begin{proof}[Proof of Theorem \ref{thm:main-wulff-freidlin}]
By \eqref{eq:wulff-shape-definition}, the Wulff shape is obtained by setting $t=1$.
Assertions (i)--(iv) therefore follow from Propositions \ref{prop:fast-receding-front}--\ref{prop:extremely-fast-front}.
\end{proof}


We next establish uniqueness of the optimal trajectories for the min-max formula \eqref{eq:stopping-value-inf-action} at the propagation front in the Freidlin regimes.

\begin{proposition}
\label{prop:freidlin-front-trajectories}
Assume that $c\leqslant-2\sqrt{r_{2}}$ or $c\geqslant2\alpha$.
For every $t\in \mathbb{R}^{+}$ and $\mathbf{x}\in t\partial\mathcal{W}$, the unconstrained optimal control problem \eqref{eq:backward-unconstrained-value} and the min-max formula \eqref{eq:stopping-value-inf-action} have the same unique minimizing trajectory $\gamma_{*}$ on $[0,t]$.
Its expression is obtained from \eqref{eq:front-straight-optimal-path} in the homogeneous branches and \eqref{eq:freidlin-front-optimal-path} in the intermediate branch by rescaling from time $1$ to time $t$.
It remains on $(t-s)\partial\mathcal{W}$ for $0\leqslant s<t$ in the homogeneous branches and on the critical cylinder.
In the intermediate branch $2\lambda t<\xi<ct$, set $\sigma_*=\frac{ct-\xi}{2\alpha}\in(0,t)$.
Then
\begin{align}
  \label{eq:freidlin-front-value-profile}
  W(\gamma_{*}(s),t-s)=\frac{2\lambda\alpha}{t}
  \begin{cases}
    (t-\sigma_*)s, & s\in[0,\sigma_*),\\
    \sigma_*(t-s), & s\in[\sigma_*,t).
  \end{cases}
\end{align}
\end{proposition}

\begin{proof}
Fix $t\in \mathbb{R}^{+}$ and $\mathbf{x}\in t\partial\mathcal{W}$.
By Theorem~\ref{thm:main-freidlin-dichotomy}, $W=\max\{0,J\}$ and $\mathcal Z(t)=\{J(\cdot,t)\leqslant0\}$.
The continuity of $J$ therefore gives $J(\mathbf{x},t)=0$.
Proposition \ref{prop:backward-sharp-value-function} gives the unique trajectory minimizing the unconstrained optimal control problem.

To compare with the unit-time formulas, we use the scaling $\gamma_t(s)=t\gamma_1(\frac{s}{t})$, which maps trajectories from $\frac{\mathbf{x}}{t}$ to $\mathbf{0}$ on $[0,1]$ bijectively onto trajectories from $\mathbf{x}$ to $\mathbf{0}$ on $[0,t]$.
Since $R(t\mathbf{y},t\tau)=R(\mathbf{y},\tau)$, it multiplies both the total action and the min-max objective by $t$.
In the two homogeneous branches, this trajectory is $\gamma_{*}(s)=(1-\frac{s}{t})\mathbf{x}$, the rescaled form of \eqref{eq:front-straight-optimal-path}.
The self-similarity of the front gives $\gamma_{*}(s)\in(t-s)\partial\mathcal{W}$ for $0\leqslant s<t$.
Theorem \ref{thm:main-wulff-freidlin}-(i) and (iv) show that all front points lie in these branches if $c\geqslant2(\alpha+\sqrt{r_{1}})$ or $c\leqslant-2\sqrt{r_{2}}$.

In the intermediate branch, $\sigma_*=t-\tau_*$, where $\tau_*$ is given by \eqref{eq:tau-star}.
The trajectory in Proposition \ref{prop:backward-sharp-value-function} is the rescaled form of \eqref{eq:freidlin-front-optimal-path}.
Formula \eqref{eq:backward-front-value-profile}, with $\nu=1-\frac{\sigma_*}{t}$, and the reduction $W=\max\{0,J\}$ give \eqref{eq:freidlin-front-value-profile}.
This proves strict positivity at every intermediate time when $\lambda>0$.
If $\lambda=0$, the cylindrical boundary in Theorem~\ref{thm:main-wulff-freidlin} gives $\rho=2\sqrt{r_{1}}t$.
The rescaled formula \eqref{eq:freidlin-front-optimal-path} then yields
\begin{align*}
  |\gamma_{*}(s)-(\gamma_{*}(s)\cdot\mathbf{e})\mathbf{e}|=2\sqrt{r_{1}}(t-s) \quad \text{and} \quad
  0\leqslant\gamma_{*}(s)\cdot\mathbf{e}\leqslant c(t-s) \quad \text{for }s\in[0,t).
\end{align*}
Thus the trajectory remains on the cylindrical front.
Substitution of $\sigma_*=t$ and $\sigma_*=0$ proves agreement with the straight trajectories at the junctions.

It remains to identify the minimizers of the min-max formula.
For the unique unconstrained minimizer $\gamma_{*}$, its time-shifted tail also minimizes the unconstrained optimal control problem, as in the proof of Theorem~\ref{thm:freidlin-sharp}.
Consequently,
\begin{align*}
  A_{\gamma_{*}}(s)=J(\mathbf{x},t)-J(\gamma_{*}(s),t-s)=-J(\gamma_{*}(s),t-s)\leqslant0, \quad s\in[0,t),
\end{align*}
and $A_{\gamma_{*}}(t)=J(\mathbf{x},t)=0$.
Hence $\gamma_{*}$ attains the min-max value $W(\mathbf{x},t)=0$.
Conversely, every minimizing trajectory $\gamma$ for the min-max formula satisfies
\begin{align*}
  0=J(\mathbf{x},t)\leqslant A_\gamma(t)\leqslant\max_{s\in[0,t]}A_\gamma(s)=W(\mathbf{x},t)=0.
\end{align*}
Thus $\gamma$ also minimizes the unconstrained optimal control problem and agrees with $\gamma_{*}$ on $[0,t]$ by Proposition~\ref{prop:backward-sharp-value-function}. This completes the proof.
\end{proof}

\subsection{Wulff shapes and optimal trajectories in the non-Freidlin regime}

We next evaluate the convex hull representation obtained in Theorem \ref{thm:zero-set-slow-interface}.

\begin{proposition}
\label{prop:explicit-wulff-slow-interface}
For $\mathbf{x}\in\mathbb{R}^{N}$, as before, we denote $\xi=\mathbf{x}\cdot\mathbf{e}$ and $\rho=\left|\mathbf{x}-\xi\mathbf{e}\right|$.
Then, for every $t\in\mathbb{R}^{+}$, the propagation region $\mathcal{Z}(t)$ and its boundary are given as follows.
\begin{enumerate}[label=(\roman*),font=\normalfont, leftmargin=*]
  \item If $-2\sqrt{r_{2}}<c\leqslant-2\sqrt{r_{1}}$, then
  \begin{align*}
  \mathcal Z(t)=\Sigma_{2}^{c}(t).
  \end{align*}
Moreover,
\begin{align*}
  \partial\mathcal Z(t)=\left\{\mathbf{x}\in\mathbb{R}^{N} \mid \xi^{2}+\rho^{2}=4r_{2}t^{2}, \; \xi\geqslant ct \right\} \cup \left\{\mathbf{x}\in\mathbb{R}^{N} \,\middle|\, \xi=ct, \; 0\leqslant\rho\leqslant t\sqrt{4r_{2}-c^{2}} \right\}.
\end{align*}

  \item If $-2\sqrt{r_{1}}<c<2\alpha$, define $\xi_{*}$ and $\rho_{*}$ as in \eqref{eq:contact-point}. Then
  \begin{align}
\label{eq:zt-non-freidlin-2}
  \mathcal Z(t)=\Sigma_{1}^{\xi_{*}}(t) \cup \left\{\mathbf{x}\in\mathbb{R}^{N} \,\middle|\, \xi_{*}t\leqslant\xi\leqslant ct, \; \xi_{*}\xi+\rho_{*}\rho\leqslant4r_{1}t \right\} \cup \Sigma_{2}^{c}(t).
  \end{align}
  Moreover, $\partial\mathcal{Z}(t)=\Gamma_{1}(t)\cup\Gamma_{\mathrm{m}}(t)\cup\Gamma_{2}(t)$, where
  \begin{align*}
  &\Gamma_{1}(t)= \left\{\mathbf{x}\in\mathbb{R}^{N} \mid \xi^{2}+\rho^{2}=4r_{1}t^{2}, \; -2\sqrt{r_{1}}t\leqslant\xi\leqslant\xi_{*}t \right\},\\
  &\Gamma_{\mathrm{m}}(t)= \left\{\mathbf{x}\in\mathbb{R}^{N} \mid \xi_{*}\xi+\rho_{*}\rho=4r_{1}t, \; \xi_{*}t\leqslant\xi\leqslant ct \right\},\\
  &\Gamma_{2}(t)= \left\{\mathbf{x}\in\mathbb{R}^{N} \mid \xi^{2}+\rho^{2}=4r_{2}t^{2}, \; ct\leqslant\xi\leqslant2\sqrt{r_{2}}t \right\}.
  \end{align*}
\end{enumerate}
\end{proposition}

\begin{proof}
By \eqref{eq:propagation-region-scaling} and the definitions of $\Sigma_{j}^{q}(t)$ and $\Gamma_{j}(t)$, it suffices to prove the assertions for $t=1$.
Theorem \ref{thm:zero-set-slow-interface} gives
\begin{align}
\label{eq:wulff-slow-convex-hull}
  \mathcal{Z}(1)=\mathrm{co} \left(\overline {B_{2\sqrt{r_{1}}}(\mathbf{0})} \cup \Sigma_{2}^{c}(1) \right).
\end{align}

Suppose first that $-2\sqrt{r_{2}}<c\leqslant-2\sqrt{r_{1}}$.
Every $\mathbf{x}\in\overline {B_{2\sqrt{r_{1}}}(\mathbf{0})}$ satisfies $\xi\geqslant-|\mathbf{x}|\geqslant-2\sqrt{r_{1}}\geqslant c$.
Hence $\overline {B_{2\sqrt{r_{1}}}(\mathbf{0})} \subset\Sigma_{2}^{c}(1)$.
Since $\Sigma_{2}^{c}(1)$ is convex, \eqref{eq:wulff-slow-convex-hull} yields $\mathcal{W}=\Sigma_{2}^{c}(1)$.
Its boundary consists of the spherical part $\xi^{2}+\rho^{2}=4r_{2}$ with $\xi\geqslant c$ and the planar part $\xi=c$ with $0\leqslant\rho\leqslant\sqrt{4r_{2}-c^{2}}$.
This proves assertion (i) for $t=1$.

Assume now that $-2\sqrt{r_{1}}<c<2\alpha$.
Since the sets in \eqref{eq:wulff-slow-convex-hull} depend on $\mathbf{x}$ only through $\xi$ and $\rho$, we work in the variables $(\xi,\rho)$. Set
\begin{align*}
  P=\left(c,\sqrt{4r_{2}-c^{2}} \right), \qquad T=(\xi_{*},\rho_{*}).
\end{align*}
A direct calculation implies
\begin{align}
\label{eq:tangent-condition}
  \xi_{*}^{2}+\rho_{*}^{2}=4r_{1}, \qquad c\xi_{*}+\sqrt{4r_{2}-c^{2}}\,\rho_{*}=4r_{1},
\end{align}
together with $-2\sqrt{r_{1}}<\xi_{*}<c$ and $\rho_{*}>0$. Thus the line
\begin{align*}
  \ell: \quad \xi_{*}\xi+\rho_{*}\rho=4r_{1}
\end{align*}
is tangent to $\xi^{2}+\rho^{2}=4r_{1}$ at $T$ and passes through $P$.

We now claim that
\begin{align}
\label{eq:large-region-below-tangent}
  \left\{(\xi,\rho)\in\mathbb{R}\times[0,\infty) \mid \xi^{2}+\rho^{2}\leqslant4r_{2}, \; \xi\geqslant c \right\} \subset \left\{(\xi,\rho)\in\mathbb{R}\times[0,\infty) \mid \xi_{*}\xi+\rho_{*}\rho\leqslant4r_{1} \right\}.
\end{align}
For $\xi\in[c,2\sqrt{r_{2}}]$, let
\begin{align*}
  h(\xi)=\frac{4r_{1}-\xi_{*}\xi}{\rho_{*}}, \qquad f(\xi)=\sqrt{4r_{2}-\xi^{2}}.
\end{align*}
The second identity in \eqref{eq:tangent-condition} gives $h(c)=f(c)$, while
\begin{align*}
  h'(\xi)-f'(c)=-\frac{\xi_{*}}{\rho_{*}}+\frac{c}{\sqrt{4r_{2}-c^{2}}}=\frac{4\alpha\sqrt{r_{1}}} {\rho_{*}\sqrt{4r_{2}-c^{2}}} >0.
\end{align*}
Since $f$ is strictly concave, $f'(\xi)\leqslant f'(c)<h'(\xi)$ for $\xi\in[c,2\sqrt{r_{2}})$.
Thus $h-f$ is strictly increasing on this interval.
Using $h(c)=f(c)$ and continuity at $\xi=2\sqrt{r_{2}}$, we obtain
\begin{align*}
  f(\xi)\leqslant h(\xi) \quad\text{for } \xi\in[c,2\sqrt{r_{2}}],
\end{align*}
which proves \eqref{eq:large-region-below-tangent}.

Let $\mathcal C$ denote the right-hand side of \eqref{eq:zt-non-freidlin-2} at $t=1$.
We first show that $\mathcal Z(1)\subset\mathcal C$. Extend $h$ to $[\xi_{*},c]$ by the same formula and define
\begin{align*}
  F(\xi)=
  \begin{cases}
    \sqrt{4r_{1}-\xi^{2}}, &\text{if } -2\sqrt{r_{1}}\leqslant\xi\leqslant\xi_{*},\\
    h(\xi), &\text{if }  \xi_{*}\leqslant\xi\leqslant c,\\
    f(\xi), &\text{if }  c\leqslant\xi\leqslant2\sqrt{r_{2}}.
  \end{cases}
\end{align*}
Then $F$ is continuous and piecewise concave with $C^{1}$ matching at $\xi_{*}$ and a concave kink at $c$.
Hence $F$ is concave and $\mathcal C$ is convex.
The tangent line lies above the small upper semicircle in the $(\xi,\rho)$-plane, and the large ball contains the small one.
Thus $\mathcal C$ contains both $\overline{B_{2\sqrt{r_{1}}}(\mathbf{0})}$ and $\Sigma_{2}^{c}(1)$, and hence contains their convex hull. Therefore, \eqref{eq:wulff-slow-convex-hull} implies $\mathcal Z(1)\subset\mathcal C$.

We now show $\mathcal C\subset\mathcal Z(1)$.
\eqref{eq:wulff-slow-convex-hull} shows that $\mathcal Z(1)$ is convex and contains both $\overline{B_{2\sqrt{r_{1}}}(\mathbf{0})}$ and $\Sigma_{2}^{c}(1)$.
Thus the left and right portions of $\mathcal C$ already belong to $\mathcal Z(1)$.
For the intermediate portion, consider any $\mathbf{x}=\xi\mathbf{e}+{\boldsymbol\zeta}\in\mathcal C$ with $\xi_{*}<\xi<c$. We have
${\boldsymbol\zeta}\perp\mathbf{e}$ and $|{\boldsymbol\zeta}|\leqslant h(\xi)$.
Set
\begin{align*}
  \theta=\frac{\xi-\xi_{*}}{c-\xi_{*}},\qquad  d_{c}=\sqrt{4r_{2}-c^{2}}.
\end{align*}
By \eqref{eq:tangent-condition}, we have $h(\xi)=(1-\theta)\rho_{*}+\theta d_{c}>0$. Moreover, we can verify that
\begin{align*}
  \mathbf{x}
  =(1-\theta)\left(\xi_{*}\mathbf{e}+\frac{\rho_{*}}{h(\xi)}{\boldsymbol\zeta}\right)
  +\theta\left(c\mathbf{e}+\frac{d_{c}}{h(\xi)}{\boldsymbol\zeta}\right).
\end{align*}
Since $\frac{|{\boldsymbol\zeta}|}{h(\xi)}\leqslant1$, the two points on the right belong to $\overline{B_{2\sqrt{r_{1}}}(\mathbf{0})}$ and $\Sigma_{2}^{c}(1)$, respectively.
Therefore, \eqref{eq:wulff-slow-convex-hull} yields $\mathbf{x}\in\mathcal Z(1)$.
Hence \eqref{eq:zt-non-freidlin-2} holds for $t=1$, and $\partial{\mathcal{Z}}(1)=\Gamma_{1}(1)\cup\Gamma_{\mathrm{m}}(1)\cup\Gamma_{2}(1)$. This proves assertion (ii), and the proof is complete.
\end{proof}

\begin{proof}[Proof of Theorem \ref{thm:main-wulff-nonfreidlin}]
It follows from Proposition \ref{prop:explicit-wulff-slow-interface} by setting $t=1$ and using \eqref{eq:wulff-shape-definition}.
\end{proof}


In the non-Freidlin regime, we characterize the optimal trajectories for the min-max formula \eqref{eq:stopping-value-inf-action} at the propagation front and determine where uniqueness holds.


\begin{proposition}
\label{prop:nonfreidlin-front-trajectories}
Assume that $-2\sqrt{r_{2}}<c<2\alpha$, fix $t\in\mathbb{R}^{+}$, and let $\mathbf{x}\in t\partial\mathcal{W}$.
On the spherical and conical parts of the front, including their junctions with the other front pieces, the min-max formula \eqref{eq:stopping-value-inf-action} has a unique minimizing trajectory $\gamma_{*}$ on $[0,t]$.
In both cases, $\gamma_{*}(s)=t\gamma_1(\frac{s}{t})$ for $s\in[0,t]$, where $\gamma_1$ is given by \eqref{eq:front-straight-optimal-path} on the spherical parts and by \eqref{eq:cone-front-optimal-path} on the conical part, with $\frac{\mathbf{x}}{t}$ in place of $\mathbf{x}$.
These trajectories have zero running cost almost everywhere and satisfy $\gamma_{*}(s)\in(t-s)\partial\mathcal{W}$ for $0\leqslant s<t$.
In the relative interior of the planar part, namely when
\begin{align*}
  -2\sqrt{r_{2}}<c\leqslant-2\sqrt{r_{1}}, \qquad \xi=ct, \qquad \rho<t\sqrt{4r_{2}-c^{2}},
\end{align*}
there are infinitely many minimizing trajectories.
\end{proposition}

\begin{proof}
Fix $t\in \mathbb{R}^{+}$ and $\mathbf{x}\in\partial\mathcal Z(t)=t\partial\mathcal{W}$.
By Proposition~\ref{prop:inf-action-characterization}, a trajectory $\gamma\in\mathbb X^{\mathbf{x},t}$ minimizes the min-max formula if and only if
\begin{align}
\label{eq:front-path-nonpositive-action}
  A_\gamma(\sigma)\leqslant0 \quad\text{for every }\sigma\in[0,t].
\end{align}
All paths below are considered on $[0,t]$ and extended by zero after time $t$.
We prove the assertions in four steps.

\textbf{Step 1.} We consider the spherical part in the higher-growth region.
If $|\mathbf{x}|=2\sqrt{r_{2}}t$ and $\xi\geqslant ct$, one can verify that the straight trajectory $\gamma_{*}(s)=(1-\frac{s}{t})\mathbf{x}$ stays in the higher-growth region or on the interface and has zero running cost.
It therefore satisfies \eqref{eq:front-path-nonpositive-action}.
For any minimizing trajectory, $R\leqslant r_{2}$ and the Cauchy-Schwarz inequality give
\begin{align*}
  0\geqslant A_\gamma(t)\geqslant\frac{1}{4}\int_0^t|\dot\gamma(s)|^{2}\,ds-r_{2}t\geqslant\frac{|\mathbf{x}|^{2}}{4t}-r_{2}t=0.
\end{align*}
Equality forces $\dot\gamma=-\frac{\mathbf{x}}{t}$ almost everywhere, so the straight trajectory is the unique minimizer on $[0,t]$.
This includes the junctions of the spherical and planar parts. 

\textbf{Step 2.} We consider the lower spherical part when $-2\sqrt{r_{1}}<c<2\alpha$.
Let $|\mathbf{x}|=2\sqrt{r_{1}}t$ and $\xi\leqslant\xi_*t$.
The straight trajectory lies in the lower-growth region before reaching the origin and has zero running cost.
To prove uniqueness, let $\gamma$ be any minimizing trajectory and let $\tau$ be its first contact time with the shifting interface, as defined in \eqref{eq:contact-time}.
Since $\xi\leqslant\xi_*t<ct$, we have $0<\tau\leqslant t$.
We claim that $\tau=t$.

Suppose, to the contrary, that $\tau<t$.
Use the notation $B(\tau)$ and $\rho_{\tau}(\xi,t)$ of Lemma \ref{lem:fixed-contact-optimization}, and choose $\mathbf{n}=\frac{{\boldsymbol\zeta}}{\rho}$ if $\rho>0$, or any unit vector orthogonal to $\mathbf{e}$ if $\rho=0$.
We compare $\mathbf{x}$ with the endpoint attaining the bound $\rho_\tau(\xi,t)$ in this lemma.
To describe this endpoint, set $d_c=\sqrt{4r_{2}-c^{2}}$ and
\begin{align*}
  P=c\mathbf{e}+d_c\mathbf{n}, \quad
  q_{\tau}=\frac{\xi-c(t-\tau)}{\tau}\mathbf{e}+\frac{\sqrt{B(\tau)}}{\tau}\mathbf{n}.
\end{align*}
Then $|q_{\tau}|=2\sqrt{r_{1}}$, and the sharp endpoint in that lemma satisfies
\begin{align}
  \label{eq:front-sharp-endpoint-decomposition}
  \mathbf{x}_{\tau}=\xi\mathbf{e}+\rho_{\tau}(\xi,t)\mathbf{n}
  =(t-\tau)P+\tau q_{\tau}, \quad \text{and} \quad \rho\leqslant\rho_{\tau}(\xi,t).
\end{align}
Denote $v=\frac{\mathbf{x}}{t}$.
The decomposition in \eqref{eq:front-sharp-endpoint-decomposition} expresses $\frac{\mathbf{x}_{\tau}}{t}$ as a convex combination of $P$ and $q_{\tau}$.
We compare $v\cdot\mathbf{x}_\tau$ with $v\cdot\mathbf{x}=4r_{1}t$.
We now show that
\begin{align}
  \label{eq:lower-sphere-interface-support}
  v\cdot P\leqslant4r_{1},
  \qquad\text{with equality if and only if }\xi=\xi_*t.
\end{align}
Indeed, define
\begin{align*}
  f(q)=cq+d_c\sqrt{4r_{1}-q^2},
  \qquad q\in[-2\sqrt{r_{1}},\xi_*].
\end{align*}
Since $v=\frac{\xi}{t}\mathbf{e}+\frac{\rho}{t}\mathbf{n}$ and $\xi^2+\rho^2=4r_{1}t^2$, we have
\begin{align*}
  v\cdot P
  =c\frac{\xi}{t}+d_c\frac{\rho}{t}
  =c\frac{\xi}{t}+d_c\sqrt{4r_{1}-\left(\frac{\xi}{t}\right)^2}
  =f\left(\frac{\xi}{t}\right).
\end{align*}
Direct differentiation gives $f''<0$ on $(-2\sqrt{r_{1}},\xi_*)$.
Using \eqref{eq:contact-point}, we also obtain
\begin{align*}
  f'(\xi_*)=\frac{c\rho_{*}-d_c\xi_*}{\rho_{*}}
  =\frac{4\alpha\sqrt{r_{1}}}{\rho_{*}}>0.
\end{align*}
Therefore, $f'(q)\geqslant f'(\xi_*)>0$ for $-2\sqrt{r_{1}}<q\leqslant\xi_*$.
Furthermore, $f$ is strictly increasing on $[-2\sqrt{r_{1}},\xi_*]$.
Moreover, \eqref{eq:tangent-condition} yields $f(\xi_*)=c\xi_*+d_c\rho_*=4r_{1}$.
Since $\frac{\xi}{t}\leqslant\xi_*$, we obtain \eqref{eq:lower-sphere-interface-support}.

Recall that $|q_{\tau}|=2\sqrt{r_{1}}$. Hence, $v\cdot q_{\tau}\leqslant4r_{1}$. Moreover, \eqref{eq:front-sharp-endpoint-decomposition} and \eqref{eq:lower-sphere-interface-support} imply
\begin{align*}
  4r_{1}t=v\cdot\mathbf{x}
  \leqslant v\cdot\mathbf{x}_{\tau}
  =(t-\tau)v\cdot P+\tau v\cdot q_{\tau}
  \leqslant4r_{1}t.
\end{align*}
If $\xi<\xi_*t$, then $v\cdot P<4r_{1}$, and the assumption $\tau<t$ makes the last inequality strict, a contradiction.
If $\xi=\xi_*t$, equality in the Cauchy-Schwarz inequality gives $q_\tau=v=\xi_*\mathbf{e}+\rho_*\mathbf{n}$.
Comparing the $\mathbf{e}$-components in \eqref{eq:front-sharp-endpoint-decomposition} yields $t\xi_*=(t-\tau)c+\tau\xi_*$.
Since $\tau<t$, this implies $\xi_*=c$, contradicting $\xi_*<c$.
Hence $\tau=t$, as claimed.
Thus every minimizing trajectory has $\tau=t$ and lies in the lower-growth region for $0\leqslant s<t$.
The Cauchy-Schwarz inequality gives
\begin{align*}
  0\geqslant A_\gamma(t)=\frac{1}{4}\int_0^t|\dot\gamma(s)|^2\,ds-r_{1}t\geqslant\frac{|\mathbf{x}|^2}{4t}-r_{1}t=0.
\end{align*}
Equality in the Cauchy-Schwarz inequality, together with the endpoint conditions, yields $\gamma(s)=(1-\frac{s}{t})\mathbf{x}$ for $s\in[0,t]$.
This proves uniqueness on the lower spherical part, including its junction with the conical part.

\textbf{Step 3.} We consider the strict conical part, where $\xi_*t<\xi<ct$ and $\xi_*\xi+\rho_{*}\rho=4r_{1}t$.
Use $P$ and $Q$ from \eqref{eq:cone-front-path-parameters}, with $\mathbf{n}=\frac{{\boldsymbol\zeta}}{\rho}$, and set $\tau_*=\frac{ct-\xi}{c-\xi_*}$.
The tangency relations \eqref{eq:tangent-condition} give
\begin{align}
  \label{eq:cone-front-tangency-identities}
  |Q|^{2}=Q\cdot P=4r_{1}, \qquad |P|^{2}=4r_{2}, \qquad |P-Q|^{2}=4\alpha^{2}.
\end{align}
By the definition of $\tau_*$, we have
\begin{align*}
  \xi=ct-\tau_*(c-\xi_*)=\tau_*\xi_*+(t-\tau_*)c.
\end{align*}
Using the equation $\xi_*\xi+\rho_{*}\rho=4r_{1}t$ and \eqref{eq:tangent-condition}, we obtain
\begin{align*}
  \rho
  =\frac{4r_{1}t-\xi_*\xi}{\rho_{*}}
  =\frac{\tau_*(4r_{1}-\xi_*^{2})+(t-\tau_*)(4r_{1}-c\xi_*)}{\rho_{*}}
  =\tau_*\rho_{*}+(t-\tau_*)d_c.
\end{align*}
Thus $\mathbf{x}=\tau_*Q+(t-\tau_*)P$, and hence $\mathbf{x}-\tau_*Q=(t-\tau_*)P$.
The rescaled trajectory $\gamma_{*}(s)=\mathbf{x}-sQ$ for $s\leqslant\tau_*$ and $\gamma_{*}(s)=(t-s)P$ for $s\geqslant\tau_*$ is therefore continuous at $s=\tau_*$.
For $s<\tau_*$,
\begin{align*}
  \gamma_{*}(s)\cdot\mathbf{e}-c(t-s)
  =-(c-\xi_*)(\tau_*-s)<0,
\end{align*}
while the second segment lies on the interface.
In particular, $|Q|=2\sqrt{r_{1}}$ and $|P|=2\sqrt{r_{2}}$.
Thus $A_{\gamma_{*}}(\sigma)=0$ for every $\sigma\in[0,t]$.
Since every admissible trajectory $\gamma$ satisfies $\max_{\sigma\in[0,t]}A_\gamma(\sigma)\geqslant A_\gamma(0)=0$, the trajectory $\gamma_{*}$ attains the minimum $0$ in \eqref{eq:stopping-value-inf-action}.


Let $\gamma$ be any minimizing trajectory, and let $\tau$ be its first contact time as in \eqref{eq:contact-time}.
If $\tau=t$, \eqref{eq:front-path-nonpositive-action} would imply $|\mathbf{x}|\leqslant2\sqrt{r_{1}}t$.
However, \eqref{eq:cone-front-tangency-identities} together with $\mathbf{x}=\tau_*Q+(t-\tau_*)P$ gives
\begin{align*}
  |\mathbf{x}|^{2}=4r_{1}t^{2}+4\alpha^{2}(t-\tau_*)^{2}>4r_{1}t^{2}.
\end{align*}
Hence $0<\tau<t$, and Lemma \ref{lem:fixed-contact-optimization} applies.
Define $q_{\tau}$ and $\mathbf{x}_{\tau}$ as in Step 2, so that \eqref{eq:front-sharp-endpoint-decomposition} holds.
Using $Q\cdot q_{\tau}\leqslant4r_{1}$, \eqref{eq:cone-front-tangency-identities} and $\rho_{*}>0$, we obtain
\begin{align*}
  4r_{1}t=Q\cdot\mathbf{x}
  \leqslant Q\cdot\mathbf{x}_{\tau}
  =(t-\tau)Q\cdot P+\tau Q\cdot q_{\tau}
  \leqslant4r_{1}t.
\end{align*}
Consequently, $\rho_{\tau}(\xi,t)=\rho$ and $q_{\tau}=Q$.
Then $\mathbf{x}_\tau=\mathbf{x}$.
Together with $q_\tau=Q$, \eqref{eq:front-sharp-endpoint-decomposition} gives
\begin{align}\label{eq:x-decompose}
  \mathbf{x}=(t-\tau)P+\tau Q.
\end{align}
Taking the inner product with $\mathbf{e}$, we obtain
\begin{align*}
  \xi=(t-\tau)c+\tau\xi_*=ct-\tau(c-\xi_*).
\end{align*}
Since $c>\xi_*$, it follows that $\tau=\frac{ct-\xi}{c-\xi_*}=\tau_*$.

It remains to prove that $\gamma(s)=\gamma_{*}(s)$ for all $s\in[0,t]$.
By the definition of $\tau$, we may write $\gamma(\tau)=\mathbf{p}+c(t-\tau)\mathbf{e}$ with $\mathbf{p}\cdot\mathbf{e}=0$, as in \eqref{eq:gamm-deco}.
Let $\tilde{\mathbf{p}}$ be the closest point to $\mathbf{p}$ on $[\mathbf{0},{\boldsymbol\zeta}]$.
We claim that $\tilde{\mathbf{p}}=\mathbf{p}$, or equivalently, $\mathbf{p}\in[\mathbf{0},{\boldsymbol\zeta}]$.
Indeed, as in the proof of Lemma \ref{lem:fixed-contact-optimization}, set $\rho_1=|{\boldsymbol\zeta}-\tilde{\mathbf{p}}|$ and $\rho_2=|\tilde{\mathbf{p}}|$.
The estimate in that proof implies
\begin{align*}
  \rho\leqslant G_{\tau}(\rho_1)
  \leqslant G_{\tau}(\sqrt{B(\tau)})
  =\rho_{\tau}(\xi,t)=\rho.
\end{align*}
Since $G_{\tau}$ is strictly increasing on $[0,\sqrt{B(\tau)}]$, we have $\rho_1=\sqrt{B(\tau)}$ and $\rho_2=(t-\tau)d_c$.
Set $\delta=ct-\xi$.
Using $\mathbf{x}=\xi\mathbf{e}+{\boldsymbol\zeta}$ and $\gamma(\tau)=\mathbf{p}+c(t-\tau)\mathbf{e}$, we obtain
\begin{align*}
  \mathbf{x}-\gamma(\tau)=({\boldsymbol\zeta}-\mathbf{p})+(\xi-ct+c\tau)\mathbf{e}=({\boldsymbol\zeta}-\mathbf{p})-(\delta-c\tau)\mathbf{e}.
\end{align*}
It follows from the fact ${\boldsymbol\zeta}\cdot\mathbf{e}=\mathbf{p}\cdot\mathbf{e}=0$ that
\begin{align}\label{eq:diff-x-to-gam}
  |\mathbf{x}-\gamma(\tau)|^{2}
  =|{\boldsymbol\zeta}-\mathbf{p}|^{2}-2(\delta-c\tau)({\boldsymbol\zeta}-\mathbf{p})\cdot\mathbf{e}+(\delta-c\tau)^{2}
  =|{\boldsymbol\zeta}-\mathbf{p}|^{2}+(\delta-c\tau)^{2}.
\end{align}
Since $\rho_1=\tau\rho_*>0$ and $\rho_2=(t-\tau)d_c>0$, $\tilde{\mathbf{p}}$ must lie strictly between $\mathbf{0}$ and ${\boldsymbol\zeta}$.
If $\mathbf{p}$ lies on the line through $\mathbf{0}$ and ${\boldsymbol\zeta}$, then $\mathbf{p}=\tilde{\mathbf{p}}$; otherwise, $(\mathbf{p}-\tilde{\mathbf{p}})\cdot({\boldsymbol\zeta}-\tilde{\mathbf{p}})=0$.
In either case, we have
\begin{align*}
  |{\boldsymbol\zeta}-\mathbf{p}|^{2}
  =|{\boldsymbol\zeta}-\tilde{\mathbf{p}}|^{2}+|\mathbf{p}-\tilde{\mathbf{p}}|^{2}.
\end{align*}
Using $|{\boldsymbol\zeta}-\mathbf{p}|\geqslant\rho_1$ and the Cauchy-Schwarz inequality, we deduce
\begin{align}\label{eq:action-chain}
  0\geqslant A_\gamma(\tau)
  \geqslant\frac{|\mathbf{x}-\gamma(\tau)|^{2}}{4\tau}-r_{1}\tau
  \geqslant\frac{\rho_1^{2}+(\delta-c\tau)^{2}}{4\tau}-r_{1}\tau=0.
\end{align}
Therefore, it holds that $|\mathbf{x}-\gamma(\tau)|^{2}=\rho_1^{2}+(\delta-c\tau)^{2}$.
Substituting it into \eqref{eq:diff-x-to-gam}, we obtain $|{\boldsymbol\zeta}-\mathbf{p}|^{2}=\rho_1^{2}=|{\boldsymbol\zeta}-\tilde{\mathbf{p}}|^{2}$.
Hence
\begin{align*}
  |\mathbf{p}-\tilde{\mathbf{p}}|^{2}=|{\boldsymbol\zeta}-\mathbf{p}|^{2}-|{\boldsymbol\zeta}-\tilde{\mathbf{p}}|^{2}=\rho_1^{2}-\rho_1^{2}=0.
\end{align*}
Thus $\mathbf{p}=\tilde{\mathbf{p}}$, and the claim is proved.

Since $\mathbf{p}=\tilde{\mathbf{p}}\in[\mathbf{0},{\boldsymbol\zeta}]$ and ${\boldsymbol\zeta}=\rho\mathbf{n}$, we have $\mathbf{p}=\rho_2\mathbf{n}=(t-\tau)d_c\mathbf{n}$.
Thus $\gamma(\tau)=\mathbf{p}+c(t-\tau)\mathbf{e}=(t-\tau)P$.
Using \eqref{eq:x-decompose}, we obtain $\gamma(\tau)-\mathbf{x}=-\tau Q$.

The bounds \eqref{eq:action-chain} imply that $A_\gamma(\tau)=0$.
Since $\tau$ is the first contact time, $R(\gamma(s),t-s)=r_{1}$ for $0\leqslant s<\tau$.
The definition of $A_\gamma$ therefore gives
\begin{align*}
  \int_0^\tau|\dot\gamma(s)|^2\,ds=4\left(A_\gamma(\tau)+r_{1}\tau\right)&=4r_{1}\tau=\tau|Q|^2
  =\frac{|\gamma(\tau)-\mathbf{x}|^2}{\tau}.
\end{align*}
Equality in the Cauchy-Schwarz inequality implies $\dot\gamma(s)=\frac{\gamma(\tau)-\mathbf{x}}{\tau}=-Q$ for almost every $s\in(0,\tau)$.
Together with $\gamma(0)=\mathbf{x}$, this yields $\gamma(s)=\mathbf{x}-sQ$ for $s\in[0,\tau]$.

On $[\tau,t]$, using $A_\gamma(t)\leqslant0$, $A_\gamma(\tau)=0$, and $R\leqslant r_{2}$, we obtain
\begin{align*}
  0\geqslant A_\gamma(t)-A_\gamma(\tau)\geqslant\frac{1}{4}\int_\tau^t|\dot\gamma(s)|^2\,ds-r_{2}(t-\tau)\geqslant\frac{|\gamma(\tau)|^2}{4(t-\tau)}-r_{2}(t-\tau)=0.
\end{align*}
The last equality follows from $\gamma(\tau)=(t-\tau)P$ and $|P|^2=4r_{2}$.
Equality in the Cauchy-Schwarz inequality gives $\dot\gamma(s)=-\frac{\gamma(\tau)}{t-\tau}=-P$ for almost all $s\in(\tau,t)$.
Since $\gamma(t)=\mathbf{0}$, we obtain $\gamma(s)=(t-s)P$ for $s\in[\tau,t]$.
Together with $\tau=\tau_*$, this proves that $\gamma=\gamma_{*}$ on $[0,t]$, and hence establishes uniqueness on the strict conical part.
The formulas at its two junctions follow from Steps 1 and 2.

We next verify that these optimal trajectories satisfy $\gamma_{*}(s)\in(t-s)\partial\mathcal{W}$ for every $0\leqslant s<t$.
For $\tau_*\leqslant s<t$, we have $\frac{\gamma_{*}(s)}{t-s}=P$. Since $P\in\partial\mathcal{W}$, the assertion follows immediately.
For $0\leqslant s\leqslant\tau_*$, using $\gamma_{*}(s)=\mathbf{x}-sQ$ and $\mathbf{x}=\tau_*Q+(t-\tau_*)P$, we obtain
\begin{align*}
  \frac{\gamma_{*}(s)}{t-s}=\frac{\mathbf{x}-sQ}{t-s}=\frac{(\tau_*-s)Q+(t-\tau_*)P}{t-s}=\frac{\tau_*-s}{t-s}Q+\frac{t-\tau_*}{t-s}P.
\end{align*}
Then $\frac{\gamma_{*}(s)}{t-s}=\frac{\mathbf{x}-sQ}{t-s}$ is a convex combination of $P$ and $Q$, so it belongs to $[Q,P]$, which lies on the conical part of $\partial\mathcal{W}$.
Hence $\gamma_{*}(s)\in(t-s)\partial\mathcal{W}$ also holds on the first segment.
For the spherical trajectories, $\gamma_{*}(s)=(1-\frac{s}{t})\mathbf{x}$ gives
\begin{align*}
  \frac{\gamma_{*}(s)}{t-s}=\frac{\mathbf{x}}{t}\in\partial\mathcal{W}\quad \text{for } s\in[0,t),
\end{align*}
which proves the same assertion.

\textbf{Step 4.} We consider the relative interior of the planar part when $-2\sqrt{r_{2}}<c\leqslant-2\sqrt{r_{1}}$.
Here $\xi=ct$ and $|\mathbf{x}|<2\sqrt{r_{2}}t$.
For each $\epsilon$ satisfying $0<\epsilon<2\sqrt{r_{2}}-\frac{|\mathbf{x}|}{t}$, define
\begin{align*}
  \gamma_{\epsilon}(s)=\left(1-\frac{s}{t}\right)\mathbf{x}+\epsilon\frac{s(t-s)}{t}\mathbf{e}, \quad s\in[0,t].
\end{align*}
These trajectories have the required endpoints and satisfy
\begin{align*}
|\dot\gamma_{\epsilon}(s)|\leqslant\frac{|\mathbf{x}|}{t}+\epsilon<2\sqrt{r_{2}}, \quad\text{and} \quad
  \gamma_{\epsilon}(s)\cdot\mathbf{e}-c(t-s)
  =\epsilon\frac{s(t-s)}{t}>0 \quad\text{for }s\in(0,t).
\end{align*}
Thus their running costs are negative almost everywhere, so $A_{\gamma_{\epsilon}}(\sigma)\leqslant0$ for every $\sigma\in[0,t]$.
Since $A_{\gamma_{\epsilon}}(0)=0$, we have
\begin{align*}
  \max_{\sigma\in[0,t]}A_{\gamma_{\epsilon}}(\sigma)=0=W(\mathbf{x},t).
\end{align*}
Thus each trajectory constructed above attains the minimum in \eqref{eq:stopping-value-inf-action}.
They are distinct minimizing trajectories, which proves nonuniqueness at every relative interior point of the planar part.
Together with the front decomposition in Proposition \ref{prop:explicit-wulff-slow-interface}, the four steps cover all non-Freidlin regime.
\end{proof}

\begin{proof}[Proof of Theorem \ref{thm:main-front-trajectories}]
The assertions in the Freidlin regimes follow from Proposition \ref{prop:freidlin-front-trajectories} with $t=1$.
The assertions in the non-Freidlin regime follow from Proposition \ref{prop:nonfreidlin-front-trajectories} with $t=1$.
\end{proof}

\subsection{Directional spreading speeds}
We now derive the directional spreading speeds from the Wulff boundaries obtained above.

\begin{proof}[Proof of Corollary \ref{cor:wulff-directional-spreading}]
The explicit Wulff shape descriptions established above show that $\mathcal{W}$ is convex.
Moreover, \eqref{eq:wulff-contained-ball} and \eqref{eq:wulff-contains-ball} imply that
$\overline{B_{2\sqrt{r_{1}}}(\mathbf{0})} \subset \mathcal{W}\subset \overline{B_{2\sqrt{r_{2}}}(\mathbf{0})}$.
Since $\mathcal{W}$ is closed, it is compact and $\mathbf{0}\in\textrm{Int}\mathcal{W}$.

For $\mathbf{n}\in\mathbb S^{N-1}$, set $I_\mathbf{n}=\{s\geqslant0\mid s\mathbf{n}\in\mathcal{W}\}$.
Closedness and convexity imply that $I_\mathbf{n}$ is a closed interval. Together with the ball inclusions $\overline{B_{2\sqrt{r_{1}}}(\mathbf{0})} \subset \mathcal{W}\subset \overline{B_{2\sqrt{r_{2}}}(\mathbf{0})}$, we have
\begin{align*}
  I_\mathbf{n}=[0,c_\mathrm{dir}(\mathbf{n})], \qquad 2\sqrt{r_{1}} \leqslant c_\mathrm{dir}(\mathbf{n}) \leqslant 2\sqrt{r_{2}}.
\end{align*}
The endpoint $c_\mathrm{dir}(\mathbf{n})\mathbf{n}$ cannot belong to $\textrm{Int}\mathcal{W}$, since otherwise $I_\mathbf{n}$ could be extended.
Hence it belongs to $\partial\mathcal{W}$.

Moreover, $q\mathbf{n}\in\textrm{Int}\mathcal{W}$ whenever $0\leqslant q<c_\mathrm{dir}(\mathbf{n})$.
Indeed, if $B_{\epsilon}(\mathbf{0})\subset\mathcal{W}$ and $\theta=\frac{q}{c_\mathrm{dir}(\mathbf{n})}<1$, convexity gives $B_{(1-\theta)\epsilon}(q\mathbf{n})\subset\mathcal{W}$.

Fix any nonempty compact set $K\subset\mathbb{R}^{N}$ and denote $R_{K}=\max_{\mathbf{x}\in K}|\mathbf{x}|$.
If $0\leqslant q<c_\mathrm{dir}(\mathbf{n})$, let $\bar{d}=\textrm{dist}(q\mathbf{n},\overline{\mathcal{W}}^{\,c})>0$.
Choose $\eta>0$ so small that Theorem \ref{thm:wulff} applies and $\eta<\frac{\bar{d}}{2}$.
For large $t$,
\begin{align*}
  \textrm{dist} \left(q\mathbf{n}+\frac{\mathbf{x}}{t}, \overline{\mathcal{W}}^{\,c} \right) \geqslant \bar{d}-\frac{R_{K}}{t}>\eta \quad\text{for every }\mathbf{x}\in K.
\end{align*}
The persistence assertion of Theorem \ref{thm:wulff} gives the first conclusion.

If $q>c_\mathrm{dir}(\mathbf{n})$, then $d_{+}=\textrm{dist}(q\mathbf{n},\mathcal{W})>0$.
Choosing $\eta>0$ sufficiently small with $\eta<\frac{d_{+}}{2}$, the same argument implies, for all sufficiently large $t$,  that
\begin{align*}
  \textrm{dist} \left(q\mathbf{n}+\frac{\mathbf{x}}{t}, \mathcal{W}\right)>\eta \quad\text{for every }\mathbf{x}\in K.
\end{align*}
The extinction assertion of Theorem \ref{thm:wulff} gives the second conclusion.
\end{proof}

\begin{proof}[Proof of Corollary \ref{cor:main-directional-speed}]
By Corollary \ref{cor:wulff-directional-spreading}, $c_\mathrm{dir}(\mathbf{n})\mathbf{n}$ is the unique intersection of the ray $\{s\mathbf{n}\mid s\geqslant0\}$ with $\partial\mathcal{W}$.
It therefore remains to compute this intersection from the boundary descriptions obtained above.
Writing $\mathbf{x}=s\mathbf{n}$, we have $\xi=s\beta$ and $\rho=s\sqrt{1-\beta^{2}}$.

\textbf{(i)} Assume that $c\geqslant2(\alpha+\sqrt{r_{1}})$.
By \eqref{eq:wulff-shape-definition} and \eqref{eq:extremely-fast-front} with $t=1$,
\begin{align*}
  \partial\mathcal{W}= \left\{\mathbf{x}\in\mathbb{R}^{N} \,\middle|\, |\mathbf{x}|=2\sqrt{r_{1}} \right\}.
\end{align*}
Since $|\mathbf{n}|=1$, it follows that $c_\mathrm{dir}(\mathbf{n})=2\sqrt{r_{1}}$.

\textbf{(ii)} Assume that $2\sqrt{r_{2}}<c<2(\alpha+\sqrt{r_{1}})$.
By \eqref{eq:intermediate-fast-front}, $\partial\mathcal{W}=\Gamma_{1}(1)\cup\Gamma_{\mathrm{m}}(1)$.
On $\Gamma_{1}(1)$, we have $s=2\sqrt{r_{1}}$, and $\xi\leqslant2\lambda$ is equivalent to $\beta\leqslant \frac{\lambda}{\sqrt{r_{1}}}$.
On $\Gamma_{\mathrm{m}}(1)$, substitution into \eqref{eq:intermediate-fast-middle} implies
\begin{align}
\label{eq:directional-middle-quadratic}
  Q_{\beta}(s) &:= (1-\beta^{2})s^{2} +4\lambda\beta s -4(r_{1}+\lambda^{2})=0.
\end{align}
Since $\beta>0$, $Q_{\beta}$ is strictly increasing on $(0,\infty)$, and its unique positive root is $c_{\lambda}(\beta)$. Moreover,
\begin{align*}
  Q_{\beta}\left(\frac{2\lambda}{\beta}\right)<0, \qquad Q_{\beta}\left(\frac{c_\mathrm{nlp}}{\beta}\right)\geqslant0,
\end{align*}
where $\lambda c_\mathrm{nlp}=r_{1}+\lambda^{2}$, so this root satisfies the restrictions defining $\Gamma_{\mathrm{m}}(1)$.
At $\beta=\frac{\lambda}{\sqrt{r_{1}}}$, the two branches agree.
This proves assertion (ii).

\textbf{(iii)} Assume that $2\alpha\leqslant c\leqslant2\sqrt{r_{2}}$.
By \eqref{eq:front-decomposition}, $\partial\mathcal{W}=\Gamma_{1}(1)\cup\Gamma_{\mathrm{m}}(1)\cup\Gamma_{2}(1)$. In this case,
\begin{align*}
  0 \leqslant \frac{\lambda}{\sqrt{r_{1}}}<\frac{c}{2\sqrt{r_{2}}} \leqslant1.
\end{align*}
On $\Gamma_{1}(1)$, the condition $\xi\leqslant2\lambda$ is equivalent to $\beta\leqslant\frac{\lambda}{\sqrt{r_{1}}}$, and hence $s=2\sqrt{r_{1}}$.
Similarly, on $\Gamma_{2}(1)$, the condition $\xi\geqslant c$ is equivalent to $\beta\geqslant \frac{c}{2\sqrt{r_{2}}}$, and hence $s=2\sqrt{r_{2}}$.

For $\frac{\lambda}{\sqrt{r_{1}}}<\beta<\frac{c}{2\sqrt{r_{2}}}$, the same calculation as in case (ii) shows that the intersection with $\Gamma_{\mathrm{m}}(1)$ satisfies \eqref{eq:directional-middle-quadratic}.
As proved there, $Q_{\beta}$ is strictly increasing on $(0,\infty)$ and its unique positive root is $c_{\lambda}(\beta)$. Moreover,
\begin{align*}
  Q_{\beta}\left(\frac{2\lambda}{\beta}\right)=4\left(\frac{\lambda^{2}}{\beta^{2}}-r_{1} \right) <0,
  \qquad Q_{\beta}\left(\frac{c}{\beta}\right)=\frac{c^{2}}{\beta^{2}}-4r_{2} >0.
\end{align*}
Thus this root satisfies $2\lambda<s\beta<c$ and lies on $\Gamma_{\mathrm{m}}(1)$.
At both endpoints, the adjacent branch formulas agree. This proves assertion (iii).

\textbf{(iv)} Assume that $-2\sqrt{r_{1}}<c<2\alpha$.
By Proposition \ref{prop:explicit-wulff-slow-interface} with $t=1$, $\partial\mathcal{W}=\Gamma_{1}(1)\cup\Gamma_{\mathrm{m}}(1)\cup\Gamma_{2}(1)$.
Consider first the intersection of the ray $s\mathbf{n}$ with the sphere of radius $2\sqrt{r_{1}}$.
At $s=2\sqrt{r_{1}}$, we have $\xi=2\sqrt{r_{1}}\beta$.
Hence this point belongs to $\Gamma_{1}(1)$ precisely when $\beta\leqslant\frac{\xi_{*}}{2\sqrt{r_{1}}}$.
In this case, $c_\mathrm{dir}(\mathbf{n})=2\sqrt{r_{1}}$.
Similarly, at the sphere of radius $2\sqrt{r_{2}}$, we have $\xi=2\sqrt{r_{2}}\beta$.
Hence this point belongs to $\Gamma_{2}(1)$ precisely when $\beta\geqslant\frac{c}{2\sqrt{r_{2}}}$.
In this case, $c_\mathrm{dir}(\mathbf{n})=2\sqrt{r_{2}}$.

Before considering the remaining directions, we claim that $\frac{c}{2\sqrt{r_{2}}}>\frac{\xi_{*}}{2\sqrt{r_{1}}}$. Indeed,
\begin{align*}
  \frac{c}{2\sqrt{r_{2}}}-\frac{\xi_{*}}{2\sqrt{r_{1}}}=\frac{\alpha\sqrt{4r_{2}-c^{2}}+c(\sqrt{r_{2}}-\sqrt{r_{1}})}{2r_{2}}>0,
\end{align*}
where the last inequality is immediate if $c\geqslant0$; if $c<0$, then $c^{2}<4r_{1}$ and $-c<2\sqrt{r_{1}}$ imply
\begin{align*}
  \alpha\sqrt{4r_{2}-c^{2}}>2(r_{2}-r_{1})>-c(\sqrt{r_{2}}-\sqrt{r_{1}}).
\end{align*}
Consequently, the remaining directions are precisely those satisfying $\tfrac{\xi_{*}}{2\sqrt{r_{1}}}<\beta<\tfrac{c}{2\sqrt{r_{2}}}$.
For such a direction, the boundary point belongs to neither $\Gamma_{1}(1)$ nor $\Gamma_{2}(1)$.
Since every ray $s\mathbf{n}$ meets $\partial\mathcal{W}$ at a unique point, together with $\partial\mathcal{W}=\Gamma_{1}(1)\cup\Gamma_{\mathrm{m}}(1)\cup\Gamma_{2}(1)$, it must meet $\Gamma_{\mathrm{m}}(1)$.
Substituting $\xi=s\beta$ and $\rho=s\sqrt{1-\beta^{2}}$ into $\xi_{*}\xi+\rho_{*}\rho=4r_{1}$ yields $s(\xi_{*}\beta+\rho_{*}\sqrt{1-\beta^{2}})=4r_{1}$.
Since $s>0$, we have $\xi_{*}\beta+\rho_{*}\sqrt{1-\beta^{2}}>0$, and hence
\begin{align*}
  s=\frac{4r_{1}} {\xi_{*}\beta+\rho_{*}\sqrt{1-\beta^{2}}}.
\end{align*}
By \eqref{eq:tangent-condition}, this expression agrees with the adjacent spherical formulas at the two threshold values. This proves assertion (iv).

\textbf{(v)} Assume that $-2\sqrt{r_{2}}<c\leqslant-2\sqrt{r_{1}}$.
By Proposition \ref{prop:explicit-wulff-slow-interface} with $t=1$,
\begin{align*}
  s\mathbf{n}\in\mathcal{W}\;\iff\; 0\leqslant s\leqslant2\sqrt{r_{2}} \quad\text{and}\quad s\beta\geqslant c.
\end{align*}
If $\beta<\frac{c}{2\sqrt{r_{2}}}$, then $\beta<0$ and $s\beta\geqslant c$ if and only if $s\leqslant\frac{c}{\beta}$.
Since $\frac{c}{\beta}<2\sqrt{r_{2}}$, the largest admissible value of $s$ is $\frac{c}{\beta}$.
If $\beta\geqslant\frac{c}{2\sqrt{r_{2}}}$, then $2\sqrt{r_{2}}\beta\geqslant c$.
Hence $2\sqrt{r_{2}}\mathbf{n}\in\mathcal{W}$, while no point $s\mathbf{n}$ with $s>2\sqrt{r_{2}}$ belongs to $\mathcal{W}$.
Therefore, the largest admissible value of $s$ is $2\sqrt{r_{2}}$.
At the threshold value, the two expressions agree.
This proves assertion (v).

\textbf{(vi)} Assume that $c\leqslant-2\sqrt{r_{2}}$.
By \eqref{eq:wulff-shape-definition} and \eqref{eq:fast-receding-front} with $t=1$,
\begin{align*}
  \partial\mathcal{W}= \left\{\mathbf{x}\in\mathbb{R}^{N} \,\middle|\, |\mathbf{x}|=2\sqrt{r_{2}} \right\}.
\end{align*}
Hence, $c_\mathrm{dir}(\mathbf{n})=2\sqrt{r_{2}}$.
This proves assertion (vi) and completes the proof.
\end{proof}

\section{Discussion and extensions}
\label{section:discussion}

The results above describe the full propagation geometry for a monotone transition between two favorable regions in dimensions $N\geqslant2$.
The variational approach also applies in one dimension and to favorable--unfavorable environments, but the conditions for its reductions and the interpretation of its zero set require separate consideration.
We discuss these differences below.

\subsection{The one-dimensional problem}
\label{subsec:discussion-1d}

Consider \eqref{eq:rd-nd-main} with $N=1$ and $\mathbf{e}=1$, retaining the assumptions $0<r_{1}<r_{2}$ and compactly supported initial data.
The heat-semigroup estimates and comparison argument for the geometric optics limit remain valid in one dimension.
The stopping-game representation and min-max formula also remain valid, as noted in Remark \ref{remark:minmax-lam-conjecture}.

The one-dimensional spreading speeds can be obtained from the axial speeds of the multidimensional problem.
Take $\mathbf{x}=\xi\mathbf{e}$ and write an admissible trajectory as $\gamma(s)=Y(s)\mathbf{e}+{\boldsymbol\eta}(s)$, where ${\boldsymbol\eta}(s)\cdot\mathbf{e}=0$.
The trajectory $Y\mathbf{e}$ preserves the endpoints and the growth term, so
\begin{align}\label{eq:aux-traj-1D}
  A_\gamma(\sigma)=A_{Y\mathbf{e}}(\sigma)+\frac{1}{4}\int_0^\sigma|\dot{\boldsymbol\eta}(s)|^2\,ds\geqslant A_{Y\mathbf{e}}(\sigma),\quad \sigma\in[0,t].
\end{align}
By \eqref{eq:aux-traj-1D}, replacing $\gamma$ by $Y\mathbf{e}$ cannot increase its maximum action, so the one-dimensional min-max value does not exceed the multidimensional value.
The reverse inequality follows because axial trajectories are also admissible in the multidimensional problem, so the two values agree on the axis.
The same argument at $\sigma=t$ gives equality for $J$.
Hence the one-dimensional expansion set is
\begin{align*}
  \mathcal W_{\mathrm{1D}}:=\{x\in\mathbb{R}\mid W(x,1)=0\}=[-s_{-},s_{+}],
\end{align*}
where $s_{-}$ and $s_{+}$ are given by Corollary \ref{cor:main-directional-speed} with $\beta=-1$ and $\beta=1$, respectively.
This recovers the spreading speeds obtained by Lam and Yu \cite{lam2022a} for compactly supported initial data.

However, the set of speeds for which Freidlin's condition holds differs between one and higher dimensions.
This difference already occurs at $c=-2\sqrt{r_{2}}$, where the condition holds for $N\geqslant2$ but fails for $N=1$.
Indeed, the axial agreement of $J$ justifies using \eqref{eq:backward-V-piecewise} with ${\boldsymbol\zeta}=\mathbf0$, giving, for $N=1$,
\begin{align}\label{eq:traject-zero-set}
  J(ct,t)=0,
  \qquad
  \{x\in\mathbb{R}\mid J(x,t)>0\}=(-\infty,ct)\cup(-ct,\infty)\quad \text{for } t>0.
\end{align}
Let $\gamma\in\mathbb{X}^{ct,t}$ satisfy $J(\gamma(s),t-s)>0$ for all $s\in(0,t)$.
Since $J(0,t-s)=-r_{2}(t-s)<0$ for $s\in(0,t)$, the trajectory cannot pass through $0$ before time $t$, as this would contradict $J(\gamma(s),t-s)>0$.
By continuity and $\gamma(0)=ct<0$, we therefore have $\gamma(s)<0$ for all $s\in(0,t)$.
Applying \eqref{eq:traject-zero-set} at time $t-s$ gives either $\gamma(s)<c(t-s)$ or $\gamma(s)>-c(t-s)$.
The second inequality contradicts $\gamma(s)<0$, since $-c(t-s)>0$.
Thus $\gamma(s)<c(t-s)$ and $R(\gamma(s),t-s)=r_{1}$ for all $s\in(0,t)$.
Using the endpoints $\gamma(0)=ct$, $\gamma(t)=0$ and the identity $c^2=4r_{2}$, the Cauchy-Schwarz inequality gives
\begin{align*}
  A_\gamma(t)=\frac{1}{4}\int_{0}^{t}|\dot\gamma(s)|^{2}\,ds-r_{1}t\geqslant\frac{|\gamma(t)-\gamma(0)|^{2}}{4t}-r_{1}t=(r_{2}-r_{1})t>0=J(ct,t).
\end{align*}
This positive lower bound is independent of $\gamma$, so these actions cannot approach $J(ct,t)=0$ and Freidlin's condition fails at $c=-2\sqrt{r_{2}}$ in one dimension.

For $N\geqslant2$, Freidlin's condition holds at the critical speed $c=-2\sqrt{r_{2}}$ by Theorem \ref{thm:main-freidlin-dichotomy}(ii).
At $ct\mathbf{e}$, the orthogonal perturbations in Corollary \ref{cor:backward-freidlin-perturbation} keep the trajectories on the shifting interface, preserve $R=r_{2}$, and ensure $J>0$ at every intermediate time with total action tending to zero.
Nevertheless, the agreement of the one- and multidimensional values on the axis yields $W=\max\{0,J\}$ in one dimension at this speed, despite the failure of Freidlin's condition.

\subsection{Favorable--unfavorable environments}

Suppose that $r\in C(\mathbb{R};\mathbb{R})$ is nondecreasing and satisfies $r(-\infty)=r_{1}<0<r_{2}=r(+\infty)$, with the same assumptions on the initial datum.
The heat-semigroup estimates still yield the growth bound \eqref{eq:growth-bound}.
The upper bound \eqref{eq:stopping-value-upper-bound} becomes $W(\mathbf{x},t)\leqslant\tfrac{|\mathbf{x}|^2}{4t}-r_{1}t$.
The bound in Lemma \ref{lem:bounded-data-value-regularity} becomes $0\leqslant W_0(\mathbf{x},t)\leqslant w_0(\mathbf{x})-r_{1}t$, so the bounded-data arguments apply on finite time intervals.
The contribution of $-R$ over the two connecting intervals in \eqref{eq:traj-eta} is at most $2(-r_{1})\delta$, which is absorbed into $C_K\delta$ in \eqref{eq:path-perturbation-bound}.
The remaining error terms involving $R$ are controlled by $|R|\leqslant\max\{-r_{1},r_{2}\}$.
With these adjustments, the geometric optics limit, stopping-game representation and min-max formula extend to this setting.

Since the min-max minimum is attained, $W(\mathbf{x},t)=0$ gives a minimizing trajectory $\gamma_{*}$ with $A_{\gamma_{*}}(\sigma)\leqslant0$ for all $\sigma\in[0,t]$.
This forces $\mathbf{x}\cdot\mathbf{e}\geqslant ct$, since otherwise $A_{\gamma_{*}}(\sigma)\geqslant-r_{1}\sigma>0$ for small $\sigma>0$.
Besides, $0\geqslant A_{\gamma_{*}}(t)\geqslant \frac{|\mathbf{x}|^2}{4t}-r_{2}t$ gives $|\mathbf{x}|\leqslant2\sqrt{r_{2}}t$.
Conversely, under these two constraints, the straight trajectory $\gamma(s)=(1-\frac{s}{t})\mathbf{x}$ stays in the $r_{2}$-region and satisfies $A_\gamma(\sigma)=(\tfrac{|\mathbf{x}|^2}{4t^2}-r_{2})\sigma\leqslant0$ for $\sigma\in[0,t]$, so $W(\mathbf{x},t)=0$.
Thus, for $c<2\sqrt{r_{2}}$, the zero set of $W$ at $t=1$ is
\begin{align*}
  \mathcal{D}_{c}:=\left\{\mathbf{x}\in\mathbb{R}^{N}\mid W(\mathbf{x},1)=0\right\}
  =\left\{\mathbf{x}\in\mathbb{R}^{N}\mid|\mathbf{x}|\leqslant2\sqrt{r_{2}},\quad\mathbf{x}\cdot\mathbf{e}\geqslant c\right\}.
\end{align*}
On compact subsets of the space--time interior of $\{W=0\}$, $r(\tfrac{\mathbf{x}\cdot\mathbf{e}-ct}{\varepsilon})\to r_{2}>0$ uniformly, so the persistence argument for \eqref{eq:persistence} applies.
Together with decay from $w^\varepsilon\to W$ where $W>0$ and the heat-kernel estimate, this identifies $\mathcal D_c$ as the closed expansion set.
This recovers the truncated-ball geometry in related favorable--unfavorable models \cite{yi2025,vo2026}.

For $c\geqslant2\sqrt{r_{2}}$, the solution tends to zero uniformly in space.
Indeed, choose $L,\kappa>0$ such that $r(z)\leqslant-\kappa$ for $z\leqslant-L$.
Since $r\leqslant r_{2}$ and $g$ is compactly supported, comparison with the linear equation $v_t=\Delta v+r_{2}v$ yields, for some $C_L>0$ and all sufficiently large $t$,
\begin{align}\label{eq:bad-estimate}
  \sup_{\mathbf{x}\cdot\mathbf{e}\geqslant ct-L}u(\mathbf{x},t)\leqslant C_Lt^{-N/2}\exp\left[\left(r_{2}-\frac{c^{2}}{4}\right)t\right]\longrightarrow0.
\end{align}
Given $\epsilon>0$, \eqref{eq:bad-estimate} gives $T>0$ such that $u\leqslant\epsilon$ on $\mathbf{x}\cdot\mathbf{e}=ct-L$ for all $t\geqslant T$.
With $M=\max\{\|g\|_{\infty},r_{2}\}$, comparison with the supersolution $\epsilon+Me^{-\kappa(t-T)}$ gives $u\leqslant\epsilon+Me^{-\kappa(t-T)}$ for $\mathbf{x}\cdot\mathbf{e}\leqslant ct-L$ and $t\geqslant T$.
Letting $t\to\infty$ and then $\epsilon\downarrow0$ proves uniform extinction.

When $r_{1}>0$, the lower-growth region supports propagation and, in the non-Freidlin regime, contributes the ball $D_{1}$ to the convex-hull formula in Theorem \ref{thm:main-freidlin-dichotomy}(iii).
When $r_{1}<0$, persistent expansion is confined to the favorable half-space.


\appendix

\section{Comparison principle}
\label{app:comparison-structure}

\begin{lemma}\label{lem:structural-check}
Let $\Omega= \mathbb{R}^{N}\times(0,T]$ with some given $T\in\mathbb{R}^{+}$ and $\Gamma=\{ (\mathbf{x},t) \in \mathbb{R}^{N} \times [0,T] \mid \mathbf{x}\cdot \mathbf{e}=ct \}$.
The Hamiltonian $H(\mathbf{p}, \mathbf{x}, t)=|\mathbf{p}|^{2}+R(\mathbf{x},t)$ associated with \eqref{eq:def-R-macro} satisfies the following conditions: 
\begin{enumerate}[label=(A\arabic*),font=\normalfont, leftmargin=*]
    \item $H\in C(\mathbb{R}^{N}\times(\Omega\backslash \Gamma))$.
    \item For each $(\mathbf{x}_{0},t_{0})\in (\Omega\backslash \Gamma)\cap (B_{R}(\mathbf{0})\times(0,T))$, there exists a constant $\delta_{0}>0$ such that
    \begin{align*}
  H(\mathbf{p},\mathbf{y},t)-H(\mathbf{p},\mathbf{x},t) \leqslant \omega_{R}\left(|\mathbf{x}-\mathbf{y}|\left(1+|\mathbf{p}|\right)\right)
    \end{align*}
    for $\mathbf{p},\mathbf{x},\mathbf{y},t$ satisfying $\mathbf{p}\in \mathbb{R}^{N}$ and $|(\mathbf{x},t)-(\mathbf{x}_{0},t_{0})|+|(\mathbf{y},t)-(\mathbf{x}_{0},t_{0})|<\delta_{0}$, where $\omega_{R}$ is a modulus of continuity.
    \item For each $(\mathbf{x}_{0},t_{0})\in \Omega\cap \Gamma\cap (B_{R}(\mathbf{0})\times(0,T))$, there exist a constant $\delta_{0}>0$ and a unit vector $(\mathbf{k}_{0},h_{0})\in \mathbb{R}^{N}\times\mathbb{R}$ such that
    \begin{align*}
  H^{*}(\mathbf{p},\mathbf{y},s)-H_{*}(\mathbf{p},\mathbf{x},t) \leqslant \omega_{R}((|\mathbf{x}-\mathbf{y}|+|t-s|)(1+|\mathbf{p}|))
    \end{align*}
    for all $\mathbf{p}\in \mathbb{R}^{N}$ and $\mathbf{x},\mathbf{y},t,s$ satisfying
    \begin{align*}
  |(\mathbf{x},t)-(\mathbf{x}_{0},t_{0})|+|(\mathbf{y},s)-(\mathbf{x}_{0},t_{0})| <\delta_{0}
  \quad \text{and} \quad \left\|\frac{(\mathbf{x}-\mathbf{y},t-s)}{|(\mathbf{x}-\mathbf{y},t-s)|} -(\mathbf{k}_{0},h_{0})\right\| <\delta_{0}.
    \end{align*}
    \item There exists $M\geqslant 0$ such that for each $\theta \in[0,1)$ and $\mathbf{x}_{0}\in \mathbb{R}^{N}$, there exist constants $\bar\epsilon(\theta,\mathbf{x}_{0})>0$ and $\bar C(\theta,\mathbf{x}_{0})>0$ such that
    \begin{align*}
  H\left(\theta \mathbf{p}-\frac{\epsilon(\mathbf{x}-\mathbf{x}_{0})} {|\mathbf{x}-\mathbf{x}_{0}|^{2}+1},\mathbf{x},t\right)-M
  \leqslant \theta \left(H(\mathbf{p},\mathbf{x},t)-M\right) +\epsilon \bar C(\theta,\mathbf{x}_{0})
    \end{align*}
    for all $\mathbf{p}\in \mathbb{R}^{N}$, $(\mathbf{x},t)\in \Omega$ and $\epsilon\in[0,\bar\epsilon(\theta,\mathbf{x}_{0})]$.
\end{enumerate}
\end{lemma}

\begin{proof}
The singularity set is the moving hyperplane $\Gamma=\{ (\mathbf{x},t) \in \mathbb{R}^{N} \times [0,T] \mid \mathbf{x}\cdot \mathbf{e}=ct \}$.
Since $H$ is independent of $(\mathbf{x},t)$ on each side of $\Gamma$, (A1) and (A2) hold.

To verify the local monotonicity condition (A3), fix any $(\mathbf{x}_{0}, t_{0}) \in \Gamma$ and choose $(\mathbf{k}_{0}, h_{0})=(\mathbf{e}, 0) \in \mathbb{R}^{N} \times \mathbb{R}$ and set the constant $\delta_{0}$ as $\delta_{0}=\tfrac{1}{4(1+|c|)}>0$.
Consider any two distinct points $(\mathbf{x}, t), (\mathbf{y}, s)$ satisfying $\|(\mathbf{x},t)-(\mathbf{x}_{0}, t_{0})\|+\|(\mathbf{y},s)-(\mathbf{x}_{0}, t_{0})\|<\delta_{0}$ and the directional restriction:
\begin{align}
\label{eq:directional-restriction}
  \left\| \frac{(\mathbf{x}-\mathbf{y}, t-s)}{\|(\mathbf{x}-\mathbf{y}, t-s)\|}-(\mathbf{e}, 0) \right\|<\delta_{0}.
\end{align}
Let $\Delta=\|(\mathbf{x}-\mathbf{y}, t-s)\|>0$.
\eqref{eq:directional-restriction} implies
\begin{align*}
  \mathbf{x}-\mathbf{y}=\Delta (\mathbf{e}+\epsilon_\mathbf{x}), \quad t-s=\Delta \epsilon_{t},
\end{align*}
where the perturbation vector $(\epsilon_\mathbf{x}, \epsilon_{t}) \in \mathbb{R}^{N} \times \mathbb{R}$ satisfies $\|(\epsilon_\mathbf{x}, \epsilon_{t})\|<\delta_{0}$.
This bound guarantees that $\|\epsilon_\mathbf{x}\|<\delta_{0}$ and $|\epsilon_{t}|<\delta_{0}$.

The difference in the signed distance of $(\mathbf{y}, s)$ and $(\mathbf{x}, t)$ is given by
\begin{align*}
  (\mathbf{y} \cdot \mathbf{e}-c s)-(\mathbf{x} \cdot \mathbf{e}-c t) &=-(\mathbf{x}-\mathbf{y}) \cdot \mathbf{e}+c(t-s)
 =-\Delta \left(1+\epsilon_\mathbf{x} \cdot \mathbf{e}-c \epsilon_{t} \right).
\end{align*}
By the perturbation bounds, we have
\begin{align*}
  \left| \epsilon_\mathbf{x}\cdot\mathbf{e}-c\epsilon_{t} \right| \leqslant \|\epsilon_\mathbf{x}\|+|c|\,|\epsilon_{t}|<(1+|c|)\delta_{0}=\frac{1}{4}.
\end{align*}
Hence, $1+\epsilon_\mathbf{x}\cdot\mathbf{e}-c\epsilon_{t}>\frac{3}{4}$. Since $\Delta>0$, we have
\begin{align}
\label{eq:strict-descent-exact}
  \mathbf{y} \cdot \mathbf{e}-c s<\mathbf{x} \cdot \mathbf{e}-c t.
\end{align}

With \eqref{eq:strict-descent-exact} established, for any $\mathbf{p}\in \mathbb{R}^{N}$, we consider $H^{*}(\mathbf{p}, \mathbf{y}, s)-H_{*}(\mathbf{p}, \mathbf{x}, t)=R^{*}(\mathbf{y}, s)-R_{*}(\mathbf{x}, t)$ by considering the following two cases:
\begin{itemize}[leftmargin=*]
\item \textbf{Case 1:} $R_{*}(\mathbf{x}, t)=r_{2}$.
Since $r_{2}$ is the global supremum of the function $R$, the upper envelope trivially satisfies $R^{*}(\mathbf{y}, s) \leqslant r_{2}$ for all $(\mathbf{y},s)$.
Thus, $R^{*}(\mathbf{y}, s)-R_{*}(\mathbf{x}, t) \leqslant 0$.
\item \textbf{Case 2:} $R_{*}(\mathbf{x}, t)=r_{1}$.
The lower envelope $R_{*}$ evaluates to the global infimum $r_{1}$ if and only if $\mathbf{x}\cdot \mathbf{e}\leqslant ct$.
Combining with \eqref{eq:strict-descent-exact}, we have $\mathbf{y}\cdot\mathbf{e}-cs<0$.
This further implies that $(\mathbf{y}, s)$ lies in the open interior of the domain where $R \equiv r_{1}$.
Over this open set, it holds that $R^{*}(\mathbf{y}, s)=r_{1}$.
Consequently, $R^{*}(\mathbf{y}, s)-R_{*}(\mathbf{x}, t) =0$.
\end{itemize}
In all cases, we obtain $H^{*}(\mathbf{p}, \mathbf{y}, s)-H_{*}(\mathbf{p}, \mathbf{x}, t) \leqslant 0$.
Setting the modulus of continuity $\omega_{R}(\rho) =\rho$, condition (A3) is verified.

We verify (A4) with $M=r_{2}$. Put $\mathbf{q}=\frac{\mathbf{x}-\mathbf{x}_{0}} {|\mathbf{x}-\mathbf{x}_{0}|^{2}+1}$.
Then $|\mathbf{q}|\leqslant \frac{1}{2}$.
For $\theta\in[0,1)$ and $\epsilon\in[0,1)$, by $R(\mathbf{x},t)\leqslant r_{2}$ and Young's inequality, we obtain
\begin{align*}
  H(\theta\mathbf{p}-\epsilon\mathbf{q},\mathbf{x},t)-r_{2} &= \theta^{2}|\mathbf{p}|^{2} -2\theta\epsilon\mathbf{p}\cdot\mathbf{q} +\epsilon^{2}|\mathbf{q}|^{2} +R(\mathbf{x},t)-r_{2}\\
  &\leqslant \theta|\mathbf{p}|^{2} +\frac{\epsilon^{2}|\mathbf{q}|^{2}}{1-\theta} +R(\mathbf{x},t)-r_{2}\\
  &\leqslant \theta\left(H(\mathbf{p},\mathbf{x},t)-r_{2}\right) +\frac{\epsilon^{2}}{4(1-\theta)}\\
  &\leqslant \theta\left(H(\mathbf{p},\mathbf{x},t)-r_{2}\right) +\epsilon\frac{1}{4(1-\theta)}.
\end{align*}
Thus (A4) holds with $\bar\epsilon(\theta,\mathbf{x}_{0})=1$ and $\bar C(\theta,\mathbf{x}_{0})=1/[4(1-\theta)]$. This completes the proof.
\end{proof}

\begin{corollary}[Comparison principle]
\label{coro:comparison-principle}
Let $0\leqslant \tau<T<\infty$ and set $Q_{\tau,T}=\mathbb{R}^{N}\times(\tau,T]$.
Suppose that $U\in USC(\mathbb{R}^{N}\times[\tau,T])$ and $V\in LSC(\mathbb{R}^{N}\times[\tau,T])$ are, respectively, a viscosity subsolution and a viscosity supersolution of
\begin{align}
\label{eq:main-variational-inequa}
  \min\left\{\partial_{t} w+|Dw|^{2}+R(\mathbf{x},t), w \right\}=0 \quad\text{in }Q_{\tau,T}.
\end{align}
If
\begin{align*}
  U(\mathbf{x},\tau)\leqslant V(\mathbf{x},\tau) \quad\text{for all }\mathbf{x}\in\mathbb{R}^{N},
\end{align*}
then
\begin{align*}
  U(\mathbf{x},t)\leqslant V(\mathbf{x},t) \quad\text{for all }(\mathbf{x},t)\in \mathbb{R}^{N}\times(\tau,T].
\end{align*}
\end{corollary}

\begin{proof}
By Lemma \ref{lem:structural-check}, the Hamiltonian $H(\mathbf{p},\mathbf{x},t)=|\mathbf{p}|^{2}+R(\mathbf{x},t)$ satisfies conditions (A1)--(A4) in \cite[Appendix A]{liu2021a}.
These conditions remain valid after restricting the time interval to $(\tau,T]$.
The variational inequality \eqref{eq:main-variational-inequa} is equation (A.1) in \cite{liu2021a} with the parameter $L=0$.
The conclusion therefore follows directly from \cite[Theorem A.1]{liu2021a}.
\end{proof}

\section{Technical proofs for the stopping game with bounded initial data}
\label{app:bounded-data}

This appendix provides the deferred proofs of the bounded-data regularity estimate, the Dynamic Programming Principle, and the viscosity characterization of the bounded-data stopping value.

\subsection{Bounds and continuity of the value function}
\label{app:bounded-data-regularity}

\begin{proof}[Proof of Lemma \ref{lem:bounded-data-value-regularity}]
We complete the proof by the following steps.

\textbf{Step 1.} We show the bounds for $W_{0}$.
The admissible stopping rule $\vartheta\equiv0$ ensures $W_{0}\geqslant0$.
We now prove the pointwise upper bound.
Fix $\vartheta\in\Theta$ and test the inner infimum by the admissible constant trajectory $\gamma_\mathbf{x}(s)\equiv\mathbf{x}$.
If $\vartheta[\gamma_\mathbf{x}]<t$, then the terminal cost does not appear and
\begin{align*}
  \mathcal P_{0}(\mathbf{x},t;\vartheta,\gamma_\mathbf{x})=-\int_{0}^{\vartheta[\gamma_\mathbf{x}]}R(\mathbf{x},t-s)\,ds \leqslant0 \leqslant w_{0}(\mathbf{x}).
\end{align*}
If $\vartheta[\gamma_\mathbf{x}]\geqslant t$, then
\begin{align*}
  \mathcal P_{0}(\mathbf{x},t;\vartheta,\gamma_\mathbf{x})=-\int_{0}^{t}R(\mathbf{x},t-s)\,ds+w_{0}(\mathbf{x}) \leqslant w_{0}(\mathbf{x}).
\end{align*}
Thus, for all fixed $\vartheta\in\Theta$,
\begin{align*}
  \inf_{\eta\in\mathbb{A}^{\mathbf{x}}} \mathcal P_{0}(\mathbf{x},t;\vartheta,\eta)
  \leqslant \mathcal P_{0}(\mathbf{x},t;\vartheta,\gamma_\mathbf{x}) \leqslant w_{0}(\mathbf{x}).
\end{align*}
Taking the supremum over $\vartheta\in\Theta$, we obtain
\begin{align}
\label{eq:W0-upper-bound}
  W_{0}(\mathbf{x},t)\leqslant w_{0}(\mathbf{x}), \quad (\mathbf{x},t)\in\mathbb{R}^{N}\times\mathbb{R}^{+}.
\end{align}

\textbf{Step 2.} We next prove the continuity of $W_{0}$.
Fix a compact cylinder $K=\overline{B_{A}(\mathbf{0})}\times[\tau,T]$, $0<\tau<T<\infty$, and let $(\mathbf{x},t),(\mathbf{y},\sigma)\in K$. Set $d=|\mathbf{x}-\mathbf{y}|+|t-\sigma|$.
For $0<d<\min\{1,\tau^{2}/16\}$, set $\delta=\sqrt d$.
Use the transformation $\mathcal T$ constructed in \eqref{eq:traj-transform}, together with the corresponding pullback of stopping times
\begin{align*}
  \vartheta[\gamma]=\Pi\left(\tilde\vartheta[\mathcal T\gamma]\right), \quad \gamma\in\mathbb{X},
\end{align*}
where $\Pi$ is given by \eqref{eq:pullback-time}.
For every $\tilde\vartheta\in\Theta$, the preceding argument implies $\vartheta\in\Theta$ and $\vartheta[\gamma]\in[0,t]$ for every $\gamma\in\mathbb{X}$.

Let $\epsilon\in(0,1)$ and choose an $\epsilon$-optimal trajectory $\gamma_{\epsilon}\in\mathbb{A}^{\mathbf{x}}$ for the inner minimization associated with the stopping time $\vartheta$, i.e.,
\begin{align*}
  \mathcal P_{0} (\mathbf{x},t;\vartheta,\gamma_{\epsilon}) \leqslant \inf_{\gamma\in\mathbb A^\mathbf{x}} \mathcal P_{0}(\mathbf{x},t;\vartheta,\gamma) +\epsilon
  \leqslant W_{0}(\mathbf{x},t)+\epsilon \leqslant \|w_{0}\|_{L^{\infty}}+\epsilon.
\end{align*}
Using $R\leqslant r_{2}$ and $w_{0}\geqslant 0$, it holds that
\begin{align}
\label{eq:W0-near-minimizer-energy}
  \int_{0}^{t\wedge\vartheta[\gamma_{\epsilon}]} |\dot\gamma_{\epsilon}(s)|^{2}\,ds \leqslant 4\left(\|w_{0}\|_{L^{\infty}}+r_{2}T+1\right).
\end{align}

Set $\eta_{\epsilon}=\mathcal T\gamma_{\epsilon}$.
Since $\gamma_{\epsilon}(0)=\mathbf{x}$, we have $\eta_{\epsilon}(0)=\mathbf{y}$ and $\eta_{\epsilon}\in\mathbb{A}^{\mathbf{y}}$.
Moreover, $\eta_{\epsilon}(\sigma)=\gamma_{\epsilon}(t)$, and the estimates established in the proof of Proposition \ref{prop:stopping-game-basic-properties} for the three nonconstant branches on $[0,\sigma]$ remain valid.
The constant extension on $(\sigma,\infty)$ does not enter the payoff. Define $\tau_{\epsilon}=\tilde\vartheta[\eta_{\epsilon}]$.
If $\tau_{\epsilon}\geqslant\sigma-\delta$, then the definition of $\Pi$ implies $\vartheta[\gamma_{\epsilon}] =\Pi(\tau_{\epsilon})=t$, and \eqref{eq:W0-near-minimizer-energy} bounds the kinetic energy term of $\gamma_{\epsilon}$ on $[0,t]$.
If $\tau_{\epsilon}<\sigma-\delta$, then the payoff contains no contribution from the connecting interval $[\sigma-\delta,\sigma]$.
Therefore, repeating the estimates leading to \eqref{eq:path-perturbation-bound}, with \eqref{eq:W0-near-minimizer-energy} in place of \eqref{eq:energy-bound-for-near-minimizer}, and using the preceding case distinction, we obtain
\begin{align*}
  \int_{0}^{\sigma\wedge\tau_{\epsilon}} L(\dot\eta_{\epsilon}(s),\eta_{\epsilon}(s),\sigma-s)\,ds
  \leqslant& \int_{0}^{\vartheta[\gamma_{\epsilon}]} L(\dot\gamma_{\epsilon}(s),\gamma_{\epsilon}(s),t-s)\,ds+C_{K}\left(d+\delta+\frac{d^{2}+\delta^{2}}{\delta} \right)\\
  \leqslant& \int_{0}^{\vartheta[\gamma_{\epsilon}]} L(\dot\gamma_{\epsilon}(s),\gamma_{\epsilon}(s),t-s)\,ds +C_{K}\sqrt d,
\end{align*}
where the constant $C_{K}>0$ depends only on $K,r_{2},c$, $N$ and $\|w_{0}\|_{L^{\infty}}$. Moreover,
\begin{align*}
  \mathbbm{1}_{\{\tilde\vartheta[\eta_{\epsilon}]\geqslant\sigma\}} w_{0}(\eta_{\epsilon}(\sigma))
  \leqslant \mathbbm{1}_{\{\vartheta[\gamma_{\epsilon}]\geqslant t\}} w_{0}(\gamma_{\epsilon}(t)).
\end{align*}
Indeed, if $\tilde\vartheta[\eta_{\epsilon}]\geqslant\sigma$, then $\vartheta[\gamma_{\epsilon}]=t$ and $\eta_{\epsilon}(\sigma)=\gamma_{\epsilon}(t)$; if $\tilde\vartheta[\eta_{\epsilon}]<\sigma$, the left-hand side is zero, whereas the right-hand side is nonnegative since $w_{0}\geqslant0$. Consequently,
\begin{align*}
  \mathcal P_{0}(\mathbf{y},\sigma; \tilde\vartheta,\eta_{\epsilon}) \leqslant \mathcal P_{0}(\mathbf{x},t; \vartheta,\gamma_{\epsilon})+C_{K}\sqrt d.
\end{align*}
Since $\eta_{\epsilon}\in\mathbb{A}^{\mathbf{y}}$, for each fixed $\tilde\vartheta\in\Theta$ we have
\begin{align*}
  \inf_{\eta\in\mathbb{A}^{\mathbf{y}}} \mathcal P_{0}(\mathbf{y},\sigma;\tilde\vartheta,\eta)
  \leqslant \mathcal P_{0}(\mathbf{y},\sigma; \tilde\vartheta,\eta_{\epsilon}).
\end{align*}
Combining this with the preceding estimate, using the $\epsilon$-optimality of $\gamma_{\epsilon}$ and the fact that the pullback stopping rule $\vartheta$ belongs to $\Theta$, and then taking the supremum over $\tilde\vartheta\in\Theta$, we obtain
\begin{align*}
  W_{0}(\mathbf{y},\sigma) \leqslant W_{0}(\mathbf{x},t)+\epsilon+C_{K}\sqrt d.
\end{align*}
Letting $\epsilon\downarrow0$ yields
\begin{align*}
  W_{0}(\mathbf{y},\sigma) \leqslant W_{0}(\mathbf{x},t)+C_{K}\sqrt d.
\end{align*}
Interchanging $(\mathbf{x},t)$ and $(\mathbf{y},\sigma)$ proves that $W_{0}\in C(\mathbb{R}^{N}\times\mathbb{R}^{+})$.
\end{proof}

\subsection{Dynamic Programming Principle}
\label{app:bounded-data-dpp}

The proof uses concatenation of trajectories and the property of stopping times; compare this with the strategy concatenation argument for differential games in \cite[Section 3]{evans1984}.

\begin{proof}[Proof of Proposition \ref{prop:bounded-data-dpp}]
Denote the right-hand side of \eqref{eq:bounded-data-dpp} by $\mathcal D_{h}W_{0}(\mathbf{x},t)$.

\textbf{Step 1.} We show first that
\begin{align}
\label{eq:inequality-1}
  W_{0}(\mathbf{x},t)\leqslant\mathcal D_{h}W_{0}(\mathbf{x},t), \quad (\mathbf{x},t)\in \mathbb{R}^{N}\times\mathbb{R}^{+}.
\end{align}
Fix $\vartheta\in\Theta$.
For $\gamma^{0}\in\mathbb{X}$ with $\gamma^{0}(0)=\mathbf{x}$ and $\gamma^{1}\in\mathbb{A}^{\gamma^{0}(h)}$, define the connecting trajectory:
\begin{align*}
  (\gamma^{0}*\gamma^{1})(s) =
  \begin{cases}
    \gamma^{0}(s), &s\in[0,h],\\
    \gamma^{1}(s-h), &s\in(h,\infty).
  \end{cases}
\end{align*}
Then $\gamma^{0}*\gamma^{1}\in\mathbb{A}^{\mathbf{x}}$.
For fixed $\gamma^{0}$, the stopping-time property implies that $h\wedge\vartheta[\gamma^{0}*\gamma^{1}]$ and the event $\{\vartheta[\gamma^{0}*\gamma^{1}]>h\}$ depend only on $\gamma^{0}$.
Define the tail rule $\widehat\vartheta^{\gamma^{0}}:\mathbb{X}\to[0,\infty]$ by
\begin{align*}
  \widehat\vartheta^{\gamma^{0}}[\gamma^{1}] =
  \begin{cases}
    \max\{\vartheta[\gamma^{0}*\gamma^{1}]-h,0\}, &\gamma^{1}(0)=\gamma^{0}(h),\\
    0, &\gamma^{1}(0)\ne\gamma^{0}(h).
  \end{cases}
\end{align*}

We first verify that $\widehat\vartheta^{\gamma^{0}}$ is a stopping time.
Let $\gamma^{1},\widehat\gamma^{1}\in\mathbb{X}$ agree on $[0,r]$, and assume that $\widehat\vartheta^{\gamma^{0}}[\gamma^{1}]\leqslant r$.
If $\gamma^{1}(0)\ne\gamma^{0}(h)$, then also $\widehat\gamma^{1}(0)\ne\gamma^{0}(h)$, and both tail stopping times are zero.
It remains to consider the case $\gamma^{1}(0)=\gamma^{0}(h)$.
Then $\widehat\gamma^{1}(0)=\gamma^{0}(h)$, and the trajectories $\gamma^{0}*\gamma^{1}$ and $\gamma^{0}*\widehat\gamma^{1}$ agree on $[0,h+r]$.
If $\vartheta[\gamma^{0}*\gamma^{1}]\leqslant h$, the stopping-time property of $\vartheta$ implies
\begin{align*}
  \vartheta[\gamma^{0}*\widehat\gamma^{1}]=\vartheta[\gamma^{0}*\gamma^{1}],
\end{align*}
and hence both tail stopping times are zero.
If $\vartheta[\gamma^{0}*\gamma^{1}]>h$, then
\begin{align*}
  \vartheta[\gamma^{0}*\gamma^{1}]=h+\widehat\vartheta^{\gamma^{0}}[\gamma^{1}] \leqslant h+r.
\end{align*}
Applying the stopping-time property of $\vartheta$ on $[0,h+r]$, we obtain
\begin{align*}
  \vartheta[\gamma^{0}*\widehat\gamma^{1}]=\vartheta[\gamma^{0}*\gamma^{1}],
\end{align*}
and therefore $\widehat\vartheta^{\gamma^{0}}[\widehat\gamma^{1}] =\widehat\vartheta^{\gamma^{0}}[\gamma^{1}]$.
Thus $\widehat\vartheta^{\gamma^{0}}\in\Theta$.

If $\vartheta[\gamma^{0}*\gamma^{1}]>h$, then, by the definition of $\widehat\vartheta^{\gamma^{0}}$,
\begin{align*}
  \widehat\vartheta^{\gamma^{0}}[\gamma^{1}] =\vartheta[\gamma^{0}*\gamma^{1}]-h,
  \quad t\wedge \vartheta[\gamma^{0}*\gamma^{1}]=h+\left((t-h)\wedge \widehat\vartheta^{\gamma^{0}}[\gamma^{1}]\right).
\end{align*}
Moreover,
\begin{align*}
  \vartheta[\gamma^{0}*\gamma^{1}]\geqslant t \;\iff\; \widehat\vartheta^{\gamma^{0}}[\gamma^{1}] \geqslant t-h
  \quad \text{and} \quad (\gamma^{0}*\gamma^{1})(t)=\gamma^{1}(t-h).
\end{align*}
Splitting the integral at $h$ and using the change of variables $r=s-h$ on the second interval, we obtain
\begin{align}
\label{eq:infimum-1-P}
  \mathcal P_{0}(\mathbf{x},t;\vartheta,\gamma^{0}*\gamma^{1})
 =\int_{0}^{h} L(\dot\gamma^{0}(s),\gamma^{0}(s),t-s)\,ds+\mathcal P_{0} \left(\gamma^{0}(h),t-h; \widehat\vartheta^{\gamma^{0}},\gamma^{1} \right).
\end{align}
If $\vartheta[\gamma^{0}*\gamma^{1}]\leqslant h$, then $\vartheta[\gamma^{0}*\gamma^{1}]<t$ and
\begin{align}
\label{eq:infimum-2-P}
  \mathcal P_{0}(\mathbf{x},t;\vartheta,\gamma^{0}*\gamma^{1})=\int_{0}^{\vartheta[\gamma^{0}*\gamma^{1}]} L(\dot\gamma^{0}(s),\gamma^{0}(s),t-s)\,ds.
\end{align}
Therefore, by the definition of $W_{0}$, together with \eqref{eq:infimum-1-P} and \eqref{eq:infimum-2-P}, we have
\begin{equation}
\label{eq:infimum-sum-P}
  \begin{aligned}
    &\inf_{\tilde\gamma^{1}\in\mathbb{A}^{\gamma^{0}(h)}} \mathcal P_{0}(\mathbf{x},t;\vartheta,\gamma^{0}*\tilde\gamma^{1})\\
    \leqslant& \int_{0}^{h\wedge\vartheta[\gamma^{0}*\gamma^{1}]} L(\dot\gamma^{0}(s),\gamma^{0}(s),t-s)\,ds+\mathbbm{1}_{\{\vartheta[\gamma^{0}*\gamma^{1}]>h\}} W_{0}(\gamma^{0}(h),t-h).
  \end{aligned}
\end{equation}
Since the set of all connecting trajectories $\gamma^{0}*\tilde\gamma^{1}$ with $\gamma^{0}\in\mathbb A^\mathbf{x}$ and $\tilde\gamma^{1}\in\mathbb A^{\gamma^{0}(h)}$ coincides with $\mathbb A^\mathbf{x}$, for the fixed stopping time $\vartheta\in \Theta$, it follows from \eqref{eq:infimum-sum-P} that
\begin{align*}
  &\inf_{\gamma\in\mathbb A^\mathbf{x}} \mathcal P_{0}(\mathbf{x},t;\vartheta,\gamma)\\
  =& \inf_{\gamma^0\in\mathbb A^\mathbf{x}} \inf_{\tilde\gamma^{1}\in \mathbb A^{\gamma^{0}(h)}} \mathcal P_{0} (\mathbf{x},t;\vartheta,\gamma^{0}*\tilde\gamma^{1})\\
  \leqslant& \inf_{\gamma^0\in\mathbb A^\mathbf{x}} \left\{\int_{0}^{h\wedge\vartheta[\gamma^{0}]} L(\dot\gamma^{0}(s),\gamma^{0}(s),t-s)\,ds+\mathbbm{1}_{\{\vartheta[\gamma^{0}]>h\}} W_{0}(\gamma^{0}(h),t-h) \right\}\\
  =& \inf_{\gamma\in\mathbb A^\mathbf{x}} \left\{\int_{0}^{h\wedge\vartheta[\gamma]} L(\dot\gamma(s),\gamma(s),t-s)\,ds+\mathbbm{1}_{\{\vartheta[\gamma]>h\}} W_{0}(\gamma(h),t-h) \right\}.
\end{align*}
Taking the supremum over $\vartheta\in\Theta$ yields \eqref{eq:inequality-1}.

\textbf{Step 2.} We now prove
\begin{align}
\label{eq:inequality-2}
  W_{0}(\mathbf{x},t) \geqslant \mathcal D_{h}W_{0}(\mathbf{x},t), \quad (\mathbf{x},t)\in \mathbb{R}^{N}\times\mathbb{R}^{+}.
\end{align}
Let $\epsilon>0$ and choose $\vartheta^{0}\in\Theta$ such that
\begin{align}
\label{eq:sub-optimal-stopping-obstacle-hje}
  \mathcal D_{h}^{\vartheta^{0}}W_{0}(\mathbf{x},t) \geqslant \mathcal D_{h}W_{0}(\mathbf{x},t)-\epsilon,
\end{align}
where $\mathcal D_{h}^{\vartheta}W_{0}$ denotes the value in \eqref{eq:bounded-data-dpp} with the stopping time $\vartheta$ fixed.
For every $\mathbf{y}\in\mathbb{R}^{N}$, choose $\vartheta^\mathbf{y}\in\Theta$ such that
\begin{align}
\label{eq:stopping-control}
  W_{0}(\mathbf{y},t-h) \leqslant \inf_{\gamma^{1}\in\mathbb{A}^{\mathbf{y}}} \mathcal P_{0}(\mathbf{y},t-h;\vartheta^\mathbf{y},\gamma^{1})+\epsilon.
\end{align}
Given $\gamma\in\mathbb{X}$, set $\gamma_{h}(\sigma)=\gamma(h+\sigma)$ for $\sigma\geqslant0$ and define
\begin{align}
\label{eq:new-defined-stopping}
  \vartheta^{\epsilon}[\gamma] =
  \begin{cases}
    \vartheta^{0}[\gamma], &\vartheta^{0}[\gamma]\leqslant h,\\
    h+\vartheta^{\gamma(h)}[\gamma_{h}], &\vartheta^{0}[\gamma]>h.
  \end{cases}
\end{align}

We verify that $\vartheta^{\epsilon}\in\Theta$.
Let $\gamma,\widehat\gamma\in\mathbb{X}$ agree on $[0,s]$, and assume that $\vartheta^{\epsilon}[\gamma]\leqslant s$.
If $\vartheta^{0}[\gamma]\leqslant h$, then $\vartheta^{\epsilon}[\gamma]=\vartheta^{0}[\gamma]\leqslant s$.
Since $\vartheta^{0}$ is a stopping time, we have $\vartheta^{0}[\widehat\gamma]=\vartheta^{0}[\gamma]$.
In particular, $\vartheta^{0}[\widehat\gamma]\leqslant h$, and hence $\vartheta^{\epsilon}[\widehat\gamma]=\vartheta^{\epsilon}[\gamma]$.
If $\vartheta^{0}[\gamma]>h$, then
\begin{align*}
  h+\vartheta^{\gamma(h)}[\gamma_{h}]=\vartheta^{\epsilon}[\gamma] \leqslant s.
\end{align*}
Thus $s\geqslant h$ and $\gamma(h)=\widehat\gamma(h)$.
We claim that $\vartheta^{0}[\widehat\gamma]>h$.
Indeed, if $\vartheta^{0}[\widehat\gamma]\leqslant h$, then the stopping-time property of $\vartheta^{0}$, applied on $[0,h]$, implies
\begin{align*}
  \vartheta^{0}[\gamma]=\vartheta^{0}[\widehat\gamma]\leqslant h,
\end{align*}
which is a contradiction.
Therefore the claim holds, so the second branch in \eqref{eq:new-defined-stopping} applies to both $\gamma$ and $\widehat\gamma$.
Moreover, by \eqref{eq:new-defined-stopping} and $\vartheta^{\epsilon}[\gamma]\leqslant s$, we have
\begin{align}
\label{eq:stopping-tail-bound}
  \vartheta^{\gamma(h)}[\gamma_{h}] \leqslant s-h.
\end{align}
Since $\gamma$ and $\widehat\gamma$ agree on $[0,s]$ and $s\geqslant h$, their tail trajectories $\gamma_{h}$ and $\widehat\gamma_{h}$ agree on $[0,s-h]$.
Combining this with \eqref{eq:stopping-tail-bound} and $\vartheta^{\gamma(h)}\in\Theta$, the stopping-time property implies
\begin{align*}
  \vartheta^{\gamma(h)}[\widehat\gamma_{h}]=\vartheta^{\gamma(h)}[\gamma_{h}].
\end{align*}
Finally, since $\widehat\gamma(h)=\gamma(h)$, the second branch in \eqref{eq:new-defined-stopping} yields
\begin{align*}
  \vartheta^{\epsilon}[\widehat\gamma]=h+\vartheta^{\widehat\gamma(h)}[\widehat\gamma_{h}]=h+\vartheta^{\gamma(h)}[\gamma_{h}]=\vartheta^{\epsilon}[\gamma].
\end{align*}
Therefore $\vartheta^{\epsilon}\in\Theta$.

We now claim that, for all $\gamma\in\mathbb{A}^{\mathbf{x}}$,
\begin{align}
\label{eq:sub-optimal-W0}
  \mathcal P_{0}(\mathbf{x},t;\vartheta^{\epsilon},\gamma)
  \geqslant \int_{0}^{h\wedge\vartheta^{0}[\gamma]} L(\dot\gamma(s),\gamma(s),t-s)\,ds+\mathbbm{1}_{\{\vartheta^{0}[\gamma]>h\}} W_{0}(\gamma(h),t-h)-\epsilon.
\end{align}
Indeed, if $\vartheta^{0}[\gamma]\leqslant h$, then it follows from $h<t$ that $\vartheta^{\epsilon}[\gamma]=\vartheta^{0}[\gamma]\leqslant h<t$ and
\begin{align*}
  \mathcal P_{0}(\mathbf{x},t;\vartheta^{\epsilon},\gamma)=\int_{0}^{h\wedge\vartheta^{0}[\gamma]} L(\dot\gamma(s),\gamma(s),t-s)\,ds.
\end{align*}
Since $\mathbbm{1}_{\{\vartheta^{0}[\gamma]>h\}}=0$, the claim follows in this case.
If $\vartheta^{0}[\gamma]>h$, then $\vartheta^{\epsilon}[\gamma]=h+\vartheta^{\gamma(h)}[\gamma_{h}]$. Thus
\begin{align*}
  \vartheta^{\epsilon}[\gamma]\geqslant t \;\iff\; \vartheta^{\gamma(h)}[\gamma_{h}]\geqslant t-h
  \quad \text{and} \quad t\wedge\vartheta^{\epsilon}[\gamma]=h+\left((t-h)\wedge\vartheta^{\gamma(h)}[\gamma_{h}]\right).
\end{align*}
Using the change of variables $\sigma=s-h$ on the interval after $h$, we obtain
\begin{align*}
  \mathcal P_{0}(\mathbf{x},t;\vartheta^{\epsilon},\gamma) =& \int_{0}^{h} L(\dot\gamma(s),\gamma(s),t-s)\,ds\\
  &+ \int_{0}^{(t-h)\wedge\vartheta^{\gamma(h)}[\gamma_{h}]} L(\dot\gamma_{h}(\sigma),\gamma_{h}(\sigma),t-h-\sigma)\,d\sigma+ \mathbbm{1}_{\{\vartheta^{\gamma(h)}[\gamma_{h}]\geqslant t-h\}}w_{0}(\gamma_{h}(t-h))\\
  =& \int_{0}^{h} L(\dot\gamma(s),\gamma(s),t-s)\,ds+\mathcal P_{0}(\gamma(h),t-h;\vartheta^{\gamma(h)},\gamma_{h}).
\end{align*}
Since the choice of $\vartheta^{\gamma(h)}$ satisfies \eqref{eq:stopping-control}, we have
\begin{align*}
  \mathcal P_{0}(\gamma(h),t-h;\vartheta^{\gamma(h)},\gamma_{h}) \geqslant W_{0}(\gamma(h),t-h)-\epsilon.
\end{align*}
Because in the present case $h\wedge\vartheta^{0}[\gamma]=h$ and $\mathbbm{1}_{\{\vartheta^{0}[\gamma]>h\}}=1$, the claim follows.

Taking the infimum over $\gamma\in\mathbb{A}^{\mathbf{x}}$ on both sides of \eqref{eq:sub-optimal-W0} and using \eqref{eq:sub-optimal-stopping-obstacle-hje}, we obtain
\begin{align*}
  \inf_{\gamma\in\mathbb{A}^{\mathbf{x}}} \mathcal P_{0}(\mathbf{x},t;\vartheta^{\epsilon},\gamma)
  \geqslant& \inf_{\gamma\in\mathbb{A}^{\mathbf{x}}} \left\{\int_{0}^{h\wedge\vartheta^{0}[\gamma]} L(\dot\gamma(s),\gamma(s),t-s)\,ds+\mathbbm{1}_{\{\vartheta^{0}[\gamma]>h\}} W_{0}(\gamma(h),t-h) \right\}-\epsilon\\
  =& \mathcal D_{h}^{\vartheta^{0}}W_{0}(\mathbf{x},t)-\epsilon\\
  \geqslant& \mathcal D_{h}W_{0}(\mathbf{x},t)-2\epsilon.
\end{align*}
Since $\vartheta^{\epsilon}\in\Theta$, it is admissible in the definition of $W_{0}(\mathbf{x},t)$. Hence
\begin{align*}
  W_{0}(\mathbf{x},t)\geqslant \mathcal D_{h}W_{0}(\mathbf{x},t)-2\epsilon.
\end{align*}
Letting $\epsilon\downarrow0$, we have \eqref{eq:inequality-2}. This completes the proof.
\end{proof}

\subsection{Viscosity characterization of the value function}
\label{app:bounded-data-viscosity}

\begin{proof}[Proof of Proposition \ref{prop:bounded-data-stopping-game}]
\textbf{Step 1.} We verify the viscosity subsolution property.
Let $\phi\in C^{1}$ touch $W_{0}$ from above at $(\mathbf{x}_{0},t_{0})\in\mathbb{R}^{N}\times\mathbb{R}^{+}$ and $W_{0}(\mathbf{x}_{0},t_{0})>0$.
Fix any $\mathbf{q}\in\mathbb{R}^{N}$ and set
\begin{align*}
  \gamma_\mathbf{q}(s)=\mathbf{x}_{0}-s\mathbf{q}, \quad s\in[0,\infty).
\end{align*}
For $h\in(0,t_{0})$, let $\epsilon_{h}=h^{2}$.
By Lemma \ref{lem:bounded-data-value-regularity} and Proposition \ref{prop:bounded-data-dpp}, choose $\vartheta_{h}\in\Theta$ such that
\begin{align*}
  W_{0}(\mathbf{x}_{0},t_{0}) \leqslant& \inf_{\gamma\in\mathbb{A}^{\mathbf{x}_{0}}} \left\{\int_{0}^{h\wedge\vartheta_{h}[\gamma]} L(\dot\gamma(s),\gamma(s),t_{0}-s)\,ds+\mathbbm{1}_{\{\vartheta_{h}[\gamma]>h\}} W_{0}(\gamma(h),t_{0}-h) \right\}+\epsilon_{h} .
\end{align*}
Testing the inner infimum by $\gamma_\mathbf{q}$, it holds that
\begin{align}
\label{eq:subsolution-dpp-test}
  W_{0}(\mathbf{x}_{0},t_{0}) \leqslant \int_{0}^{h\wedge\vartheta_{h}[\gamma_\mathbf{q}]} L(\dot\gamma_\mathbf{q}(s),\gamma_\mathbf{q}(s),t_{0}-s)\,ds+\mathbbm{1}_{\{\vartheta_{h}[\gamma_\mathbf{q}]>h\}} W_{0}(\gamma_\mathbf{q}(h),t_{0}-h)+\epsilon_{h} .
\end{align}
We claim that $\vartheta_{h}[\gamma_\mathbf{q}]>h$ for all sufficiently small $h$.
Otherwise, along a subsequence $\{h_{n}\}_{n}$ satisfying $h_{n}\downarrow0$ as $n\rightarrow\infty$ and $\vartheta_{h_{n}}[\gamma_\mathbf{q}]\leqslant h_{n}$, \eqref{eq:subsolution-dpp-test} and $R\geqslant0$ imply
\begin{align*}
  W_{0}(\mathbf{x}_{0},t_{0}) \leqslant \frac{1}{4}|\mathbf{q}|^{2}h_{n}+h_{n}^{2},
\end{align*}
which contradicts $W_{0}(\mathbf{x}_{0},t_{0})>0$.
Hence, for all sufficiently small $h$, \eqref{eq:subsolution-dpp-test} is reduced to
\begin{align*}
  W_{0}(\mathbf{x}_{0},t_{0}) \leqslant \int_{0}^{h} L(\dot\gamma_\mathbf{q}(s),\gamma_\mathbf{q}(s),t_{0}-s)\,ds+W_{0}(\gamma_\mathbf{q}(h),t_{0}-h)+h^{2} .
\end{align*}
Since $\phi$ touches $W_{0}$ from above, for all $h\ll1$,
\begin{align*}
  W_{0}(\gamma_\mathbf{q}(h),t_{0}-h)-W_{0}(\mathbf{x}_{0},t_{0}) \leqslant \phi(\gamma_\mathbf{q}(h),t_{0}-h)-\phi(\mathbf{x}_{0},t_{0}).
\end{align*}
Therefore
\begin{align*}
  0 \leqslant \int_{0}^{h} \left[ \frac{1}{4}|\mathbf{q}|^{2}-R(\gamma_\mathbf{q}(s),t_{0}-s)+\frac{d}{ds} \phi(\mathbf{x}_{0}-s\mathbf{q},t_{0}-s) \right]ds+h^{2} .
\end{align*}
Dividing by $h$ implies
\begin{equation}
\label{eq:W0-sub-test}
  \begin{aligned}
    0 \leqslant& \frac{1}{4}|\mathbf{q}|^{2}-\frac{1}{h}\int_{0}^{h} R(\gamma_\mathbf{q}(s),t_{0}-s)\,ds\\
    &- \frac{1}{h}\int_{0}^{h} \mathbf{q}\cdot D\phi(\mathbf{x}_{0}-s\mathbf{q},t_{0}-s)\,ds-\frac{1}{h}\int_{0}^{h} \phi_{t}(\mathbf{x}_{0}-s\mathbf{q},t_{0}-s)\,ds+h .
  \end{aligned}
\end{equation}
For $s\in[0,h]$, we have $|\gamma_\mathbf{q}(s)-\mathbf{x}_{0}|+|(t_{0}-s)-t_{0}|\leqslant (1+|\mathbf{q}|)h$. Hence, we have
\begin{align*}
  A_{h}:= \frac{1}{h}\int_{0}^{h} R(\gamma_\mathbf{q}(s),t_{0}-s)\,ds \geqslant \inf_{\substack{ |\mathbf{z}-\mathbf{x}_{0}|+|\tau-t_{0}|<\rho_{h} }} R(\mathbf{z},\tau),
\end{align*}
where $\rho_{h}=2(1+|\mathbf{q}|)h$.
Since $\rho_{h}\downarrow0$ as $h\downarrow0$, the definition of the l.s.c. envelope implies
\begin{align*}
  \liminf_{h\downarrow0} \frac{1}{h}\int_{0}^{h} R(\gamma_\mathbf{q}(s),t_{0}-s)\,ds \geqslant R_{*}(\mathbf{x}_{0},t_{0}).
\end{align*}
Taking the upper limit as $h\downarrow0$ in \eqref{eq:W0-sub-test}, using the continuity of $\phi_t$ and $D\phi$, and recalling that $\underset{h\downarrow0}\limsup\,(-A_{h}) =-\underset{h\downarrow0}\liminf\,(A_{h})\leqslant -R_{*}(\mathbf{x}_{0},t_{0})$, we obtain
\begin{align*}
  \phi_{t}(\mathbf{x}_{0},t_{0})+\mathbf{q}\cdot D\phi(\mathbf{x}_{0},t_{0})-\frac{1}{4}|\mathbf{q}|^{2}+R_{*}(\mathbf{x}_{0},t_{0}) \leqslant0.
\end{align*}
Taking the supremum over $\mathbf{q}\in\mathbb{R}^{N}$ yields
\begin{align*}
  \phi_{t}(\mathbf{x}_{0},t_{0})+|D\phi(\mathbf{x}_{0},t_{0})|^{2} +R_{*}(\mathbf{x}_{0},t_{0})\leqslant0.
\end{align*}
Hence $W_{0}$ is a viscosity subsolution.

\textbf{Step 2.} We verify the viscosity supersolution property.
Let $\phi\in C^{1}$ touch $W_{0}$ from below at $(\mathbf{x}_{0},t_{0})$.
Since $W_{0}\geqslant0$, it suffices to show
\begin{align*}
  \phi_{t}(\mathbf{x}_{0},t_{0})+|D\phi(\mathbf{x}_{0},t_{0})|^{2} +R^{*}(\mathbf{x}_{0},t_{0})\geqslant0.
\end{align*}
Assume by contradiction that, for some $\eta>0$,
\begin{align*}
  \phi_{t}(\mathbf{x}_{0},t_{0})+|D\phi(\mathbf{x}_{0},t_{0})|^{2} +R^{*}(\mathbf{x}_{0},t_{0}) \leqslant -4\eta .
\end{align*}
By the definition of $R^{*}$ and the continuity of $\phi_{t}$ and $D\phi$, there exists $\delta\in(0,t_{0})$ such that, whenever $(\mathbf{z},\tau)\in B_{\delta}(\mathbf{x}_{0})\times(t_{0}-\delta,t_{0}+\delta)$,
\begin{align}
\label{eq:super-local-negative}
  \phi_{t}(\mathbf{z},\tau)+|D\phi(\mathbf{z},\tau)|^{2}+R(\mathbf{z},\tau) \leqslant -2\eta .
\end{align}
We also take $\delta$ small enough such that
\begin{align}
\label{eq:lower-touch}
  W_{0}(\mathbf{z},\tau)\geqslant\phi(\mathbf{z},\tau) \quad \text{in }B_{\delta}(\mathbf{x}_{0})\times(t_{0}-\delta,t_{0}+\delta).
\end{align}
Fix $h\in(0,\delta)$.
Using the DPP \eqref{eq:bounded-data-dpp} with the stopping rule $\vartheta\equiv t_{0}$, we obtain
\begin{align}
\label{eq:super-dpp-constant-stop}
  \phi(\mathbf{x}_{0},t_{0})=W_{0}(\mathbf{x}_{0},t_{0})
  \geqslant \inf_{\gamma\in\mathbb{A}^{\mathbf{x}_{0}}} \left\{\int_{0}^{h} L(\dot\gamma(s),\gamma(s),t_{0}-s)\,ds+W_{0}(\gamma(h),t_{0}-h) \right\}.
\end{align}
We claim that the right-hand side is strictly larger than $\phi(\mathbf{x}_{0},t_{0})$ for $h$ small, which contradicts \eqref{eq:super-dpp-constant-stop}.

Indeed, let $\gamma\in\mathbb{A}^{\mathbf{x}_{0}}$.
If $\gamma$ leaves $B_{\delta}(\mathbf{x}_{0})$ at some $s_{h}\in(0,h]$, then $|\gamma(s_{h})-\gamma(0)|=|\gamma(s_{h})-\mathbf{x}_{0}|\geqslant \delta$.
By the Cauchy-Schwarz inequality, we have
\begin{align*}
  \int_{0}^{h}|\dot\gamma(s)|^{2}\,ds\geqslant \frac{\delta^{2}}{s_{h}}\geqslant \frac{\delta^{2}}{h}.
\end{align*}
Since $W_{0}\geqslant0$ and $R\leqslant r_{2}$, it follows that
\begin{align*}
  \int_{0}^{h} L(\dot\gamma(s),\gamma(s),t_{0}-s)\,ds+W_{0}(\gamma(h),t_{0}-h) \geqslant \frac{\delta^{2}}{4h}-r_{2}h .
\end{align*}
For $h$ sufficiently small, $\frac{\delta^{2}}{4h}-r_{2}h$ is larger than $\phi(\mathbf{x}_{0},t_{0})+\eta h$.
It remains to consider trajectories satisfying $\gamma(s)\in B_{\delta}(\mathbf{x}_{0})$ for $s\in(0,h]$.
Then by \eqref{eq:lower-touch}
\begin{align*}
  &\int_{0}^{h} L(\dot\gamma(s),\gamma(s),t_{0}-s)\,ds+W_{0}(\gamma(h),t_{0}-h)-\phi(\mathbf{x}_{0},t_{0})\\
  \geqslant& \int_{0}^{h} \left[ \frac{1}{4}|\dot\gamma(s)|^{2}-R(\gamma(s),t_{0}-s) \right]ds+\phi(\gamma(h),t_{0}-h)-\phi(\mathbf{x}_{0},t_{0})\\
  =& \int_{0}^{h} \left[ \frac{1}{4}|\dot\gamma(s)|^{2}+D\phi(\gamma(s),t_{0}-s)\cdot\dot\gamma(s)-\phi_{t}(\gamma(s),t_{0}-s)-R(\gamma(s),t_{0}-s) \right]ds .
\end{align*}
Using Fenchel's inequality and \eqref{eq:super-local-negative}, we obtain
\begin{align*}
  &\int_{0}^{h} L(\dot\gamma(s),\gamma(s),t_{0}-s)\,ds+W_{0}(\gamma(h),t_{0}-h)-\phi(\mathbf{x}_{0},t_{0})\\
  \geqslant& -\int_{0}^{h} \left[ \phi_{t}+|D\phi|^{2}+R \right](\gamma(s),t_{0}-s)\,ds \geqslant 2\eta h .
\end{align*}
Combining the two cases, for all admissible trajectories and all sufficiently small $h$, it holds that
\begin{align}
\label{eq:contrad-super}
  \int_{0}^{h} L(\dot\gamma(s),\gamma(s),t_{0}-s)\,ds+W_{0}(\gamma(h),t_{0}-h) \geqslant \phi(\mathbf{x}_{0},t_{0})+\eta h .
\end{align}
Taking the infimum over $\gamma$ on both sides of \eqref{eq:contrad-super} contradicts \eqref{eq:super-dpp-constant-stop}. Therefore
\begin{align*}
  \phi_{t}(\mathbf{x}_{0},t_{0})+|D\phi(\mathbf{x}_{0},t_{0})|^{2} +R^{*}(\mathbf{x}_{0},t_{0})\geqslant0.
\end{align*}
Thus $W_{0}$ is a viscosity supersolution of \eqref{eq:bounded-data-obstacle}.

\textbf{Step 3.} We verify the initial data and prove uniqueness.
Lemma \ref{lem:bounded-data-value-regularity} together with Steps 1--2 implies that $W_{0}$ is a bounded continuous viscosity solution of \eqref{eq:bounded-data-obstacle} in $\mathbb{R}^{N}\times\mathbb{R}^{+}$. We claim that
\begin{align}
\label{eq:bounded-data-initial-data}
  W_{0}(\mathbf{x},t)\to w_{0}(\mathbf{x}) \quad\text{as }t\downarrow0
\end{align}
locally uniformly in $\mathbf{x}$.

Let $K\subset\subset\mathbb{R}^{N}$ be compact and set $M_{0}=\|w_{0}\|_{L^{\infty}}$.
Fix any $\mathbf{x}\in K$ and $t\in\mathbb{R}^{+}$.
We first establish the lower bound.
Take the stopping time $\vartheta\equiv t$.
For all $\gamma\in\mathbb{A}^{\mathbf{x}}$, by the Cauchy-Schwarz inequality
\begin{align*}
  \int_{0}^{t}|\dot\gamma(s)|^{2}\,ds \geqslant \frac{|\mathbf{x}-\gamma(t)|^{2}}{t}.
\end{align*}
Since $R\leqslant r_{2}$, we obtain
\begin{align}
\label{eq:bounded-data-initial-lower}
  W_{0}(\mathbf{x},t) \geqslant \inf_{\mathbf{y}\in\mathbb{R}^{N}} \left\{\frac{|\mathbf{x}-\mathbf{y}|^{2}}{4t}+w_{0}(\mathbf{y}) \right\} -r_{2}t.
\end{align}
Fix $\delta\in(0,1)$, and let $\omega_{K}$ be a modulus of continuity of $w_{0}$ on the compact set
\begin{align*}
  K_{1}=\{\mathbf{z}\in\mathbb{R}^{N}\mid\textrm{dist}(\mathbf{z},K)\leqslant1\}.
\end{align*}
If $|\mathbf{x}-\mathbf{y}|<\delta$, then $\mathbf{y}\in K_{1}$ and $w_{0}(\mathbf{y}) \geqslant w_{0}(\mathbf{x})-\omega_{K}(\delta)$.
If $|\mathbf{x}-\mathbf{y}|\geqslant\delta$, then
\begin{align*}
  \frac{|\mathbf{x}-\mathbf{y}|^{2}}{4t}+w_{0}(\mathbf{y}) \geqslant \frac{\delta^{2}}{4t}-M_{0}.
\end{align*}
For $t\in(0,t_{\delta}]$ with $t_{\delta}=\frac{\delta^{2}}{4(2M_{0}+1)}$, we have $\frac{\delta^{2}}{4t}-M_{0}\geqslant w_{0}(\mathbf{x})-\omega_{K}(\delta)$. Therefore
\begin{align*}
  \inf_{\mathbf{y}\in\mathbb{R}^{N}} \left\{\frac{|\mathbf{x}-\mathbf{y}|^{2}}{4t}+w_{0}(\mathbf{y}) \right\} \geqslant w_{0}(\mathbf{x})-\omega_{K}(\delta)
\end{align*}
for every $\mathbf{x}\in K$ and $t\in(0,t_{\delta}]$.
Combining this with \eqref{eq:bounded-data-initial-lower}, we get
\begin{align*}
  \liminf_{t\downarrow0}\inf_{\mathbf{x}\in K} \left(W_{0}(\mathbf{x},t)-w_{0}(\mathbf{x})\right) \geqslant -\omega_{K}(\delta).
\end{align*}
Letting $\delta\downarrow0$ yields
\begin{align}
\label{eq:initial-condition-lower-bound}
  \liminf_{t\downarrow0}\inf_{\mathbf{x}\in K} \left(W_{0}(\mathbf{x},t)-w_{0}(\mathbf{x})\right) \geqslant0.
\end{align}
On the other hand, \eqref{eq:W0-upper-bound} implies
\begin{align}
\label{eq:initial-condition-upper-bound}
  \limsup_{t\downarrow0} \sup_{\mathbf{x}\in K} \left(W_{0}(\mathbf{x},t)-w_{0}(\mathbf{x}) \right) \leqslant0.
\end{align}
Combining \eqref{eq:initial-condition-lower-bound} and \eqref{eq:initial-condition-upper-bound}, \eqref{eq:bounded-data-initial-data} follows immediately by the arbitrariness of $K\subset\subset\mathbb{R}^{N}$.
Since $W_{0}$ is continuous in $\mathbb{R}^{N}\times\mathbb{R}^{+}$ and converges locally uniformly to $w_{0}\in C(\mathbb{R}^{N})$ as $t\downarrow0$, it extends continuously to $\mathbb{R}^{N}\times[0,\infty)$ by setting $W_{0}(\mathbf{x},0)=w_{0}(\mathbf{x})$.
The comparison principle for the variational inequality with bounded continuous initial data implies uniqueness.
\end{proof}

\end{document}